\documentclass{amsart}

\usepackage{amsmath,amsfonts,amssymb,amscd,verbatim,delarray}
\newcommand{\F}{{\mathbb F}}

\newcommand{\Q}{{\mathbb Q}}
\newcommand{\Z}{{\mathbb Z}}

 \DeclareMathOperator{\Spec}{Spec}

\begin{document}
\numberwithin{equation}{section}

\newtheorem{theorem}{Theorem}[section]
\newtheorem{lemma}[theorem]{Lemma}
\newtheorem{prop}[theorem]{Proposition}
\newtheorem{proposition}[theorem]{Proposition}
\newtheorem{corollary}[theorem]{Corollary}
\newtheorem{corol}[theorem]{Corollary}
\newtheorem{conj}[theorem]{Conjecture}

\theoremstyle{definition}
\newtheorem{defn}[theorem]{Definition}
\newtheorem{example}[theorem]{Example}
\newtheorem{examples}[theorem]{Examples}
\newtheorem{remarks}[theorem]{Remarks}
\newtheorem{remark}[theorem]{Remark}
\newtheorem{algorithm}[theorem]{Algorithm}
\newtheorem{question}[theorem]{Question}
\newtheorem{problem}[theorem]{Problem}
\newtheorem{subsec}[theorem]{}


\def\toeq{{\stackrel{\sim}{\longrightarrow}}}
\def\into{{\hookrightarrow}}


\def\alp{{\alpha}}  \def\bet{{\beta}} \def\gam{{\gamma}}
 \def\del{{\delta}}
\def\eps{{\varepsilon}}
\def\kap{{\kappa}}                   \def\Chi{\text{X}}
\def\lam{{\lambda}}
 \def\sig{{\sigma}}  \def\vphi{{\varphi}} \def\om{{\omega}}
\def\Gam{{\Gamma}}   \def\Del{{\Delta}}
\def\Sig{{\Sigma}}   \def\Om{{\Omega}}
\def\ups{{\upsilon}}
\def\defis{-}

\def\F{{\mathbb{F}}}
\def\Q{{\mathbb{Q}}}
\def\Ql{{\overline{\Q }_{\ell }}}
\def\CC{{\mathbb{C}}}
\def\R{{\mathbb R}}
\def\V{{\mathbf V}}
\def\D{{\mathbf D}}

\def\XX{\mathbf{X}^*}
\def\xx{\mathbf{X}_*}

\def\AA{\Bbb A}
\def\HH{\mathbb H}
\def\PP{\Bbb P}

\def\Gm{{{\mathbb G}_{\textrm{m}}}}
\def\Gmk{{{\mathbb G}_{\textrm m,k}}}
\def\GmL{{\mathbb G_{{\textrm m},L}}}
\def\Ga{{{\mathbb G}_a}}

\def\Fb{{\overline{\F }}}
\def\Kb{{\overline K}}
\def\Yb{{\overline Y}}
\def\Xb{{\overline X}}
\def\Tb{{\overline T}}
\def\Bb{{\overline B}}
\def\Gb{{\bar{G}}}
\def\Ub{{\overline U}}
\def\Vb{{\overline V}}
\def\Hb{{\bar{H}}}
\def\kb{{\bar{k}}}
\def\HOM{{\Hom(W,G)}}

\def\Th{{\hat T}}
\def\Bh{{\hat B}}
\def\Gh{{\hat G}}

\def\cF{{\mathfrak{F}}}
\def\cC{{\mathcal C}}
\def\cU{{\mathcal U}}
\def\cX{{\mathcal X}}

\def\Xt{{\widetilde X}}
\def\Gt{{\widetilde G}}

\def\gg{{\mathfrak g}}
\def\hh{{\mathfrak h}}
\def\Iks{{\mathfrak X}}
\def\MM{{\mathfrak M}}
\def\min{^{-1}}

\def\textrm#1{\text{\textnormal{#1}}}

\def\GL{\textrm{GL}}            \def\Stab{\textrm{Stab}}
\def\Gal{\textrm{Gal}}          \def\Aut{\textrm{Aut\,}}
\def\Lie{\textrm{Lie\,}}        \def\Ext{\textrm{Ext}}
\def\PSL{\textrm{PSL}}          \def\SL{\textrm{SL}}
\def\loc{\textrm{loc}}
\def\coker{\textrm{coker\,}}    \def\Hom{\textrm{Hom}}
\def\im{\textrm{im\,}}           \def\int{\textrm{int}}
\def\inv{\textrm{inv}}           \def\can{\textrm{can}}
\def\id{\textrm{id}}
\def\Cl{\textrm{Cl}}
\def\Sz{\textrm{Sz}}
\def\End{\textrm{End}}
\def\Ker{\textrm{Ker}}
\def\HA{\textrm{HA}}
\def\Sub{\textrm{Sub}}
\def\at{\textrm{at\ }}
\def\Var{\textrm{Var}}
\def\Rad{\textrm{Rad}}
\def\Ka{\textrm{K}}

\def\tors{_{\textrm{tors}}}      \def\tor{^{\textrm{tor}}}
\def\red{^{\textrm{red}}}         \def\nt{^{\textrm{ssu}}}

\def\sss{^{\textrm{ss}}}          \def\uu{^{\textrm{u}}}
\def\ad{^{\textrm{ad}}}           \def\mm{^{\textrm{m}}}
\def\tm{^\times}                  \def\mult{^{\textrm{mult}}}

\def\uss{^{\textrm{ssu}}}         \def\ssu{^{\textrm{ssu}}}
\def\comp{_{\textrm{c}}}
\def\ab{_{\textrm{ab}}}

\def\et{_{\textrm{\'et}}}
\def\nr{_{\textrm{nr}}}

\def\nil{_{\textrm{nil}}}
\def\sol{_{\textrm{sol}}}

\def\til{\;\widetilde{}\;}
\title{{\bf ALGEBRAIC LOGIC, VARIETIES OF ALGEBRAS AND ALGEBRAIC
VARIETIES }}
\bigskip
\author[ B.Plotkin]{{\bf B. Plotkin \\
\\ Institute of Mathematics, Hebrew University, Jerusalem } }
\address{ Institute of
Mathematics, Hebrew University, Jerusalem, ISRAEL}

\email{ plotkin@macs.biu.ac.il, borisov@math.huji.ac.il}

\date{\today}

\thispagestyle{empty} \vspace{1.0cm} \setcounter{page}{0}

\newpage


\begin{abstract}
The aim of the paper is discussion of connections between the
three kinds of objects named in the title. In a sense, it is a
survey of such connections; however, some new directions are also
considered. This relates, especially, to sections 3, 4 and 5,
where we consider a field that could be understood as an universal
algebraic geometry. This geometry is parallel to universal
algebra.

In the monograph \cite{BPl3} algebraic logic was used for
building up a model of
a database. Later on, the structures arising there turned out to be useful for
solving
several problems from algebra. This is the position which the present paper is
written from.
\end{abstract}

\maketitle

\begingroup
    \begin{center}
    \bf Contents
    \end{center}
\begin{list}{}{}
\item[~]  Introduction
\item[\S1.] Algebraic logic
    \begin{enumerate}
     \item[1.1.] Some problems
     \item[1.2.] $\Theta$-logic
     \item[1.3.] Algebraic $\Theta$-logic
    \end{enumerate}
\item[\S2.] Some applications
    \begin{enumerate}
     \item[2.1.] Closures of formula systems
     \item[2.2.] Quasigroups
     \item[2.3.] Algebraic logic in group representations
     \item[2.4.] Databases
    \end{enumerate}
\item[\S3.] Algebraic varieties and varieties of algebras
    \begin{enumerate}
     \item[3.1.] Basic concepts
     \item[3.2.] Hilbert's Nullstellensatz
     \item[3.3.] Verbal varieties
     \item[3.4.] Geometric equivalence of algebras
     \item[3.5.] Generalized equations
     \item[3.6.] Algebra and topology in connection with varieties
     \item[3.7.] Relation to the $\Theta$-structure of algebras
     \item[3.8.] Additional remarks
    \end{enumerate}
\item[\S4.] Varying the set of variables, the base variety and the base
algebra
    \begin{enumerate}
      \item[4.1.] Changing $X$
      \item[4.2.] Changing $\Theta$
      \item[4.3.] Changing the algebra $G$
      \item[4.4.] The category of varieties
    \end{enumerate}
\item[\S5.] Algebraic logic and algebraic varieties
    \begin{enumerate}
      \item[5.1.] Basic concepts
      \item[5.2.] Generalized varieties, topology and other matters
      \item[5.3.] Passing to submodels
      \item[5.4.] Geometry on the level of quantifier-free logic
      \item[5.5.] Halmos algebras and Boolean algebras of varieties.
          Conclusion
    \end{enumerate}
\item[~] Bibliography
\end{list}
\endgroup

\section*{Introduction}

Essentially, the paper consists of two parts. The first part includes \S1
and \S2. In the first section, preliminary information on algebraic logic is
given (see also \cite{BPl3}),
while the second one contains a survey of certain applications.

\def\authorslast{B.I. Plotkin}

All the rest makes up the second part. Its main subject is equations and
identities in arbitrary algebraic structures. The very notion of an equation is
treated very widely, and a solution of an equation is a point regarded as
an algebra homomorphism. Here the connections existing between algebraic
logic and universal algebra work. The investigation is carried out after the
pattern of algebraic geometry, and geometry is considered on three levels. These
are the equational logic level(\S3 and \S4), and the levels of
quantifier-free logic and first-order logic (\S5). Everywhere, what we have in
mind is,in fact, $\Theta$-logic, where $\Theta$ is some variety of algebras.

On the equational level, a certain general statement of the
Hilbert Nullstellensatz (i.e. theorem on zeros) is given; it is
applicable in all cases and admits one more view on the classical
theorem. This general formulation has also a linkage with the
notion of geometric equivalence of algebras. For every algebra $G
\in \Theta$ a category $K_G$ of algebraic varieties related with
$G$ is defined. It is contravariantly embedded in the category of
algebras from $\Theta$. If algebras $G_1$ and $G_2$ are
geometrically equivalent, then the categories $K_{G_!}$ and
$K_{G_2}$ are also equivalent.

What can be said in general about the algebras $G_1$ and $G_2$ for which
the categories $K_{G_1}$ and $K_{G_2}$ are equivalent? The solution of
this problem we consider in the separate paper.

It seems to us that the material contained in the second part is only the
beginning of a vast theme.
The part is two-aimed. On the one hand, it is the wish for seeing how
fare the ideas
of algebraic geometry are applicable in universal algebra, and on the other
one--the needs and interests of the algebra itself, desire to look at its
problems at a new "viewing angle".

I would like to express my gratitude to the Algebra Department of
the Institute of mathematics at the Hebrew University in
Jerusalem, where I have excellent  conditions for my work.

The moral support by Professors A.~Lubotzky, A.~Mann, I.~Rips, A.~Shalev and
A.~Reznikov was of great importance,
and I shall ever remember the support by S.~Amitsur. I am
especially grateful to Ilya Rips; regular discussions with him are always
exceptionally interesting. A very essential help was rendered me by Prof.
J.~C\={\i}rulis from University of Latvia in Riga. I am thankful to him both for
his help
in my work and as to a representative of the Riga algebraic seminar
which I have
been linked with for years.

\section*{\S1. Algebraic logic}
\setcounter{paragraph}{0}

\paragraph{1. Some problems.}
\stepcounter{paragraph}
Let $\Theta$ be a variety of algebras, $W$ a free algebra of
countable rank in $\Theta$, $u$ a formula of the first-order
calculus specialized in $\Theta$. It will be clear from the
following four examples what we mean by a "specialized" formula.
\begin{enumerate}
    \item[1.] $ w \equiv {w'} $,
    \item[2.] $ w_1 \equiv w_1'\vee \cdots \vee w_n \equiv w_n' $,
    \item[3.] $ w_1 \equiv w_1'\wedge \cdots \wedge w_n \equiv w_n'
        \rightarrow w \equiv {w'} $,
    \item[4.] $ w_1 \equiv w_1'\vee \cdots \vee w_n
        \equiv w_n' \vee
         v_1 \not\equiv v_1'\vee \cdots \vee v_m \not\equiv v_m' $.
\end{enumerate}
All $w$ and $v$ are  elements of $W$, and if we speak of many-sorted algebras
(which also are possible), then
the sign "$\equiv$" connects  elements of the same sort. $\Theta$
can be the variety of all groups, semigroups, quasigroups,
rings, associative or Lie rings, automata, etc.

Let further $G$ be an algebra from $\Theta$. On $G$, also some relations could
be given. How should one  understand that a formula $u$ of the
corresponding language is satisfied in $G$, or more generally,
what is to be meant by the value of  $u$ in the algebra
or model $G$?

If $u$ is
$ w \equiv {w'}$, the value of $u$ in $G$ is the set of those
homomorphisms $ \mu \colon W \to G$ for which the equality
$ w^\mu = {{w'}}^\mu$ holds in $G$. If it is the whole set $\Hom(W,G)$, then $u$
is said to be valid in $G$, or an identity of the algebra $G$.

We have a similar situation also in the next example. The value of
$u$ in $G$ is the subset of $\Hom(W,G)$ consisting of those
$\mu \colon W \to G$ for which
at last one of the equalities
$ w_1^\mu = {w_1'}^\mu$, \ldots, $w_n^\mu = {w_n'}^\mu$ holds
in $G$. If it is so for all $\mu$, then $u$ is a pseudoidentity of $G$.

The third example is connected with quasi-identities, and the
fourth one--with universal formulas.

For any $u$, the value of $u$ in $G$ is defined inductively, and it is a
subset of $\Hom(W,G)$ both for algebras and models.

We have defined the value of a formula semantically, and the problem now is to
formalize this semantics using any natural
structure. So, it is convenient to turn to algebraic logic.

Let us consider the following class of problems. Let $T$ be some
set of axioms (formulas in a given specialized language), and let $K=T'$ be the
class of models it specifies. All formulas from $T$
are valid in these models. The class $K$ is axiomatic,  and if $K'=T''$
is the set of formulas valid in algebras of $K$, then $T''$
is the closure of $T$. How could we get all formulas from $T''$
syntactically  if we proceed from $T$? This problem is well solved in
logic, but its solution in algebraic logic is more natural.

It is of some interest to speak of formulas of a specified kind, for example,
identities, quasi-identities, pseudoidentities, universal
formulas, etc., and to consider the respective closures. Here we
deal with algebras without additional relations. The problem is
easily solved for identities by using the free algebra $W$, but in other cases
it is natural to turn our attention to algebraic logic.

We can also speak of the closure of a class. A set of identities specifies a
variety. The closure of an arbitrary class of algebras from $\Theta$ up to
the variety (a subvariety of $\Theta$) is constructed in accordance with the
Birkhoff theorem. There are similar theorems for pseudovarieties,
which are defined by pseudoidentities, for quasivarieties and  universal
classes.

The system of formulas valid in $G$ is an important characteristic of an
algebra $G\in\Theta$ . It is especially true of identities.
Varieties are also of interest because each variety has free algebras,
and each variety isolates the corresponding
verbal congruence on every algebra $G$.

All subvarieties of $\Theta$ are controlled by the free algebra $W$. For
other axiomatic classes, an algebra $U = U(\Theta)$ constructed from the
corresponding specialized first-order calculus using  $W$
plays the role of such controlling object. Of course, varieties
could also be controlled by this algebra.

\paragraph{2. $\bf\Theta$-logic.}
\stepcounter{paragraph}
$\Theta$-logic is built up on the ground of some
variety of algebras $\Theta$. Since we have in mind many-sorted algebras as
well, let us recall some related concepts (see \cite{Mal1,Hig,BL,BW}).

First of all, we fix a set of sorts $\Gamma$. Correspondingly, we consider a
many-sorted set $G = (G_i,\ i \in \Gamma)$.  Each $G_i$ is the domain of the
sort $i$.  Denote by $\Omega$ a set of operation symbols.  Each
$\omega \in \Omega$ has a  definite type
    $\tau = \tau (\omega) = (i_1, \ldots, i_n; j)$,
where $i, j \in \Gamma$.  Such a symbol $\omega$ is realized in $G$ as
an {\it operation\/}, i.e. a map
$$
\omega\colon G_{i_1} \times \cdots \times G_{i_n} \to G_j.
$$
In each algebra $G$, all symbols $\omega \in \Omega$ are assumed to be
realized this way. So, $G$ is an $\Omega$--algebra.

For given $\Gamma$ and $\Omega$, an {\it algebra morphism\/} $G \to G'$ has the
form
$$
\mu = ( \mu_i,\ i \in \Gamma) \colon
    G = (G_i,\ i \in \Gamma) \to G' = (G'_i,\ i \in \Gamma),
$$
where each $\mu_i \colon G_i \to G'_i$ is a map of sets. {\it Algebra
 homomorphisms\/} are morphisms that preserve operations.  For
    $\mu = (\mu_i,\ i \in \Gamma)$  and  $\omega \in \Omega$
this means that if
    $\tau (\omega) = (i_1, \ldots, i_n; j)$
and
    $a_1 \in G_{i_1}, \ldots, a_n \in G_{i_n}$,
then
$$
(a_1 \cdots a_n \omega)^{\mu_j} = a^{\mu_{i_1}}_1 \cdots a^{\mu_{i_n}}_n \omega.
$$
The multiplication of such many-sorted maps is defined componentwise, and the
product of homomorphisms is a homomorphism. Sometimes we write $a^\mu$ instead
of $a^{\mu_i}$.  If all $\mu_i$ are surjections (injections), then $\mu$ is
also a surjection (injection).
If all $\mu_i$ are bijections, then $\mu$ is a bijection, and
    $\mu^{-1} = (\mu^{-1}_i,\ i \in \Gamma)$.
A bijective homomorphism $\mu$ is an isomorphism.

A {\it kernel\/} of a homomorphism $\mu\colon G \to G'$ has the form
$\rho = (\rho_i,\ i \in \Gamma)$, where each $\rho_i$ is the kernel
equivalence of the map $\mu_i\colon  G_i \to G'_i$.  A {\it many-sorted
equivalence\/} $\rho$ is a {\it congruence\/} if it preserves all
operations $\omega \in \Omega$.  This means that if
    $\tau (\omega) = (i_1, \ldots, i_n; j)$      and if
    $a_1, {a'}_1 \in G_{i_1}$, \ldots, $a_n, {a'}_n \in G_{i_n}$,
then
$$
a_1 \rho_{i_1} {a'}_1, \ldots, a_n \rho_{i_n} {a'}_n\, \Rightarrow\,
    (a_1 \ldots a_n \omega) \rho_j ({a'}_1 \ldots {a'}_n \omega).
$$
If $\rho$ is a congruence of $G$, then one can consider the quotient algebra
    $G / \rho = (G_i / \rho_i ,\ i \in \Gamma)$.
The notions of a subalgebra and a Cartesian product of algebras are defined in
a natural way.  If, for example,
    $G^\alpha = (G^\alpha_i ,\ i \in \Gamma)$
are $\Omega$-algebras, $\alpha \in I$, then
    $\prod G^\alpha = (\prod\limits_\alpha  G^\alpha_i,\ i \in\Gamma)$,
and if
    $a_1 \in \prod\limits_\alpha G^\alpha_{i_1}$, \ldots,
    $a_n \in \prod\limits_\alpha G^\alpha_{i_n}$
then
    $(a_1 \cdots a_n) \omega (\alpha) = a_1(\alpha) \cdots
        a_n (\alpha) \omega$,
provided $\tau (\omega) = (i_1, \ldots, i_n; j)$.

Now let us make some notes on varieties of algebras.

 Varieties of algebras are specified by identities.
Let us consider the many-sorted case in detail.
Let $X = (X_i,\ i \in \Gamma)$ be a
many-sorted set with all sets $X_i$ countable.  Starting from the set of
operation symbols $\Omega $, we can construct the
terms over $X$. Denote the system of terms by $W = (W_i,\ i \in\Gamma)$.
The inductive definition of $W$ runs as follows.
  Every set $X_i$ is included in $W_i$, the set of terms of the sort $i$.
  If  $\omega \in \Omega$, $\tau (\omega)
  = (i_1, \ldots, i_n;j)$ and $w_1, \ldots, w_n$ are terms of the sorts $i_1,
  \ldots, i_n$ respectively, then $w_1 \cdots w_n \omega$ is a term of the
sort $j$.
  If there are nullary operation symbols (with $n = 0$) of the sort $i$ in
  $\Omega$, then they also belong to $W_i$.  An algebra $W$ is, naturally,
an $\Omega$-algebra, and it is called the
{\it absolutely free $\Omega$-algebra}.

An {\it identity\/} is a  formula of the kind $w \equiv {w'}$,
where $w$ and ${w'}$ are terms from $W$ of the same sort, say $i$.  Such
a formula is {\it valid in an $\Omega$-algebra}
$G = (G_i,\ i \in \Gamma)$, if $w^{\mu_i} = {{w'}}^{\mu_i}$
for every homomorphism $\mu\colon  W \to G$.

A set of identities determines a {\it variety of $\Omega$-algebras}, i.e. the
class of algebras satisfying all identities from a given set.
In every variety $\Theta$ the set $X = (X_i,\ i \in \Gamma)$
picks out a {\it free algebra}, which is a quotient algebra of the absolutely free
algebra $W$.

The Birkhoff's theorem holds also in the many-sorted case.

{\it A class $\Theta$ is a variety if and only if it is closed under Cartesian
products, subalgebras and homomorphic images.}

Any variety $\Theta$ can be taken for the  initial variety, and one can
classify various subvarieties and other axiomatizable classes in
$\Theta$.  If $X = (X_i,\ i \in \Gamma)$ is a many-sorted set, then
the free algebra in $\Theta$ associated with  $X$ is also denoted by
$W = (W_i,\ i \in \Gamma)$.  Subvarieties of $\Theta$ are defined
by  identities $w \equiv {w'}$, where $w$ and ${w'}$ are terms of the
same sort in this new $W$.  Fully characteristic
congruences of $W$ correspond to closed sets of identities.

For an arbitrary algebra $G = (G_i,\ i \in \Gamma)$ from the variety $\Theta$,
we can consider the set of homomorphisms Hom$(W,G)$.

Now we pass to $\Theta$-logic.

We again consider the many-sorted case with the set of sorts $\Gamma$.  Let $X$
be the set of variables,
with the stratification map $n: X \to \Gamma$. This map is surjective and
divides $X$ into sets $X_i$, $ i \in \Gamma$.  Each $X_i$ consists of the
variables of the sort $i$, where the sort of $x$ is $i = n (x)$.
So, we have a many-sorted set $X = (X_i,\ i \in \Gamma)$.

We fix the set of operation symbols $\Omega$ (so a certain $\Gamma$-type
$\tau$ corresponds to each $\omega \in \Omega$), and select the variety
$\Theta$ of $\Omega$-algebras.  With $\Theta$ we associate a logic which is
called {\it $\Theta$-logic}.  Let $W = (W_i,\ i \in \Gamma)$ be the free
over $X$ algebra in $\Theta$.  We also fix the set of relation symbols  $\Phi$;
each $\varphi \in \Phi$ has a
type $\tau = \tau (\varphi) = (i_1, \ldots, i_n)$,
realized in the algebra $G = (G_i,\ i \in \Gamma)$
as a subset of the Cartesian product $G_{i_1} \times \cdots \times G_{i_n}$.

Now we can construct the set of formulas of $\Theta$-logic.
First, we define the elementary formulas. These are of  the form
$$ \varphi (w_1, \ldots, w_n), $$
where $\tau (\varphi) = (i_1, \ldots, i_n)$,
$ w_1 \in W_{i_1}, \ldots, w_n \in W_{i_n}$.
Denote the set of all elementary formulas by $\Phi W$.  We let
$L = \{ \vee, \wedge, \neg, \exists x\}_{ x \in X }$ be the
signature of logical symbols, and we construct the absolutely
free algebra over $\Phi W$ in this signature.  The corresponding formula
algebra  is denoted by $L \Phi W$.

The set  $L \Phi W$ has a part  whose elements  are logical axioms.

They are the usual axioms of  calculus    with functional symbols.
(Warning: the axiom set can be not effective in general.)
The terms of this calculus are $\Theta$-terms.  The rules of inference are
standard:

\begin{enumerate}
    \item[1.] {\it Modus ponens}: $u$ and $u \to v$ imply $v$,
    \item[2.] {\it Generalization\/}: $u$ implies $\forall x u$.
\end{enumerate}

\noindent
Here $u, v \in L \Phi W$, $ u \to v$ stands for $\neg u \vee v$, and
$\forall x u$ means $\neg \exists x \neg u$.

{\it Formulas together with axioms and rules of inference constitute
the (first-order) $\Theta$-logic}.  In particular, we can speak
about the logic of group theory, the logic of ring theory, etc.

A set of formulas is said to be {\it closed\/} if it contains all the axioms
and is invariant with respect to the rules of inference.

With each $\varphi \in \Phi$ of type $\tau = (i_1, \ldots, i_n)$ we
associate a set of mutually distinct variables
$x^\varphi_{i_1}, \ldots, x^\varphi_{i_n}$.  Then the formula
$\varphi (x^\varphi_{i_1}, \ldots, x^\varphi_{i_n})$ is called a basic
one.  The variables occurring in basic formulas are called {\it attributes}.
The set of attributes is denoted by $X_0$, the set of basic formulas is denoted by
$\Phi X_0$.  It is a small part of the set of elementary formulas $\Phi W$.

For applications, logic should also contain equalities.
An {\it equality} is a formula of the
kind $w \equiv {w'}$, where $w$ and ${w'}$  are elements of $W$ of the same sort.
Informally, such a formula can be regarded either as an equation or as an
identity. In the {\it equality logic\/} the set of logical axioms is extended
by specific equality axioms.
Any equality is considered to be an additional elementary formula.
The absolutely free formula algebra can  be constructed in equality logic, too.
We use the same notation $L \Phi W$ for it.

\paragraph{3. Algebraic $\bf\Theta$ logic.}
\stepcounter{paragraph}
Normally, quantifiers are logical symbols, but they also can be defined
to be operations of a Boolean algebra:

{\it An existential quantifier $\exists$ on a Boolean algebra $H$ is a map
$\exists\colon H \to H$ subject to the  following conditions}:
    \begin{enumerate}
\item[1.] $\exists 0 = 0$,
\item[2.] $a \le \exists a$,
\item[3.] $\exists (a \wedge\exists b) = \exists a \wedge \exists b$,
    \end{enumerate}
where $a, b \in H$ and  $0$ is the zero element of $H$.

Every  quantifier is a closure operator on $H$, and two existential quantifiers
can be non-permutable. Furthermore,
an existential quantifier is additive:
    $\exists (a \vee b) = \exists a \vee \exists b$.
The set of all $a \in H$ with $\exists a = a$ is a subalgebra of the Boolean
algebra $H$.
{\it A universal quantifier} $\forall \colon H \to H$ is defined dually:
    \begin{enumerate}
\item[1.] $\forall 1 =1$,
\item[2.] $\forall  a \le a$,
\item[3.] $\forall (a \vee \forall  b) = \forall a \vee \forall b$.
    \end{enumerate}
Analogously we have: $ \forall (a \wedge b) = \forall a \wedge \forall b$.

There is a well-known correspondence between the two species of
quantifiers which allows to switch back and forth from $\exists$ to $\forall$.

In algebraic logic, algebraic structures of logic are constructed and studied.
For example, with the classical propositional
calculus associated are Boolean algebras, and with the intuitionistic
propositional calculus Heyting algebras are connected.

There are three approaches to algebraization of first-order logic,
namely, the Tarski's cylindric algebras \cite{HMT}, the Halmos' polyadic
algebras \cite{Hal} and the
categorical approach by Lawvere \cite{Law1,Law2}. These approaches are based on
deep analysis of calculus (as a rule, we shall use the word ``calculus'' for
the first-order $\Theta$-logic).
There are also algebraic equivalents of nonclassical first-order logics
\cite{Ge}.  Similar constructions are developed for other logics
\cite{Dis1}.
See also \cite{Ci0, Ci1, Ci2, Ci3, Ci4, Ci5, Vo0, Vo1, Vo2, Ro', Ro}.

We consider the respective algebraizations of $\Theta$-logic.
This generalization is necessary for databases with the {\it data type\/}
$\Theta$,
and for algebra itself as well.  We confine ourselves with Halmos algebras.

Let us proceed from a  fixed {\it scheme} consisting of a set
$X=(X_i,\ i \in \Gamma)$, a variety $\Theta$, and an  algebra $W$ free in
$\Theta$ over $X$.  Since the latter is uniquely determined by $\Theta$ and $X$,
we do not write it out. We also take into account the semigroup $\End W$, whose
elements are considered to be additional operators.

 $\bf 3.1.$ Definition.
Suppose $(X, \Theta)$ is a scheme. An algebra  $H$ is a  Halmos algebra in
this scheme if
    \begin{enumerate} \def\ite[#1]{\item[{\rm #1}]}
 \ite[1.1.]  $H$ is a Boolean algebra.
 \ite[1.2.]  The semigroup $\End W$ acts on $H$ as a semigroup of Boolean
    endomorphisms.
 \ite[1.3.] Action of quantifiers of the form $\exists (Y)$, $Y \subset X$,
    is defined.
    \end{enumerate}
These actions are connected by the following conditions:
    \begin{enumerate} \def\ite[#1]{\item[{\rm #1}]}
 \ite[2.1.]  $\exists (\emptyset)$ acts trivially.
 \ite[2.2.]  $\exists (Y_1 \cup Y_2) = \exists (Y_1) \exists (Y_2)$.
 \ite[2.3.]  $s_1 \exists (Y) =s_2 \exists (Y)$ if $s_1, s_2 \in
    \End W$ and $s_1(x) = s_2(x)$ for $x \in X \setminus Y$.
 \ite[2.4.]  $\exists (Y) s = s \exists (s^{-1} Y)$ for $s \in\End W$ if the
    following conditions are fulfilled:
    \ite[ \ 1)] $s (x_1) = s (x_2) \in Y$ implies $x_1 = x_2$,
    \ite[ \ 2)] If $x \notin s^{-1} Y$, then
        $\Delta s (x) \cap Y = \emptyset$.
    \end{enumerate}

Here, $s^{-1} Y = \{ x,\ s (x) \in Y \}$ and $\Delta s (x)$ is the {\it
support\/} of the element $w = s (x) \in W$,
i.e., {\it the set of all $x \in X$ involved in the (representation of the)
element $w$}.
For a precise definition of a support for Halmos algebras, see \cite{BPl3}.

All Halmos algebras in the given scheme form a variety denoted by
$\HA_{\Theta}$. We will deal with such Halmos algebras and
occasionally will call them {\it Halmos algebras specialized in} $\Theta$.
In the next three subsections we shall present examples of Halmos algebras.

 $\bf 3.2$. \quad\rm Given $G \in \Theta$, consider $\Hom(W, G)$.  Denote
by $M_G$ the set of all subsets of $\Hom(W, G)$, i.e.
$M_G = \Sub(\Hom(W,G))$. In fact, $M_G$ is a Boolean algebra.  If $A \in M_G$,
$ \mu \in \Hom(W, G)$ and $s \in \End W$,
then define $\mu s$ by the rule $\mu s (x) = \mu (s (x))$.  An action
of the semigroup $\End W$ on $M_G$ is defined by
    $$ \mu \in s A \Leftrightarrow \mu s \in A. $$
Where $Y \subset X$, we let $\mu \in \exists (Y) A$ if
there is $\nu\colon W \to G$ in  $A$ such that $\mu (x) = \nu (x)$ for every
$x \in X \setminus Y$. This defines the action of a quantifier on $M_G$.
All the axioms of Halmos algebra are
fulfilled in $M_G$, and so this is the first example of a $\HA_\Theta$-algebra.

Let $H$ be a Halmos algebra and $h \in H$. We denote by $\Delta h$ the {\it
support of $h$\/}:
    $$ \Delta h = \{ x \in X,\ \exists x h \neq h \}. $$
If $\Delta h$ is finite, then the element $h$ is said to be {\it finitely
supported}. In the previous item we needed supports of elements of
the algebra $W$, while here we deal with supports in a Halmos
algebra. All finitely supported elements of $H$ constitute a
subalgebra called the {\it locally finite part of $H$}. The
algebra $H$ is locally finite if all of its elements have  finite supports.

 $\bf 3.3$. \quad\rm Denote by $V_D$ the locally finite part of the Halmos algebra
$M_G$. Here $A \in V_G$ if, for some finite subset $Y \subset X$ and elements
$\mu, \nu \in \Hom(W, D)$,
    $$ \mu \in A \Leftrightarrow \nu \in A$$
whenever $\mu (x) = \nu (x)$ for all $x \in Y$.
In other words, belonging of a row to the set $A$ is checked on a finite part of
$X$.

 $\bf 3.4$. \quad\rm  Now let us consider our main example: the
$\HA_{\Theta}$-algebra of first order calculus. Suppose that all $X_i$ in
$X = (X_i,\ i \in \Gamma)$ are infinite.  Besides that, the set $\Phi$ of
relation symbols is  added to the scheme, and we again have the formula algebra
$L \Phi W$.

Variables occur in formulas, and an occurrence of a variable is either free or
bound.  We define the action of an element $s \in \End W$ on the set $L \Phi W$
as follows.
If $u$ is a formula  and the variables $x_1, \ldots, x_n$ occur freely in it,
then we substitute them by $s x_1, \ldots, s x_n$, respectively, in all free
occurrences of them.  So we get $s u$.
For example, it follows from this definition that if
    $u = \varphi(w_1, \ldots, w_n)$ is
an elementary formula, then $s u = \varphi(s w_1, \ldots, s w_n)$.
However, the definition does not provide a representation of the semigroup
End$W$ in the set $L \Phi W$.  Simple examples show that the condition
$(s_1 \cdot s_2) u = s_1 (s_2 u)$ is not fulfilled.
Let us define an equivalence $\rho$ on the set $L \Phi W$ by the rule:
$u \rho v$ if $u$ and $v$  only differ in the names of bound variables. Take
the quotient set $L \Phi W/ \rho = \overline{L \Phi W}$, and call this passage
 {\it factorization by renaming bound variables}. We  call elements of
 $\overline{L \Phi W}$ formulas, too, but they are regarded up to
renaming bound variables.  It is easy to see that the equivalence $\rho$ is
compatible with the signature, but is not--with the action of elements
from End$W$.  Therefore, all operations from the set
    $L = \{ \vee, \wedge, \neg, \exists x, \;x \in X \}$
are defined on $\overline{L \Phi X}$, but  action of elements from
End$W$ has to be defined separately.  This is carried out as follows.

Let $u$ be a formula, $x_1, \ldots, x_n$ be the list of all of its free
variables, and take $s x_1, \ldots, s x_n$, $s \in \End W$. We say that
$u$ is {\it open\/} for $s$ if there are no bound variables in $u$
 belonging to any of the sets $\Delta s x_1, \ldots, \Delta s x_n$.
 For each $u$, we denote by $\bar u$
the corresponding class of equivalent elements.  For $s \in \End W$ we
always can
find some formula $u'$ in the class $\bar u$ which is open for $s$.
Then we set $s \bar u = \overline{s u'}$.
It is easily understood that if we have another formula $u''$ in
$\bar u$ which is open for $s$, then $s u ' \rho  s u''$ and
$\overline {su'} =\overline {s u''}$.  Hence the definition of $s\bar u$
is correct.  This rule gives the
representation of the semigroup End$W$ as a semigroup of transformations of the
formula set $\overline{L \Phi W}$.

From now on we proceed from the set of formulas $\overline{L \Phi W}$.
Axioms and  rules of inference are related to this set, too.
As before, they are standard.

Now we pass to {\it Lindenbaum-Tarski algebra}.  We have
an equivalence $\tau$ which is defined as follows:  $\bar u \tau \bar v$ if
the formula $(\bar u \to \bar v) \wedge (\bar v \to \bar u)$ is derivable.
It can be verified that $\tau$ is congruence on
$\overline{L \Phi W}$, and this $\tau$ is also compatible with the action of
the semigroup End$W$.

Denote by $U$ the result of factorization of $\overline{L \Phi W}$
by $\tau$. Define an equivalence  $\eta$  on $L \Phi W$  by the rule:
    $u \eta v \Leftrightarrow \bar u \tau \bar v$.
Since $\rho \subset \tau$,
the set $U=\overline{L \Phi W}/\tau$ can be identified with $L \Phi W/ \eta$.
It is important to emphasize that the equivalence $\eta$ can also be defined
by means of Lindenbaum-Tarski scheme,  and $U$ is the {\it Lindenbaum-Tarski
algebra}. It can be proved that:

\begin{enumerate}
\item[1.] {\it $U$ is a Boolean algebra with respect to the operations
$\vee, \wedge, \neg$};
\item[2.] {\it The semigroup \End$W$ acts on $U$ as a semigroup of endomorphisms
of this algebra};
\item[3.] {\it The operations  $\exists x$ are pairwise permutable.  This allows
us to define in $U$ quantifiers $\exists (Y)$ for all $Y \subset
X$.}
\end{enumerate}

All the above leads to the following result:

 $\bf 3.5.$ Theorem.
The algebra $U$ with the indicated  operations is an
algebra in HA$_\Theta$.

This is a syntactical approach to the definition of Halmos algebra
of first-order $\Theta$--logic.  Such an approach is realized by Z. Diskin
\cite{Dis1}.  There is also a semantical approach which is described in
\cite{BPl3}.
Both of them give  the same result.  Finally, we can use
the verbal congruence of the variety HA$_\Theta$, and obtain once more the
same algebra $
U$. Elements of $U$ are the formulas of $\Theta$-logic, now considered up to
the equivalence just defined.

Now we are going to discuss homomorphisms of Halmos algebras.
We start with a very important property of the algebra $U$, and note first
of all that this algebra is locally finite.
Take the basic set $\Phi X_0$ in the formula algebra $L \Phi W$, and
let $U_0$ be the corresponding basic set in $U$. The set $U_0$ generates the
algebra $U$.

 $\bf 3.6.$ Theorem.
Let $H$ be an arbitrary Halmos algebra, and let
    $\zeta\colon \Phi X_0 \to H$
be a map such that, for every
$u = \varphi(x_1, \ldots, x_n) \in \Phi X_0$,
$\Delta \zeta (u) \subset \{ x_1, \ldots, x_n \}$.  Such $\zeta$ gives
another map $\zeta\colon U_0 \to H$, and the latter one is uniquely extended
to a homomorphism $\zeta\colon U \to H$.

Given a model $(G, \Phi, f)$, $ G \in \Theta$, consider the particular case
when $H = V_G$. Define
    $\zeta = \hat f \colon \Phi X_0 \to V_G$
by the rule: if $u = \varphi(x_1, \ldots, x_n) \in \Phi X_0$, then
$$
 \hat f (u) = \{ \mu, \ \mu \in \Hom (W, G),\
            (\mu (x_1), \ldots, \mu (x_n)) \in f (\varphi) \}.
$$
Then $\Delta \hat f (u) = \{ x_1, \ldots, x_n \}$, and we have a homomorphism
    $$ \hat f : U \to V_G.$$
It follows from its definition that if $u = \varphi(w_1, \ldots, w_n )$ is an
elementary formula, then
$$
\hat f (u) = \{ \mu,\ \ \mu \in \Hom (W, D),
            (\mu (w_1), \ldots, \mu (w_n)) \in f (\varphi) \}.
$$
So, for every model $(G, \Phi, f)$,  $G \in \Theta$, we have the  canonical
homomorphism ${\hat f}\colon U \to V_G$, and every homomorphism
$U \to V_G$ proves to be of this sort.
Now we can say that, for every $u \in U$, $\hat f (u)$ is the value of
 $u$ in the model $G$. $f (u) = 1$ means that $u$ is valid in $G$.

We add a few  words about kernels of homomorphisms in the variety
HA$_\Theta$.  If $\sigma \colon  H \to H'$ is a homomorphism in HA$_\Theta$,
then we have two kernels: the coimage of the zero and the coimage of the unit.
The coimage of the zero is an ideal, and the  coimage of the  unit is a filter.
A subset $T$ of $H$ is a {\it filter} if
    \begin{enumerate}
\item[1.] {\it $a \wedge b \in T$ if $a$ and $b$ belong to $T$},
\item[2.] {\it $a \vee b \in T$ if $a \in T$ and $b \in H$ },
\item[3.] {\it $\forall (Y) a \in T$ if $a \in T$, $ Y \subset X$.}
    \end{enumerate}

The definition of an ideal is dual.  For every filter $T$, the quotient algebra
$H/ T$ is at the same time also the algebra $H / F$, $F$ being
the ideal defined by the rule: ${h} \in F$ if and only if $\bar h\in T$.

If $T$ is a subset of $H$, then the filter of $H$ generated by $T$ consists
of elements of the form
$$
\forall (X) a_1 \wedge \cdots \wedge \forall (X) a_n \vee b,
    \quad a_i \in T,\ b \in H.
$$
Every filter is closed under existential quantifiers, while every ideal
is closed under universal quantifiers.  Besides, it can be proved that both
ideals and filters are closed under the action of the semigroup End$W$.

 $\bf 3.7.$ Theorem.
A set $T \subset H$ is a filter of
$H$ if and only if it satisfies the conditions
    \begin{enumerate}
\item[\rm 1.] $1 \in T$,
\item[\rm 2.] if $a \in T$, $a \to b \in T$, then $b \in T$,
\item[\rm 3.] if $a \in T$, then $ \forall (Y) a \in T$,
                    \quad $Y \subset X$.
    \end{enumerate}

Here, $a \to b = \neg a \vee b$.
We will consider two rules of inference in Halmos algebras:
    \begin{enumerate}
\item[1.]  {\it from $a$ and $a \to b$, infer $b$},
\item[2.] {\it from $a$, infer $\forall (Y) a$,\ $Y \subset X$.}
    \end{enumerate}
If $H$ is locally finite, then the second rule can be replaced by
    \begin{enumerate}
\item[2'.] {\it from $a$, infer $\forall x a$,\ $x \in X$.}
    \end{enumerate}

For every set $T$, one can consider the set of elements which
are inferred from $T$.

 $\bf 3.8.$ Theorem. If $T$ is a subset of $H$ containing
the unit, then the filter generated by $T$ is the set of all those
$h \in H$ which are inferred from $T$.

The notion of inferability in Halmos algebras agrees with that for
formulas in logic. The notion of a filter of a Halmos algebra corresponds
to the notion of a {\it closed} set of formulas.

Given a model $(G, \Phi, f)$,  we have $\hat f \colon U \to V_G$.
We can consider the corresponding filter $\Ker \hat f$ as the  {\it elementary
theory ($\Theta$-theory) of the  model}.

 $\bf 3.9.$ Theorem.
Let $T$ be a nonempty subset of $U$, $T' = K$ be
an axiomatizable class of models defined by the set $T$,
and $T^{''} = K'$ be the axiom set  of $K$ in $U$  (the closure of $T$).
Then $T^{''}$ is the filter of $U$ generated by $T$.

We now make some remarks on equalities in a Halmos algebra.

An equality in a Halmos algebra $H$ is a new nullary operation
$w \equiv {w'}$, where $w$, ${w'} \in W$ and are both of the same sort.

 $\bf 3.10.$ Definition.
An algebra $H \in \HA_\Theta$ is an algebra with
equalities if  all the operations $w \equiv {w'}$ are defined as
elements of $H$ and the following axioms  hold:
    \begin{enumerate}
\item[\rm 1.] $s (w \equiv {w'}) = ( s w \equiv s {w'})$,\quad $s \in \End W,$
\item[\rm 2.] $(w \equiv w) = 1$,\quad $w \in W$,
\item[\rm 3.] $(w_1 \equiv {w'}_1) \wedge \cdots \wedge (w_n \equiv {w'}_u) \,<\,
    w_1 \cdots w_n \omega \equiv {w'}_1 \cdots {w'}_w \omega$
    if $\omega \in \Omega$ is of an appropriate type,
\item[\rm 4.] $s^x_w a \wedge (w \equiv {w'}) \,<\, s^x_{{w'}} a$,\quad
        $ a \in H$, where $s^x_w$ takes $x$ into an element $w$ of
        the same sort and leaves every $y \neq x$ fixed.
    \end{enumerate}

In the algebra $V_G$, equalities are defined by the rule: $w  \equiv
{w'}$ is the set of all those $\mu\colon W \to G$ for which $w^\mu = {{w'}}^\mu$
in $G$. When considering the algebra $U$ with equalities, the symbol $\equiv$
is supposed to be added to  $\Phi$, and the  initial axioms of $\Theta$-logic
are supplemented by the standard axioms of equality. In this case the set
$\Phi$ may be empty.  Such an algebra $U$  arises from $\Theta$-logic with
equalities.

Since equalities are regarded as nullary operations,
any subalgebra of an algebra with equalities should contain all the
elements $w \equiv {w'}$.  The same remark concerns homomorphisms between
algebras with equalities.  If $H$ is an algebra
with equalities, then so is $H/ T$, where $T$ is a filter.

When $\Phi$ is empty, the algebra $U$ with equalities is the object that controls
all axiomatizable classes of algebras in $\Theta$. Relaying on $U$, we can
solve the problems mentioned before. The algebra $U$ with $\Phi$ nonempty
plays the same role for models.

We make here some further remarks.

Let $G=(G_i,\ i \in \Gamma )$ be an algebra in the variety $ \Theta $.
In order to investigate the elementary theories  of the models with
given $ \Phi $ and realized on $G$, it will be useful to introduce this
algebra into the language and the algebra of the corresponding calculus.
We assume that  $G$ is specified by generators and defining relations.

We fix the scheme of calculus, which includes the mapping
$n\colon X \to \Gamma $,
variety $ \Theta $ with the set of operation symbols $ \Omega $, and the
set of relation symbols $ \Phi $.

We denote the set of generators of $G$ by $M= \{ M_i,\ i \in \Gamma \} $ and
the set of defining relations by $\tau$. To each $a \in M$, we attach a
variable $y=y_a$. This way with each  $M_i$ a set of variables $Y_i$ is
associated, and we obtain a many-sorted set $Y=(Y_i,\ i \in \Gamma )$.
Let $W_G$ be the free algebra in $ \Theta $ over $Y$. The correspondence
$y_a \mapsto a$ yields  an epimorphism
$ \nu\colon W_G \to G $. The kernel $\Ker \nu = \rho $ is generated by $ \tau $,
and then we have an isomorphism $W_G/ \rho \to \Gamma $.

For every $a \in M_i$ and all $i \in \Gamma $, we add to $ \Omega $
a symbol of a nullary operation  $\omega_a$.
Let $ \Omega' $ be the new set of operation symbols, and let $ \Theta '$ be
the variety of  $ \Omega '$ -algebras defined by the identities  of the variety
$ \Theta $ and defining relations of $G$. Here, if
$w(y_{a_1}, \ldots ,y_{a_n})={{w'}}(y_{a_1}, \ldots ,y_{a_n})$
is a defining relation, then we must rewrite it in the form
$$
w(\omega_{a_1}, \ldots ,\omega_{a_n})={w'}(\omega_{a_1}, \ldots ,\omega_{a_n)}.
$$
We have no variables here, and this equality is also an identity.

Let  $W$ be the free algebra in $ \Theta $ over $X=(X_i,\ i \in \Gamma )$,
and let $W'$ be the free algebra in $ \Theta '$, also over $X$.

All $\omega_a$ are elements of $W'$. Suppose that $G'$ is a subalgebra of $W'$
generated by these elements. It can be proved that there exists a canonical
isomorphism $ \nu\colon G' \to G$ and that $W'$ is the free
product of $W$ and $G'$ in $ \Theta $, i.e. $W'=W*G'$.

In the old scheme we had a Halmos algebra $U$, and the new scheme
gives rise to the  algebra $U'$. There is a  canonical injection
$U \to U'$. The kernel of it admits a good description in the
case $G$ is finitely defined in  $\Theta $.

The algebra $G$ can be taken to be an  algebra in the variety $ \Theta '$,
and we can identify the homomorphism sets $\Hom(W,G)$ and $\Hom(W',G)$.
Simultaneously, we can identify algebras $V_G$ and $V'_G$. But $V_G$
is in the old scheme with the semigroup $\End W$, and $V'_G$--in the new
scheme with the semigroup $\End W'$.

Finally, take a homomorphism
$ \mu \colon W \to G$, which we identify with $ \mu \colon W' \to G$, and set
$s= \nu ^{-1} \mu \colon W' \to G'$, where $\nu\colon G' \to G$ is the
foregoing canonical isomorphism.
Because $G'$ is a subalgebra of $W'$,  $s \in \End W'$. Take
now some $ \tau \colon W' \to G$ and  $ x \in X$. We then have
$ \mu (x)= \nu \nu ^{-1} \mu (x)= \nu s(x)$. Since
$s(x) \in G'$ is a constant, $ \nu $ and $ \tau $ act equally on
$s(x)$, and $ \mu (x)= \tau s(x)$. It is so for all  $x \in X$;
therefore,  $ \mu = \tau s$ for each $ \tau $.

Given a model $(G, \Phi , f)$, $G\in \Theta $, where G is also considered as
an algebra in $ \Theta '$, we have the homomorphism $ {\hat f}: U' \to V_G$.
We consider every $u \in U$ as an element in $U'$, and the value of $u$ in $G$ is
$ {\hat f} (u)$. In order to calculate this value, we must answer the question
when $ \mu \in {\hat f}(u)$.

Together with  $u$, we also consider the element $su$. We can interpret it as
the result of substitution of the row $ \mu $ in the formula $u$.

 $\bf 3.11.$ Theorem.
$ \mu$ belongs to  ${\hat f}(u)$ if and only if
$ {\hat f}(su)=1$.

\begin{proof}
Let $ \mu \in {\hat f}(u)$. Take an arbitrary $ \tau$, $\tau s= \mu $.
Then $ \tau s \in {\hat f}(u)$ and $ \tau \in s {\hat f}(u)= {\hat f}(su)$.
Since $ \tau $ is arbitrary, we have $ {\hat f}(su)=1$.

Now let $ {\hat f}(su)=1$. Then $ \mu \in {\hat f}(su)=s {\hat f}(u)$, and
$ \mu s= \mu \in {\hat f}(u)$.
\end{proof}

The equality $ {\hat f}(su) = 1$ means that $su \in \Ker {\hat f}$, where
$\Ker {\hat f}$ is the elementary theory of the given model.

We see that, in the extended language, the value of an arbitrary formula
in a model can be calculated from the elementary theory.

\section*{\S2. Some applications}
\setcounter{paragraph}{0}

\paragraph{1. Closures of formula systems.}
\stepcounter{paragraph}
We begin with closures of systems consisting either of identities or
pseudoidentities. The discussion will be based on the Halmos algebra of
calculus, which we denote by $U$.

Take a scheme $ X$, $\Gamma$, $n\colon X\to \Gamma$, $\Omega $, $\Theta$,
all $ X_i$
in $ X $ being infinite. $U$ is a Halmos algebra with equalities in this
scheme, $W=(W_i,\ i\in\Gamma)$ is the free algebra over $X$ in $\Theta$. The
set of relation symbols $\Phi$ is empty.
It is well known that every set of identities can be presented in
the  algebra $W$. Then its closure is a fully
characteristic congruence  on $W$, generated by the set. The same set of
identities, and its closure, can  also be presented  in the algebra $U$.
We identify each identity  $w\equiv {w'}$ with the respective equality
element of $U$.

Let $T$ be a set of equalities of $U$, which are treated as identities
specifying some variety $K$.
From the result given in terms of algebra $W$ we can
conclude that the set $T$ is closed if the following conditions are fulfilled:
\begin{enumerate}
  \item[1.] $(w\equiv w)=1\in T$ and $T$ is closed under the semigroup $\End W$,
  \item[2.] if $w\equiv {w'}\in T$ and ${w'}\equiv {w'}'\in T$, then $w\equiv {w'}'\in T$,
  \item[3.] if $\omega\in \Omega$ and the type of $\omega$ is $(i_1,\ldots,i_n; j)$,
  then it follows   from
  $ w_1\equiv w_1'\in T ,\, \ldots, \,  w_n\equiv w_n'\in T $,
  $w_k,w_k'\in W_{i_k}$, that
 $ w_1 \cdots  w_n\omega\equiv w_1' \cdots w_n'\omega\in T $.
\end{enumerate}
As we know,  the formulas derivable from $T$ form the filter
generated by $T$. This filter necessary contains all identities of
$K$. Any of them is  derivable  from $T$, but, in general, the
derivation can contain not only  equalities. However, we may use
in the derivation only the rules listed above: the list is known
to be complete for equational formulas.

Now a few remarks on pseudoidentities follow. Here, $T$ is a set of
formulas--again, elements from  $U$--of the type
$$ w_1\equiv w_1'\vee \cdots \vee w_n\equiv {w'}_n.$$

Let $T$ be a closed set. This means that if $K$ is a pseudovariety in $\Theta $
defined by $T$, then all pseudoidentities of algebras from $K$ are in $T$.
It is obvious that
    \begin{enumerate}
\item[1.] $ 1\in T$ and $T$ is closed under the semigroup $\End W$,
\item[2.] if $u\in T$ and $v$ is a pseudoidentity in $U$, then $ u\vee v\in T$.
    \end{enumerate}
We shall formulate one more condition which also has to be fulfilled.
Suppose we are given pseudoidentities $ u_1, \ldots ,u_r $  from $U$,
where
$u_k$ is
$$
 (w_1^k \equiv {w_1^{k}}') \vee \cdots \vee (w^k_{n_k} \equiv {w^k_{n_k}}'),
\quad k=1,\ldots,r.
$$
Then we  define
a new  set of pseudoidentities, which we denote by $ u_1\circ \cdots \circ u_r$,
as follows.
First, we set $W\times W$ to be the union of all $W_i\times W_i$, $i\in\Gamma$.
Next, we take, for every $u_k$,
the subset $u{_k^\ast}$ of $W \times W$ consisting of all pairs
$(w^k_i,{{w'}_i}^k)$, $i=1, \ldots ,n$,
and make up the Cartesian product $V=u{_1^*} \times \cdots \times u{_r^*}$.
If $p=(p_1, \ldots , p_r) \in V$,
then we denote by $ \rho (p)$ the congruence on $W$ generated by all
$p_1, \ldots , p_r$,  $\rho (p)=( \rho _i,\ i \in \Gamma)$.
Also let  $ \psi $ be a function on $V$ such that $ \psi (p)$ is a pair
$(w,{w'})$ contained in some $ \rho _i,\ i \in \Gamma $, from $\rho (p)$.
We denote by $u=u( \psi )$ the pseudoidentity $\bigvee (w \equiv {w'})$ where
the disjunction is taken over all
$p \in V$. The set $u_1 \circ \cdots \circ u_r$ consists of all such
$u( \psi )$ for all $ \psi $. Now we can write out the third condition.

\begin{enumerate}
\item[3.] If $u_1,\ldots ,u_r  \in T$, then
        $u_1 \circ \cdots \circ u_r$ is a subset of $T$.
\end{enumerate}

Let us verify that this condition is satisfied if $T$ is closed.
Assume that $u_1, \ldots, u_r \in T$ and that
    $u=u( \psi ) \in u_1 \circ \cdots \circ u_r$.
Take a homomorphism $ \mu \colon W \to G$, $G \in K$. Then for every
$u_k$, $k=1, \ldots ,r$, we can find  $w_i^{k} \equiv w_i^{k\prime}$ so, that
$w_i^{k\, \mu }=w{_i^{k\prime\, \mu }}$ in $G$.
Denote $(w{_i^k},w{_i^{k\prime}})$ by $p_k$ and take $p=(p_1, \ldots , p_r)$.
All $p_k$ are in $\Ker \mu $, hence $ \rho (p) \subset \Ker \mu $. By
the definition $ \psi (p)=(w,{w'})$ also lies in $\Ker \mu $, and
we have $w^ \mu =w^{' \mu }$ in $G$. This means that the pseudoidentity
$u=u( \psi )$ is valid in $G$.

The converse also can be proved, and so we have

 $\bf 1.1.$ Theorem.
The set of pseudoidentities $T$ is closed if and only if
$T$ satisfies the conditions 1, 2, 3.

We may consider the conditions 1, 2, 3 as rules of inference, so that
we can construct the closure $T'$ of every $T$.

This result belongs to A.~Kushkuley. My participation in obtaining it is not
great: I only consider the many-sorted case, and put it in terms of Halmos
algebras.
A.~Kushkuley and S.~Rosenberg (the former is now in the USA, and the latter is
in Jerusalem, both of them are from Riga) have also obtained a generalization,
which we are going to consider.

We  shall deal with closures of universal formulas.  The universal formulas
we consider are formulas of type
$$
(w_1 \equiv w_1') \vee \cdots \vee (w_n \equiv w_n') \vee
    (v_1 \neq v{_1'}) \vee \cdots \vee (v_m \neq v{_m'}),
    \quad w,v \in W.$$
A set of such formulas describes an  universal class of algebras.

Given a set $T$ of universal formulas from $U$, we must construct
the closure for $T$.

For every universal formula $u$ of the above kind, we take the set $u^+$
consisting of
all pairs $(w_1, w{_1'}), \ldots ,\-(w_n, w{_n'})$, and let $u^-$ be the
set of pairs $(v_1, v{_1'}), \ldots,$\\ $(v_m, v{_m'})$. The union of
$u^+$ and $u^-$ is denoted by $u^*$.

As before, for every set of universal formulas $u_1, \ldots , u_r$,
we construct a new set $u_1 \circ \cdots \circ u_r$ of such formulas.
At first, we take the Cartesian product
$u{_1^*} \times \cdots \times u{_r^*}=V$.
If $p=(p_1, \ldots , p_r) \in V$, then $p^+$ is the set of "positive"
pairs, used in the notation of $p$, and $p^-$ is the  set of the
"negative" pairs. If $M$ is a subset of $W \times W$, then $ \rho (M)$ is
the congruence of $W$ generated by $M$.
Now let
$u \in u_1 \circ \cdots \circ u_r$ if for every $p \in V$
$$
 \rho (u^{-} \cup p^+) \cap (u^+ \cup p^-) \neq \emptyset.
$$

This definition generalizes that given for pseudoidentities.

Now let us show that if each of $u_1, \ldots , u_r$ is valid in $G$ and
$u \in u_1 \circ \cdots \circ u_r$, then $u$ is also  valid in
$G$.

Given a homomorphism $ \mu \colon W \to G$, we construct an
element $p=(p_1, \ldots , p_n)$\\ $\in V$  depending on  $ \mu $.
 Take any $u_k$:
$$
(w_1^k \equiv w_1 ^{k\prime}) \vee \cdots
    \vee (w_{n_k}^k \equiv w_{n_k}^{k\prime}) \vee
     ( v_1^k \not\equiv v_1^{k\prime}) \vee \cdots \vee
     (v_{m_k}^k \not\equiv v_{m_k}^{k\prime}).
$$
 $k=1, \ldots, r$.
 As $u_k$ is valid in $G$,  we have, for a given $\mu $, some pair
$( w_i^k, w_i^{k\prime})$ with $ w_i^{k\,\mu}=w_i^{{k\prime}\,\mu}$
 or some other $(v_j^k,v_j^{k\prime})$ with
 $ v_j^{k\,\mu} \neq  v_j^{{k\prime}\,\mu}$.
We take any one of them to be $ p_k$, and this way we construct
$p=(p_1, \ldots, p_r)$. For every pair $(w,{w'})$ in $p^+$,
 we have $ w^\mu ={ {w'}}^\mu$, and for every $(v,v') $
 in $p^-$--$ v^\mu \neq {v'}^\mu $.
Assume now that the pair $(w,{w'})$ belongs to the intersection
$$
         \rho (u^- \cup p^+) \cap (u^+ \cup p^-).
$$
All pairs in $ p^+$ belong to $\Ker\mu $. If this is  true also of all
pairs in $ u^- ,$  then $\rho=\rho(u^-\cup p^+)\subset \Ker\mu $.
 Then $(w,{w'})$ also is in $\Ker\mu$, $ w^\mu=w^{\prime\,\mu }$ and $(w,{w'})
 \notin p^-$.
 We have $(w,{w'})\in u^+$. If some $(w,{w'})$ from $u^-$ is not in $\Ker\mu $,
 then $ w^\mu \neq w^{\prime\,\mu }$. In all cases,
 the homomorphism $\mu$ belongs
 to the value of the formula $u$ in $G$. This holds for every $\mu$, and
so $u$ is valid in $G$.

Now we can formulate the main result.

 $\bf 1.2.$ Theorem. {\bf\cite{KuR}}
The set $T$ of universal formulas is closed if and only
  if the following conditions are satisfied:
    \begin{enumerate}
\item[\rm 1.]  $1\in T$, and $T$ is closed under the action of the semigroup
    $\End W$.
\item[\rm 2.] If $u \in T$ and $v$ is universal in $U$, then $u\vee v\in T$.
\item[\rm 3.] If $  u_1, \ldots, u_r\in T $,
    then all $ u\in u_1 \circ \cdots \circ u_r $ also belong to $T$.
    \end{enumerate}

Theorem 1.1 is a particular case of this theorem; we only must everywhere
 delete "negative" parts.
  The closure of every set of universal formulas can be constructed in virtue of
Theorem 1.2.

Now we shall consider closures of sets of quasi-identities.
A quasi-identity is an element $ u\in U $ of the type
    $$(w_1 \equiv w_1')\wedge \cdots  \wedge(w_n \equiv w_n')\to (w\equiv {w'}).$$
This is a universal formula of a particular kind. We are interested in the
question how to obtain  the closure of a set of quasi-identities. The problem
was investigated
by  R.~Quackenbush \cite{Qua}. We translate his result in terms of the algebra $U$.
We rewrite $u$ in the form
  $ u_0 \to (w\equiv {w'})$, where $u_0$ is $(w_1\equiv w_1')\wedge \cdots \wedge
  (w_n \equiv w_n') $.

 $\bf 1.3.$ Theorem.
The set of quasi-identities $T$ is closed if and only if it
  satisfies the following conditions:
    \begin{enumerate}
\item[\rm 1.]  $ 1\in T$, and the set $T$ is invariant under the semigroup $\End W$.
\item[\rm 2.] If $u$ is
    $(w_1\equiv w_1')\wedge \cdots \wedge(w_i\equiv w_i')\wedge \cdots
          \wedge(w_n\equiv w_n')$,
      then $u\to (w_i\equiv w_i') \in T $.
\item[\rm 3.] If $ u_0 \to (w\equiv {w'})\in T$ and $ u_0\to ({w'}\equiv {w'}')\in T $,
  then $ u_0 \to (w\equiv {w'}')\in T $.
\item[\rm 4.] If $\omega \in \Omega $ and the type of $\omega $ is
    $(i_1, \ldots, i_n;j)$,  and if $ u_0\to (w_k\equiv w_k')\in T$,
    $k=1, \ldots, n$, $w_k, w_k'\in W_{i_k} $,
    then
    $u_0\to (w_1 \cdots  w_n\omega\equiv w_1'\cdots w_n'\omega)\in T $.
\item[\rm 5.] If $ u_0\to (w_i \equiv w_i')\in T$, $i=1, \ldots, n $, \ and if
 $(w_1\equiv w_1') \wedge \cdots \wedge (w_n\equiv w_n') \to (w\equiv {w'})\in T$,
 \ then $ u_0\to (w\equiv {w'})\in T $.
    \end{enumerate}

Necessity of this is obvious, and sufficiency in this theorem
  and in 1.1 and 1.2  is based on the following scheme. Let the set $T$
  satisfy the conditions of the theorem, and assume that the formula $ u\in U$
  is not in T.
Then an algebra $G$ satisfying $T$, but not  $u$, can be found.
  The main problem is to construct such $G$.

 In \cite{Qua} the closure problem is dealt with for universal formulas as well.
 The  result obtained there differs from that of Kushkuley and Rosenberg. (Both
 results were discovered at the same time.) Moreover,  implicative classes are
discussed in \cite{Qua}, and the rules of inference include
also this one: from $ u_0\to 0 $, infer $ u_0 \to (w\equiv {w'})$ for every $w$
and ${w'}$  of the same sort.
On this question, see also \cite{Kel,Sel,H}.

\paragraph{2. Quasigroups.}
\stepcounter{paragraph}
We shall consider here a known problem in the quasigroup theory.

A quasigroup  is a group without associativity and, of course,
without unit. If the unit is added, then we have a loop. More precisely,
a quasigroup $Q$ is a set with one binary operation of multiplication,
and  equations $ax=b$ and $ya=b$ are solved in $Q$ uniquely. We
introduce two additional operations: $x=a\backslash b$ and  $y=b/a$.  So the class of
all quasigroups is a variety with the specifying identities
$$
x(x\backslash y)=y,\quad (x/y)y=x,\quad (xy)/y=x,\quad x\backslash (xy)=y.
$$
 Adding a nullary operation 1 with identities $1x=x$ and $x1=x$,
we get the variety of loops. The variety of groups arises when we add
the associativity requirement.  This is the definition of the variety of groups
in the quasigroup signature.
So, the notion of a quasigroup, as well as that of a semigroup, generalizes the
notion of a group, but they  do this in different ways.

Quasigroups have arisen from some problems of geometry. Loops also
have applications in algebraic geometry, but these applications
are not like  to those of groups. In  particular,  we cannot speak
about  representations of quasigroups as quasigroups of
permutations. In the theory of quasigroups, along with
homomorphisms, the homotopies are used. A homotopy is a triplet
$$ \mu =( \mu _1,  \mu _2,  \mu _3) \colon Q \to Q',$$ such that
$$x^{\mu _1} y^{\mu _2}=(xy)^{\mu _3}.$$

We can also speak of the category of quasigroups with homotopies as
morphisms.
If all maps $ \mu _1, \mu _2$ and $\mu _3 $  are bijective, then $ \mu $
is an isotopy. In  geometric applications,  quasigroups are
considered up to isotopies.

In the group theory the notion of an isotopy is not of interest. For groups, an
isotopy  reduces to an isomorphism, and a homotopy--to a homomorhism. On the
other hand, a quasigroup  isotopic to a  group may be not a group.

If we take the class $K$ of all quasigroups isotopic to  groups,
then $K$ is a variety closed under isotopy. It differs from the  variety
of all quasigroups.

We have the following general result.

 $\bf 2.1.$ Theorem.
Let $ \Theta $ be a variety of groups, and $ \Theta '$ be the
class of quasigroups isotopic to groups from $ \Theta $.
Then:
    \begin{enumerate}
\item[\rm 1.] $ \Theta '$ is a variety of quasigroups,
    \end{enumerate}
\noindent and
    \begin{enumerate}
\item[\rm 2.] $ \Theta '$   is invariant under isotopy.
    \end{enumerate}

Long ago the geometrical applications prompted the following
problem, which was stated, as it seems, by V.D.~Belousow.
Under what conditions a variety of quasigroups is closed under isotopies?
Every variety is closed under
isomorphisms, but not every one is closed under isotopies--for example,
any variety of groups. Moreover,
if $ \Theta $ is a variety of quasigroups
and $ \Theta '$ is the class of quasigroups isotopic to the quasigroups from
$ \Theta $, then the class $ \Theta '$ is closed under isotopies, but
it may be not a variety.

This problem was solved some years ago by A.A. Gvaramia during his
postdoc in Riga, and even in a more general setting--for an
arbitrary axiomatizable class of quasigroups. The solution was
given using Halmos algebras \cite{Gv,GvPl}.

I bring here the sketch of the solution. First, together with the
category of quasigroups with homotopies as morphisms, we also take
the category of the three-sorted quasigroups. Its objects have the form
$A=(A_1,A_2,A_3)$, and there is defined an operation $*\colon A_1
\times A_2 \to A_3$
such that every pair of elements in $a_1*a_2=a_3$ uniquely determines
the third one. As in quasigroups, we have the inverse operations
$*^{-1} = \backslash :A_1 \times A_3 \to A_2$ and $ ^{-1}*=/:A_3
\times A_2 \to A_1$.
We call such objects also invertible automata. Morphisms in this
category are the homomorphisms
$$
 \mu=(\mu_1, \mu_2, \mu_3) \colon A=(A_1,A_2,A_3)\to A'=(A_1',A_2',A_3').
$$
They are coordinated with all three operations. We connect a
regular automaton $\at Q=(Q,Q,Q)$ with every quasigroup $Q$, and
every homotopy $ \mu=(\mu_1, \mu_2, \mu_3)\colon Q\to Q'$ gives a
homomorphism
     $$\at\mu=( \mu_1, \mu_2, \mu_3)\colon \at Q\to \at Q'.$$
So we can consider the category of quasigroups with homotopies as a subcategory
of the category of invertible automata.

It is easily proved that, for every automaton $ A=(A_1,A_2,A_3)$, there is a
quasigroup $Q$ with an isomorphism

    $$\mu=(\mu_1, \mu_2, \mu_3): A\to \at Q.$$

 Note that if $\mu=(\mu_1, \mu_2, \mu_3)\colon Q\to Q'$ is a homotopy, then
the image of $\mu$ in $Q'$ is a three-sorted subquasigroup of the
quasigroup $Q'$. It is obvious that if $\Iks$ is a variety of
invertible automata, and $q\Iks$ consists of all quasigroups $Q$
such that $\at Q \in \Iks$, then $q \Iks$ is a variety of
quasigroups, closed under isotopies.

Let, on the other hand, $\Theta$ be a class of quasigroups and $\at \Theta $
consist of automata each of which is isomorphic to some $\at Q$ with
 $ Q \in \Theta $. So we have:
    \begin{enumerate}
\item[1.] If $\Iks$ is an abstract class of automata, then
            $ \at(q \Iks)=\Iks $,
\item[2.] If $\Theta $ is closed under isotopy, then $ q(\at \Theta )=\Theta $.
    \end{enumerate}
This leads to a connection between varieties of quasigroups
closed under isotopies, and varieties of invertible automata.

We are interested in identities and arbitrary formulas which define
classes closed under isotopies.

Let $F=F(X)$ be the free quasigroup over the set $X$, and let
$X=X_1 \cup X_2 \cup X_3$ be some partition of $X$. We have also
the triple $X=(X_1,X_2,X_3)$, and it  generates an automaton $
\Phi =(\Phi _1, \Phi _2, \Phi _3)$ in $F$. All intersections $
\Phi _i \cap \Phi _j$, $i \neq j$, prove to be empty, and $ \Phi $
is the free automaton over $X=(X_1,X_2,X_3)$.

 Let now $ \Theta $ be the variety of all quasigroups and $ \Theta '$--the
 variety of all invertible automata. We can consider calculi in these varieties,
and hence we have Halmos algebras $U$ and $U'$ over $X$ and
$X=(X_1,X_2,X_3)$ respectively. The above notes allow us to
consider an injection $U' \to U$, and so we can take in $U$ some
three-sorted formulas which at the same time are thought of as
one-sorted. We call these formulas special.

Given a quasigroup $Q$, we have the canonical homomorphisms
$$
f_Q\colon U \to V_Q
$$
and
$$ f_{\at\! Q}:U' \to V_{\at\! Q} .
$$
Now take a formula $u \in U'$ and consider it simultaneously as a special
formula in $U$.
It proves that
$$
f_Q(u)=1 \in V_Q\; \Leftrightarrow\;  f_{\at\! Q}(u)=1 \in V_{\at\! Q}.
$$

Assume now that $u$ is valid in $Q$ and that $Q'$ is isotopic to $Q$.
Then $u$ is valid also
in $\at Q$ and $\at Q'$. Hence, $u$ is valid in $Q'$.
So the class of all special formulas is closed under isotopies. The converse
is proved by some additional reasoning, also in terms of Halmos algebras.
We thus have

 $\bf 2.2.$ Theorem.
A formula is closed under isotopies if and only if it is
equivalent to a special one. The same holds for sets of formulas.

Note also that, for every formula, some special derived
formula can be taken, and this leads to a construction of basis of
special formulas.

Let us make some remarks about the derived formulas.

Set $X=\{\{u\},\{v\},X'\}$, and let $\Phi =( \Phi _1, \Phi _2, \Phi _3)$ be
an automaton over  $X$. Introduce an operation $ \circ $ on $\Phi_3$
by the rule:
$$
f_1 \circ f_2=(f_1/v)(u \backslash f_2).
$$
Then $ \Phi _3$ becomes a quasigroup. Let $F=F(X)$  be the free quasigroup over
$X=\{u\} \cup \{v\} \cup X'$, and $F(X')$--the free quasigroup over $X'$.
The identity map $X' \to X'$ gives a homomorphism $F(X') \to \Phi _3 $,
where $ \Phi _3$ is a quasigroup with respect to  the operation $ \circ $.
For every $w\in F(X')$, the corresponding ${\bar w} \in \Phi _3 $
is an {\it automata element}. The identity $ w \equiv {w'}$ gives a
derived automata identity $ {\bar w} \equiv {\bar w}'$, and likewise  for
arbitrary formulas. For example,   the automata identity
corresponding to the identity $(xy)z\equiv x(yz)$ is the one
$$
((x/v)(u \backslash y)/v)(u\backslash z) \equiv (x/v)(u \backslash ((y/v)
(u \backslash z))),
$$
which specifies the variety of quasigroups isotopic to groups.
If $ \Theta $ is the  variety of groups which satisfy some set of identities
$(w_\alpha \equiv w_\alpha',\  \alpha \in I) $, then the set
${\bar w_\alpha} \equiv {\bar w}_\alpha'$,\ $ \alpha \in I $
determines the variety
$ \Theta '$ of quasigroups isotopic to groups from $ \Theta $.
The main result yields characteristic conditions for varieties,
quasivarieties, pseudovarieties, universal classes, etc. of quasigroups closed
under isotopies.

\paragraph{3. Algebraic logic in group representations.}
\stepcounter{paragraph}
A representation is considered to be a pair $\rho = (V, G)$, where $V$ is a
$K$-module, $K$ is a commutative ring with unit, and $G$ is a group acting on
$V$. Let $\mu\colon G \to \Aut V$ be the corresponding homomorphism. Varieties of
representations are studied in \cite{PlV}.

Halmos algebras can be applied in the situation when $\Theta$ is the variety of
representations over a given ring $K$. Let $X$ be an infinite set of variable
that run over $V$, and $Y$ be an infinite set of variables that run over the
acting group $G$. We then have $F = F(X)$ and $W = (XKF, F)$.

There are two types of elementary formulas:
    \begin{enumerate}
    \item[1.] $x_1 \circ u_1 + \cdots + x_n \circ u_n \equiv 0$, \quad $u_i \in
    KF$,
    \item[2.] $f \equiv 1$, \quad $f \in F$.
    \end{enumerate}
We call the formulas of the first type the action formulas. An action formula
is constructed from the elementary action formulas by means of Boolean
operations and quantifiers with variables from $X$: we do not quantify variables
ranging over the group. In particular, we can speak of action identities,
quasi-identities, pseudoidentities and so on. In general, universal formulas
have the form
\[
(w_1 \equiv {w'}_1) \lor \cdots \lor (w_n \equiv {w'}_n) \lor
    (v_1 \not\equiv v'_1) \lor \cdots \lor (v_m \not\equiv v'_m),
\]
where all $w, v \in W$.

We have mentioned in \cite{PlV} that every saturated variety of representations
can be defined by action identities. This is true also of quasivarieties and,
possibly, of pseudovarieties and universal classes of representations as well.

 $\bf 3.1.$ Proposition.
Let u be an action formula, and let $\rho = (V,G)$,
$\bar{\rho} = (V, \bar{G})$ be a representation and the corresponding
faithful representation, respectively. Then the formula $u$ is valid in $\rho$
if and only if
it is valid in $\bar{\rho}$.

The proof is direct, by using the homomorphisms
\[
f_\rho\colon U \to V_\rho \quad {\rm and} \quad f_{\bar{\rho}}\,\colon U
\to V_{\bar{\rho}}.
\]

 $\bf 3.2.$ Proposition.
Let $u$ be an action  formula, $\rho = (V,G)$ be a
representation and $\rho' = (V,H)$ be its subrepresentation, where $H$ is a
subgroup of $G$. Then  $u$ is valid in $\rho'$ whenever it is valid in
$\rho$.

The next propositions follows from the two preceding ones.

 $\bf 3.3.$ Proposition.
Let $T$ be a set of action formulas and $T'$--a class
of representations defined by $T$, $T' = \Iks$.  Then $\Iks$ is saturated and
right hereditary.

Recall that a class $\Iks$ is called saturated under the condition that
the representation $\rho = (V,G)$ belongs to $\Iks$ if and only if
the representation $\bar\rho = (V, \bar{G})$
belongs to $\Iks$. This situation can also be described as
follows. Representations $(V,G)$ and $(V',G')$ are said to be similar if the
respective faithful representations are isomorphic. An abstract class $\Iks$
of representations is saturated if and only if $\Iks$ is invariant under
passing to similar representations.

Right hereditary means here that $(V,G) \in \Iks$ implies $(V,H)
\in \Iks$ if $H$ is a subgroup of $G$ and $(V,H)$ is a
subrepresentation of $(V,G)$. Proposition 3.3 gives sufficient
conditions for determining a class by action formulas. The
following question is on necessary and sufficient conditions.

 $\bf 3.4.$ Problem.
Is it true that an abstract class of representations can be
determined by a set $T$ of action formulas if and only if $\Iks$ is
axiomatizable, saturated and right hereditary?

This question, as it seems to us, does not appear to  be difficult. We must use
the known conditions of axiomatizability, and pass to the class of all
representations of the free group $F$ of countable rank in the given $\Iks$.

It may also happen that the property of right hereditariness follows from the
first two conditions.

A set of formulas $T \subset U$ is said to be saturated if the class of
representations $\Iks = T'$ is a saturated class.

 $\bf 3.5.$ Problem.
Is it true that $T$ is saturated if and only if it is
equivalent to some set $T_1$ of action formulas?

This problem is related with the preceding one. Two sets, $T$ and $T_1$ are
equivalent when the classes $T'$ and $T'_1$ coincide. Syntactically this means
that both $T$ and $T_1$ generate the same filter.

Independently we can speak of equivalence of sets of identities,\!
quasi-identities, pseudoidentities, and universal formulas. For conditions of
such equivalences in terms of the algebra $U$, see the subsection 2.1 above.
Hence, Problem 3.5 can be specified having in mind such sets of formulas. We
can also
consider characterizations of a single saturated formula.

Let us conclude with the following result. Given a class of representations
$\Iks$, we denote by $\vec{\Iks}$ the class of groups admitting a
faithful representation in $\Iks$.

 $\bf 3.6.$ Proposition.
If $\Iks$ is a universal and saturated class of representations over $K$,
then the class of groups $\vec{\Iks}$ is a universal class.

This means that such $\vec{\Iks}$ admits a description by universal formulas
in the group theory logic. In particular, the class $\vec\Iks$,
for any pseudovariety
of representations, is characterized by universal formulas in the group theory
logic.

\paragraph{4. Databases.}
\stepcounter{paragraph}
Constructing a database model presupposes that given are a data algebra
$ G=( G_i,\ i\in\Gamma) $, a  set of relation symbols $ \Phi $, and a set of
states $ F $. Every
 $ f \in F $ is a function that realizes every $\varphi \in \Phi $ as
 a relation on $ G $. For every $f \in F$, the triple $ (G, \Phi, f) $
 is a model.
   We must also take some scheme relatively to which databases are to be
   considered.
 For this purpose we take the scheme in which Halmos algebras were defined.
 In particular, $ G\in \Theta $, and every $\varphi \in \Phi $ has a type
 $\tau=(i_1,...,i_n)$, $i\in \Gamma $.
  In this scheme the triple $(G,\Phi, F)$ presents a database, but this
  is  not yet a database model. We call the triple a passive database.

Database receives queries  and produce replies to them. The
queries are written as formulas, i.e. elements of algebra $ L\Phi
W$. The same query
 can be written out in different equivalent ways. This equivalence is the
 same which we got by the rule of Lindenbaum-Tarski. Hence, we
 must consider a query as a class of equivalent formulas, and then the algebra
 of  queries is the Halmos algebra $U$. The algebra of replies is also a Halmos
 algebra. It is the algebra $ V_G $ constructed for $G$ in the given scheme.
  To every $ f\in F$  a homomorphism ${\hat f}\colon U\to V_G $ corresponds.
${\hat f }(u) $
 is the  value of a formula $u$ in $G$, and at  the same time this is the reply
to the query $u$ in the state $f$. We also write $ {\hat f}(u)=f*u $. If $u$
 is a formula in $ L\Phi W$, and ${\bar u}$ is the corresponding element in $U$,
 then $ f*u=f*{\bar u}$. So we have a database $ (F,U,V_G) $ with an
 operation $*$.
 The algebra $U$ here does not depend on $G$ and $F$, it depends only
 on $\Phi $ and the scheme. In this sense, $U$ is the universal query algebra.
 This algebra can be compressed, and $V_G$ can be reduced.
   First of all, take a subalgebra $R$ of $V_G$  generated by all $ f*u$,
$ f\in F$, $u\in U $. This gives us $(F,U,R)$. The next step consists in
specifying the filter $T$ in $U$
by the rule: $ u\in T $ if $ f*u=1 $ for all $f\in F$.
 Let $ Q=U/T $; this way we obtain a  reduced database $(F,Q,R)$. If here
 $u\in U$, and $ q={\bar u}$ is the corresponding element in $Q$, then $
f*q=f*{\bar u}=f*u $.
   The database $(F,Q,R)$ is constructed in the given scheme from the passive
 database $(G,\Phi,F)$. We call $(F,Q,R)$ an active database, or an
 algebraic model of a database.

  Using the model $(F,Q,R)$, in which Halmos algebras
play an essential role, we can  solve various database problems.
However, $ HA_\Theta $ is very hard to be used for computer applications.
That is why we must return,  in our final conclusions, to passive
 databases.

Let us denote by $F_G$ the system of all possible states of the collection
$\Phi$ in the algebra $G$. Then we arrive at the universal database
\[
(F_G, U, V_G).
\]
If $\delta\colon\! G \to G'$ is a surjective homomorphism, then it produces
an injective database homomorphism
\[
\delta_{*}\colon (F_{G'}, U, V_{G'}) \to (F_G, U, V_G).
\]
For all $f \in F_{G'}$ and $u \in U$,
\[
(f*u)^{\delta_{*}} = f^{\delta_{*}}*u.
\]
See \cite{BPl3, BPlhand}.

All this will be in use in \S5.

\section*{\S3. Algebraic varieties and varieties of algebras}
\setcounter{paragraph}{0}
\paragraph*{1. Basic concepts.}
\stepcounter{paragraph}
 The present and next section relate to the level of equational logic, and they
 are not immediately connected with algebraic logic. In the following we,
 however, shall move to the universal logic level, and constructions related to
 algebraic logic, and even to databases, find essential applications there.

 At present, we are interested in equations and identities over arbitrary
 algebraic structures. Here, algebraic varieties correlated with arbitrary
 varieties of algebras are considered.

 We first remind some matters well-known in algebraic geometry
 ~\cite{ZS,Hart,Shap}.

 Let $P$ be a field and $K$ -- its extension. We consider the ring of
 polynomials $R=P[x_1,\ldots,x_n ]$, and take the affine point space $K^{(n)}$.
 If $T$ is a collection of polynomials from $R$, then it is attached the
 algebraic variety $T'=A$ in $K^{(n)}$ that consists of all those points
 $a=(a_1,\ldots,a_n)\in K^{(n)}$, $a_i\in K$, nullifying the polynomials from
 $T$. The same variety is specified by the ideal generated by $T$. If, on the
 other hand, $A$ is a subset of $K^{(n)}$, then it induces the ideal $T=A'$ of
 $R$ that consists just of the polynomials nullified by the points from $A$. We
 arrive at a Galois correspondence between subset of $K^{(n)}$ and collections of
 polynomials from $R$. Fore every $A$, we have $A'=T$ and $A''=T'$; $A''$ is the
 Galua closure of $A$. It coincides with the intersection of the algebraic
 varieties including $A$, as is itself an algebraic variety. In a like manner,
 $A=T'$ and $A'=T''$ for every $T$. $T''$ is the closure of the collection $T$
 and is always an ideal of $R$. If $T$ is an ideal, then links between $T$ and
 $T''$ are revealed by the Hilbert theorem on zeros. According to the theorem,
 if the fields is algebraically closed, then $T''=\sqrt{T}$, where $\sqrt{T}$
 is the radical of the ideal $T$, i.e. the set of those $\varphi\in R$ with
 $\varphi^n\in T$ for some $n$. We stress that in this case the passage from $T$
 to $T''$ does not depend on the field $K$ and is completely determined by the
 algebraic closeness of it. In general, connections with $K$ can be more
 essential, and it would be misleadingly to write $T''_{K}$. In this context two
 fields $K_1$ and $K_2$ both being extensions of $P$, could be called equivalent
 in the geometry under consideration if $T''_{K_1}=T''_{K_2}$. If $K_1$ and
 $K_2$ are algebraically closed, then they are equivalent; there is a little to
 say about these things in general case.

 The intersection of any collection of algebraic varieties is again an algebraic
 variety, and so is the sum of a finite number of varieties. This gives a topology
 on $K^{(n)}$ known as the Zariski topology.

 Now we shall discuss another point of view on the same. We proceed from the
 variety $\Theta$ of all commutative and associative algebras with unit over the
 field $P$. $R$ is the free algebra of this variety over the set of variables
 $X=\{ x_1,\ldots,x_n\}$; the field $K$ can also be considered as an algebra
 from $\Theta$. Every point $(a_1,\ldots,a_n)=a\in K^{(n)}$ specifies a mapping
 $\mu\colon X\to K$, $\mu (x_i)=a_i$, $i=1,2,\ldots,n$. The mapping determines
 an algebra homomorphism $\mu\colon R\to K$. Therefore, we may identify the
 space $K^{(n)}$ with the homomorphism set $\Hom (R,K)$. Here, the "point" $\mu$
 is a root of the polynomial $\varphi =\varphi(x_1,\ldots,x_n)$ if
 $\varphi\in\Ker\mu$. An algebraic variety is now treated as a subset of $\Hom
 (R,K)$, and the set $\Hom (R,K)$ can be thought as an affine space.

 Let us rewrite the Galois correspondence considered above in these new terms.
 It is easily seen  that now
$$T'=\{\mu,\ T\subset\Ker\mu\},$$
$$A'=\bigcap\limits_{\mu\in A}\Ker\mu.$$
 This shows that $A'$ is always an ideal.

 Now we can take up the general viewpoint.

 Assume that $\Theta$ is any variety of algebras of the signature $\Omega$.
 The algebras may be many-sorted; then $\Gamma$ stands for the set of sorts. Let
 $X$ be a set of variables, and let $W=W(X)$ be the algebra from $\Theta$ free
 over $X$. Take an algebra $G\in\Theta$ and consider the set $\HOM$, which is
 now treated as an affine space. We are going to define a Galois correspondence
 between binary relations $T$ on $W$ and subsets of $\HOM$. $T$ is the set of
 pairs $(w,{w'})$ with $w$ and ${w'}$ of the same sort, and $wT{w'}$, as usually, means
 that $(w,{w'})\in T$. The equation $w={w'}$ is related with every pair $(w,{w'})$,
 and one can also consider the formula $w\equiv {w'}$ (as an element of the Halmos
 algebra $U$).

 Every $\mu=(\mu_i,\ i\in T)\colon W=(W_i,\ i\in\Gamma)\to G=(G_i,\
 i\in\Gamma)$ has the kernel $\Ker\mu=(\Ker\mu_i,\ i\in\Gamma)$, where
 $\Ker\mu_i$ is the kernel equivalence of the mapping $\mu_i$, i.e. the set of
 pairs $(w,{w'})$, $w,{w'}\in W_i$, with $w^{\mu_i}={w'}^{\mu_i}$ or, what is  the
 same, $w^{\mu}={w'}^{\mu}$. We also consider the kernel $\Ker\mu$ as the union of
 all $\Ker\mu_i$, $i\in\Gamma$. On the other hand, $\Ker\mu$ is a congruence of $W$.

 Now let $A$ be any subset of $\HOM$. We set
$$A'=T=\bigcap_{\mu\in A}\Ker\mu.$$
 If $T$ is a binary relation on $W$, then $T'=A$ is defined by the rule
$$A=\{\mu,\ T\subset\Ker\mu\}.$$
 $A=T'$ is the "algebraic variety" in $\HOM$ specified by the collection $T$,
 while $T=A'$ is always a congruence of $W$. So we have get a Galois correspondence.
 Every set $A$ can be closed up to an algebraic variety $A''$, and every $T$ --
 up to an congruence $T''$. Where $T$ is a congruence, links between $T$ and
 $T''$ are revealed by an appropriate "Hilbert's Nullstellensatz".

 Clearly, the intersection of a collection of algebraic varieties is  a variety
 again, and if $A$ is a subset of $\HOM$, then the closure $A''$ is the least
 variety including $A$. In the present general situation, however, the union of
 two varieties can fail to be a variety, and the closure of a sum of sets
 generally differs from the sum of the closures. Later on, we shall generalize
 the notion of an algebraic variety and improve this shortcoming.

 As noticed above, a pair $(w,{w'})$ can be regarded as an equation $w={w'}$. Then
 the statement $(w,{w'})\in\Ker\mu$ means that the point $\mu$ satisfies this
 equation.
 The equation determines an algebraic variety in $\HOM$, and if the latter
 coincides with $\HOM$, then the equation becomes an identity of $G$.

 A binary relation $T$ specifies the {\it variety of algebras} $G\in\Theta$
 what the identity $w\equiv {w'}$ is valid in for every $(w,{w'})\in T$. This is a
 variety in $\Theta$. The same $T$, for every particular $G\in\Theta$, specifies
 an {\it algebraic variety} thought of as a subset of $\HOM$. This is the
 connection between varieties of algebras and algebraic varieties.

 Let $T(G)$ stands for the verbal congruence of all identities of any algebra
 $G\in\Theta$ in $W$. For an arbitrary collection $T$, the equality
 $T'=\HOM$ means that $T\subset T(G)$.

 Let us discuss another specific situation.

 Recall that the unit congruence is one in which any two elements of the same
 sort are identified. The zero congruence presupposes that elements are
 equivalent only if they are equal. The zero congruence is included in any $T$, and
 any $T$ is included in the unit congruence. Let $T_0$ stands for the zero
 congruence and $T_1$--for the unit congruence on $W$. Then, obviously,
 $T'_0=\HOM$.
 Moreover, $\mu\in T'_1$ if $T_1\subset\Ker\mu$, i.e. if $\Ker\mu$ is the unit
 congruence. Such $\mu$ need not exists, i.e. $T'_1$ may be the empty set.
 Clearly, $T'$ can be empty for other $T$ as well.

 Now we single out the case when $G$ has a one-element subalgebra $H$. If $G$ is
 many-sorted, then the subalgebra is of the form $H=(H_i,\ i\in\Gamma)$ with
 all $H_i$ singletones. $H$ determines a homomorphism $\mu_0\colon W\to G$, and
 $\Ker\mu_0=T_1$. Here, $T'_1=\{\mu_0\}$. If there is no such $H$, then $T'_1$
 is  empty. If a one-element $H$
 exists, then, for every $T$, $T\subset\Ker\mu_0$, and
 $T'$ is non-empty, $\mu_0\in T'$.

 We further observe that if $A$ is the empty algebraic variety, then $A'$ is
 defined to equal $T_1$, and if $A=\HOM$, then $A'=T(G)$.

 We consider separately the case when $G$ has a one-element subalgebra
 with the homomorphism $\mu_0$, and $A=\{\mu_0\}$. It is obvious that then
 $A'=T_1$.

 We see that the idea of an algebraic variety, originally linked with algebraic
 geometry, can be carried over to arbitrary varieties of algebras.

 Let us mention the following obvious relationships. Suppose that $\alpha$ runs
 over some set $I$. Then
    \begin{enumerate}
\item[1.] $(\bigcup A_{\alpha})'=\bigcap A'_{\alpha}$.
\item[2.] $(\bigcup T_{\alpha})'=\bigcap T'_{\alpha}$.
\item[3.] $\bigcup T'_{\alpha}\subset (\bigcap T_{\alpha})'$.
\item[4.] $\bigcup A'_{\alpha}\subset (\bigcap A_{\alpha})'$.
    \end{enumerate}
\noindent
 We also note, finally, that all constructions here are carried out with respect
 to a certain set of variables $X$. This set may be either finite or infinite.
 In the next section, of concern to us will be, among other things, the question
 what happens under changes of $X$.

\paragraph{2. Hilbert's Nullstellensatz.}
\stepcounter{paragraph}
 First of all, we comment on the structure of the "general solution" of a
 equation system $T$. We look for solutions in some algebra $H\in\Theta$, and
 assume that a surjective homomorphism $\mu_0\colon W\to G$ with kernel $T$ is
 given. In particular, it can be the natural homomorphism $\mu_0\colon W\to
 W/T$.

 We consider the set $\Hom (G,H)$, and let $\mu_0\Hom (G,H)$ be the set of all
 products $\mu_0\nu$, $\nu\in\Hom (G,H)$. Surjectivity of $\mu_0$ implies that
 $\mu_0\nu_1=\mu_0\nu_2$ if and only if $\nu_1=\nu_2$. Clearly, $\mu_0\Hom
 (G,H)$ is a subset of $\Hom (W,H)$.

 $\bf 2.1.$ Proposition.
For any $T$,
$$T'_H=T'=\mu_0\Hom (G,H).$$

\begin{proof}
We use the commutative diagram

\begin{center}
\unitlength=0.8mm
\linethickness{0.4pt}
\begin{picture}(50,39)(20,45)  
\put(29.67,79.67){\vector(1,0){30.00}}
\put(21.33,80.00){\makebox(0,0)[cc]{$W$}}
\put(69.0,79.67){\makebox(0,0)[cc]{$H$}}
\put(25.00,75.00){\vector(3,-4){18.00}}
\put(47.00,51.00){\vector(3,4){18.}}
\put(44.67,46.50){\makebox(0,0)[cc]{$G$}}
\put(44.67,82.67){\makebox(0,0)[cc]{$\mu$}}
\put(29.33,63.67){\makebox(0,0)[cc]{$\mu_0$}}
\put(61.00,63.67){\makebox(0,0)[cc]{$\nu$}}
\end{picture}
\end{center}

\noindent
 with $\nu$ uniquely determined by $\mu\in T'= A$. By the condition,
 $T=\Ker\mu_0\subset\Ker\mu$, and this implies that such $\nu$ ever exists.
 Therefore $T'\subset\mu_0\Hom (G,H)$.

 We now take any $\mu=\mu_0\nu$ and assume that $wT{w'}$. Then
 $w^{\mu_0}={w'}^{\mu_0}$ and $w^{\mu}={w'}^{\mu}$, $(w,{w'})\in\Ker\mu$. Hence,
 $T\subset\Ker\mu$, $\mu\in T'=A$. This gives the converse inclusion
 $\mu_0\Hom  (G,H)\subset T'$.
\end{proof}

 Hence, the general solution of the equation system $T$, where $T$ is a
 congruence, can be presented as follows:
$$A=T'=\mu_0\Hom (W/T,H).$$
 Let, furthermore, $G$ and $H$ be two algebras from $\Theta$. We consider $\Hom
 (G,H)$ and set
$$(H\defis\Ker)(G)=\bigcap\limits_{\nu}\Ker\nu,$$
 where the intersection is over all $\nu\colon G\to H$. So $(H\defis\Ker)(G)$ is
 a congruence on $G$ depending on $H$.

 Assume again that we are given a surjective homomorphism $\mu_0\colon W\to G$
 with kernel $T$.  Then the following theorem holds.

 $\bf 2.2.$ Theorem.
 Let $\mu_0^{-1}$ means "the inverse image under $\mu_0$". Then
$$T''_H=\mu_0^{-1}(H\defis\Ker)(G).$$

 In particular, the next theorem can be regarded as general Hilbert theorem on
 zeros.

 $\bf 2.3.$ Theorem.
For every congruence $T$ on $W$,
$$
T''_H=\mu_0^{-1}(H\defis\Ker)(W/T).
$$

\noindent {\bf Proof of Theorem 2.2.} Let $\tau =(H\defis\Ker)(G)$. We consider the composition
 homomorphism
$$W\stackrel{\mu_0}{\to}G\stackrel{\mu_1}{\to}G/\tau,$$
 where $\mu_1$ is the natural homomorphism,
and set $\widetilde T=\Ker\mu_0\mu_1$. Then $w\widetilde T{w'}$ means that
     $ w^{\mu_0}\Ker\mu_1{w'}^{\mu_0} $,
i.e.     $ w^{\mu_0}\tau {w'}^{\mu_0}$.
So $\widetilde  T=\mu_0^{-1}(\tau)$.
We shall verify that $\widetilde T=T''_H$.

 Assume that $w\widetilde T{w'}$. By the definition of the congruence $\tau$,
 $(w^{\mu_0},{w'}^{\mu_0})\in\Ker\nu$ and $w^{\mu_0\nu}={w'}^{\mu_0\nu}$
 for every $\nu\colon G\to H$. By Proposition 2.1, $\mu_0\nu$ is an element of
 $T'_H=A$,
 and $(w,{w'})\in\Ker\mu_0\nu$. Therefore, $(w,{w'})\in T''_H$, and we have make sure
 that $\widetilde T\subset T''_H$.

 Now assume that $wT''_H {w'}$. Then $w^{\mu_0\nu}={w'}^{\mu_0\nu}$ for every
 $\nu\colon G\to H$, and
 $(w^{\mu_0},{w'}^{\mu_0})\in\bigcap_{\nu}\Ker\nu=\tau$. This implies that
 $w^{\mu_0\mu_1}={w'}^{\mu_0\mu_1}$ and $w\widetilde T{w'}$. So
 $T''_H\subset\widetilde T$.

 Thus, $\widetilde T=T''_H$, and Theorem 2.2, as well as Theorem 2.3, are
 proved.

\medskip
 We now shall derive the classical Hilbert theorem from the general theorem 2.3.
 Two general facts of commutative algebra will be used; they actually are
 related to the Hilbert theorem.

 The first one says that if $G$ is a finitely generated associative and
 commutative algebra, then its Jacobson radical $\Rad G$ is, at the same time,
 the
 null-radical that coincides with the set of nilpotent elements of $G$.

 The other fact that we need consists in the following: if $T$ is a proper
 ideal of the ring $R=P[x_1,\ldots,x_n]$, and if $K$ is an algebraically closed
 extension of the field $P$, then there is a homomorphism $\mu\colon R\to K$ for
 which $T\subset\Ker\mu$. A property like this could serve as a general
 definition of algebraic closeness of arbitrary universal algebras.

 We now check that, under these conditions, the following equality holds:
$$(K\defis\Ker)(R/T)=\Rad (R/T).$$

 The radical on the right is the intersection of the maximal ideals.
 Suppose that $T_0/T$ is
 a maximal ideal of $R/T$. Then $T_0$ is a maximal ideal of $R$, and there is a
 homomorphism $\mu\colon R\to K$ for which $T_0\subset\Ker\mu$. It follows from
 the maximality condition that $T_0=\Ker\mu$. Since $T\subset\Ker\mu$, the
 homomorphism $\mu$ induces another homomorphism $\nu\colon R/T\to K$, and here
 $T_0/T=\Ker\nu$. Therefore, every maximal ideal of $R/T$ is realized as the
 kernel of some $\nu$. This means that the inclusion
$$\Rad (R/T)\supset (K\defis\Ker)(R/T)$$
 holds. Every element of $\Rad (R/T)$ is nilpotent, and every nilpotent element
 of $R/T$ belongs to the kernel of any $\nu\colon R/T\to K$. Hence
the converse
 conclusion.

 The Hilbert theorem now is an obvious consequence of the equality just proved
 and Theorem 2.3.

 Other applications of Theorem 2.3 will be discussed in what follows. In any
 particular case all reduces to calculating the corresponding
 $(H\defis\Ker)(W/T)$.

 We note, furthermore, that triviality of the kernel $(H\defis\Ker)(G)$ means
 that the algebra $G$ has a full system of representations in $H$. This, in its
 turn, means that the congruence $T$ in $W$ is closed, $T=T''$, if and only if
 the algebra $W/T$ has a full system of representations in $H$.

 Let us give one more variant of Hilbert theorem.

 $\bf 2.4.$ Theorem.
For any congruence $T$ in a free algebra $W$ and any algebra $G \in
\Theta$
the corresponding closure $T''_G$ is the intersection of all congruences
$\tau$ in $W$, containing $T$ and such that there is an injection
$W / \tau \to G$.

In particular, if $\Theta$ is a variety of all groups and $G=F$ is a
free group, then the group $W / T''_G$ is approximated by free groups.

\paragraph*{3. Verbal varieties.}
\stepcounter{paragraph}
 We shall consider the particular case when $T$ is a fully invariant, or
 verbal, congruence on $W=W(X)$. For every algebra $G\in\Theta$, we call the
 respective algebraic variety $T'$ a verbal variety. We are interested in $T''$
 in this case.

 The congruence $T$ determines a variety of algebras $\Theta_T$, which is a
 subvariety of $\Theta$. Given an algebra $G$, we consider all its subalgebras
 $H$ in $\Theta_T$. For these, equations from $T$ become
 identities, $T'_H=\Hom (W,H)$.

 We denote the system of all the subalgebras by $\Theta_T(G)$. This object is,
 to a certain extent, the dual of the verbal congruence on $G$ relatively to
 $\Theta_T$. We also denote by $\widetilde T(G)$ the congruence composed of
 the identities of  the class $\Theta_T(G)$ in the free algebra $W$; we shall
 call it the {\it congruence of identities\/} (of $\Theta_T(G)$).
 For any $T$, the congruence  $\widetilde T(G)$ is verbal.

 $\bf 3.1.$ Theorem.
If $T$ is a verbal congruence, then $T''_G=\widetilde T(G)$. Specifically,
$T''_G$ is also a verbal congruence.

\begin{proof}
 Let us compute the kernel
$$(G\defis\Ker)(W/T).$$
 First of all, we observe that the algebra $W/T$ is free in $\Theta_T$ over $X$.

 Furthermore, we make the following general remark. Let $\Theta$ be a variety of
 $\Omega$-algebras, $W(X)$ be the free algebra in $\Theta $ over $X$, and $G$ be
 an $\Omega$-algebra not necessary from $\Theta$. We shall examine the kernel
 $(G\defis\Ker)(W(X))$.

 To this end, we select in $G$ all the subalgebras $H$ belonging to $\Theta$,
 and  denote the system of these subalgebras by $\Theta(G)$. Let $T(G)$ be the
 congruence of all identities of $\Theta(G)$ in $W(X)$. Then
 $(G\defis\Ker)(W(X))=T(G)$.

 Let us demonstrate this. We denote by $M$ the system of homomorphisms
 $\nu\colon W(X)\to H$, where $H$ is a subalgebra of $G$ from $\Theta(G)$. The
 intersection of all kernels $\Ker\nu$ over all $\nu\in M$ is $T(G)$. Every
 $\nu\in M$ is at the same time a homomorphism $\nu\colon W(X)\to G$. On the
 other  hand, $\nu\colon W(X)\to G$, where $\im\nu=H$ is a subalgebra in
 $\Theta(G)$, and we also have $\nu\colon W(X)\to H$. Therefore, we may identify
 the sets $M$ and $\Hom (W(X),G)$. This leads to the needed equality.

 We apply it in the situation with $\Theta_T$  taken for $\Theta$. Then the
 kernel
$$(G\defis\Ker)(W/T)$$
 is the congruence of identities of the system $\Theta_T(G)$ in $W/T$. But
 in this case the full inverse image $T''_G=\mu_0^{-1}(G\defis\Ker)(W/T)$ is the
 congruence of the identities of $\Theta_T(G)$ in $W(X)$. It follows that
 $T''_G=\widetilde T(G)$.
\end{proof}

 There is an example.

 Suppose that the initial variety $\Theta$ is the variety of groups, the set $X$
 is infinite, $F(X)$ is the corresponding free group, and $F^1(X)=T$ is
 its commutant. Take a group $G$ and consider two cases: $G$ is of finite
 exponent and the exponent of $G$ is infinite. The commutant $T$ determines the
 variety of commutative groups, and the same variety is generated, in the second
 case, by the commutative subgroups of $G$. For this reason the commutant
  is closed in the second case, $T=T''$. In the first case the commutative subgroups
 of $G$ generate the variety of commutative groups of exponent $n$.
 Consequently, $T''$ is a verbal congruence generated by the commutant and the
 element $x^n$.

 We make one more useful remark concerning verbal varieties.

 $\bf 3.2.$ Proposition.
Suppose that $T$ is a verbal congruence on $W(X)$ and that $G$ is an
 algebra from $\Theta$. Then
$$A=T'_G=\bigcup\limits_H\Hom (W,H),$$
 where the union is taken over all $H\in\Theta_T(G)$.

\begin{proof}
 Let $\mu\in A$. Then $T\subset\Ker\mu$. As $W/T\in\Theta_T$, also
 $W/\Ker\mu\in\Theta_T$. But then $H=\im\mu\in\Theta_T(G)$, $\mu\in\Hom (W,H)$.

\smallskip
\noindent Notice that here, and below, if $H$ is a subalgebra of $G$, then $\Hom (W,H)$ is
 treated as $\HOM$. Moreover, $T'_H=T'_G\cap\Hom (W,H)$ for every $T$.

 Now assume that $H\in\Theta_T(G)$. This means that $T'_H=\Hom (W,H)$, and then
 $\Hom (W,H)$ is included into $T'_G=A$.
 \end{proof}

 Theorem 3.1 also is an easy consequence of the above remark.

\paragraph{4. Geometric equivalence of algebras.}
\stepcounter{paragraph}
\hbox{}\hskip 5pt

 $\bf 4.1.$ Definition.
 Algebras $G_1$ and $G_2$ from $\Theta$ are said to be geometrically equivalent
 if
$$T''_{G_1}=T''_{G_2}$$
 for every $T$ from $W(X)$.

It is easily understood that this condition is equivalent to the following one:
 any congruence $T$ on $W(X)$ is $G_1$-closed if and only if it is $G_2$-closed:
 $T''_{G_1}=T$ iff $T''_{G_2}=T$.

 Indeed, the equivalence of $G_1$ and $G_2$ implies the latter condition.
 Assume, on the other hand, that this condition is fulfilled for every $T$. Then
 $T\subset T''_{G_1}$ and, furthermore, $T''_{G_2}\subset T''_{G_1}$. The
 converse inclusion is proved in the same way.

 For every $G\in\Theta$, let $M_G$ stands for the system of all algebraic
 varieties in $\HOM$. By $\Cl_G(W)$ we denote the system of all $G$-closed
 congruences on $W=W(X)$. There is a natural bijection between the sets $M_G$
 and $\Cl_G(W)$.

 Equivalence of the algebras $G_1$ and $G_2$ means that the sets $\Cl_{G_1}(W)$
 and $\Cl_{G_2}(W)$ coincide. Now it is clear that the equivalence determines a
 canonic bijection $M_{G_1}\to M_{G_2}$. We once more stress that the
 definition of equivalence is related to specific $X$. Therefore, the question
 is, in fact, on $X$-equivalence. Clearly, if $G_1$ and $G_2$ are isomorphic,
 then they are equivalent relatively to every $X$.

 The main problem is to learn to recognize equivalence of algebras by their
 properties. For example, can two groups be equivalent if one of them is
 commutative while the other is not? In the classical geometry two algebraically
 closed fields $K_1$ and $K_2$, if they both are extensions of the same $P$, are
 equivalent. As we have already mentioned, the corresponding problem for fields
 that are not algebraically closed, is still open. We can only note, for
 example,
 that every $K$ is equivalent to every of its ultarpowers (cf. \S 5).

 It follows from Theorem 2.3 that algebras $G_1$ and $G_2$ are
equivalent if and
 only if
$$(G_1\defis\Ker)(W/T)=(G_2\defis\Ker)(W/T)$$
 for every $T$.

 Let us take, for example, the variety of vector spaces over a
 given field $P$ for $\Theta$. If $G$ and $H$ are vector spaces in $\Theta$,
 then
 $(H\defis\Ker)(G)=0$, whence, in this case, all $T$ are closed for every $G$,
 and any two spaces are equivalent.

 $\bf 4.2.$ Problem.
Assume that $K$ is a commutative ring with unit and that $\Theta$ is the
 variety of $K$-modules. Under what conditions are two modules $G_1$ and $G_2$
 equivalent?

 Of course, the problem has to be considered for several particular $K$, e.g.
 $K=\Z$. What two commutative groups are equivalent? The general problem also
 depends on the choice of the set of variables $X$.

 Let again the variety $\Theta$ be arbitrary, and let $X$ be fixed.

 $\bf 4.3.$ Theorem.
If algebras $G_1$ and $G_2$ are $X$-equivalent, then they have the same
 identities on $W(X)$.

\begin{proof}
 We shall apply Theorem 3.1. Let $T=T(G_1)$ be the congruence of all identities
 of $G_1$ on $W=W(X)$. Then $T'_{G_1}=\Hom (W,G_1)$, and $T''_{G_1}=T(G_1)=T$.
 Now take  $T''_{G_2}=\widetilde T(G_2)$, and let $\Theta_T(G_2)$  be the system
 of all subalgebras $H$ of $G$  belonging to the variety $\Theta_T$.
  $\widetilde T(G_2)$ is the congruence of all identities of $\Theta_T(G_2)$.
  If $G_1$ and $G_2$ are equivalent, then
 $T=\widetilde T(G_2)$. We have $\widetilde T(G_2)\supset T(G_1)$, whence
 $T(G_1)\supset T(G_2)$. Likewise, $T(G_1)\subset T(G_2)$ and, consequently,
 $T(G_1)=T(G_2)$.
\end{proof}

 For $X$ infinite, the equality $T(G_1)=T(G_2)$ implies that $\Var (G_1)=\Var
 (G_2)$; this means that the algebras $T(G_1)$ and $T(G_2)$ have the same
 equational theory and the same equational logic.

 We shall mention some consequences for $X$ infinite.

 First of all, we observe that if $G_1$ and $G_2$ are finite simple groups, then
 $\Var (G_1)=\Var (G_2)$ only if $G_1$ and $G_2$ are isomorphic. We therefore
 conclude that two finite simple groups are equivalent if and only if they are
 isomorphic.

 Moreover, we now can say that a commutative group is equivalent to no
 non-commutative group.

 It is naturally, in the case $X$ is infinite, to wonder whether equivalence of
 algebras $G_1$ and $G_2$ implies that they have the same universal theory.
 The example with vector spaces demonstrates that it is in general not so. It is
 not so also under arbitrary $\Theta$. This follows form the proposition below.

 We take any $\Theta$ and also fix arbitrary $X$.

 $\bf 4.4.$ Proposition.
For every algebra $G\in\Theta$ and every set $I$, the algebras $G$ and
 $G^I$ are equivalent.

\begin{proof}
 We confine ourselves to one-sorted algebras.

 Let us take any algebra $A\in\Theta$, and show that the equality
$$(G\defis\Ker)(A)=(G^I\defis\Ker)(A)$$
 ever holds. We denote its left-hand part by $\tau_1$ and the right-hand one --
 by $\tau_2$. Suppose that $a\tau_1{a'}$ for $a$ and ${a'}$ from $A$. This means
 that, for every $\nu\colon A\to G$, $a^{\nu}={a'}^{\nu}$.

 We need to verify that $a\tau_2{a'}$. To this end, take arbitrary $\mu\colon A\to
 G^I$. We then have to prove that $a^{\mu}={a'}^{\mu}$, i.e. that
 $a^{\mu}(\alpha)={a'}^{\mu}(\alpha)$ for every $\alpha\in I$.

 Let us consider projections $\pi_{\alpha}\colon G^I\to G$. Then
 $\mu\pi_{\alpha}=\nu_{\alpha}\colon A\to G$. Moreover,
 $a^{\mu}(\alpha)={a'}^{\mu\pi_{\alpha}}$. What remains to show is that
 $a^{\mu\pi_{\alpha}}={a'}^{\mu\pi_{\alpha}}$ or, equivalently , that
 $a^{\nu_{\alpha}}={a'}^{\nu_{\alpha}}$ for every $\alpha \in I$. Clearly,
 this is so when $a\tau_1{a'}$.  Therefore, $a\tau_1{a'}$ implies  $a\tau_2{a'}$.

 Conversely, suppose that $a\tau_2{a'}$. Given $\nu\colon A\to G$, we construct
 $\mu\colon A\to G^I$ by setting $a^{\mu}(\alpha)=a^{\nu}$ for every $\alpha$.
 For all $a\in A$, the elements $a^{\mu}$ are constants, and $a^{\mu}={a'}^{\mu}$.
 If $\alpha\in I$, then $a^{\mu}(\alpha)=a^{\nu}={a'}^{\mu}(\alpha)={a'}^{\nu}$. So
 $a\tau_1{a'}$, and $a\tau_2{a'}$ implies $a\tau_1{a'}$. We have arrived at
 $\tau_1=\tau_2$.
\end{proof}

 Now, we take a congruence $T$ on $W$ and let $A$ be  the algebra $W/T$. Then
$$(G\defis\Ker)(W/T)=(G^I\defis\Ker)(W/T).$$
 This means that the algebras $G$ and $G^I$ are equivalent.

 Generally, equivalent algebras $G$ and $G^I$ may have different
 pseudoidentities. Much more, distinct are the universal theories of $G$ and
 $G^I$.

 The concept of algebra equivalence can also be defined on the universal logic
 level. As we shall see, in this case the equivalence of algebras implies that
 they have the same universal theories.

 In conclusion of the subsection we note that if two algebras $G_1$ and $G_2$
 are $X$-equivalent and if neither of them has a one-element subalgebra, then,
 for any
 proper congruence $T$ on $W(X)$, the varieties
$A=T'_{G_1}$ and $B=T'_{G_2}$
 are either both empty or both
 nonempty.

 Indeed, assume that $G_1$ and $G_2$ are equivalent and that $A$ is empty. Then
 $T''_{G_1}=A'$ is the unit congruence $T_1$. But then $T_1=T''_{G_2}=B'$ and
 $B=B''=T'_1$ are empty varieties.

 The converse does not hold: if for any proper congruence $T$ on $W(X)$ the
 varieties $A=T'_{G_1}$ and $B=T'_{G_2}$ are both empty, or both nonempty, then
 this does not mean that $G_1$ and $G_2$ are equivalent.

\paragraph*{5. Generalized equations.}
\stepcounter{paragraph}
 In this subsection, an algebra $G$ from the variety $\Theta$ is assumed to be
 fixed, and constants occurring in the equations considered are supposed to
 belong $G$.

 Such equations are connected with the passage to another variety $\Theta'=\Theta^G$
 depending on $G$. We begin with defining the category $\Theta'$. Its objects are
 pairs $(H,h)$ where $H$ is an algebra from $\Theta$ and $h\colon G\to H$ is a
 homomorphism in $\Theta$. The representation--homomorphism $h$ makes elements of $G$ constants in $H$. The pairs $(H,h)$ are termed $G$-algebras, or
 algebras over $G$; compare with associative algebras over a field $P$.

 If $(H_1,h_1)$ and $(H_2,h_2)$ are two $G$-algebras, then a homomorphism
 $\alpha\colon H_1\to H_2$ is a $G$-algebra (homo)morphism in case the diagram

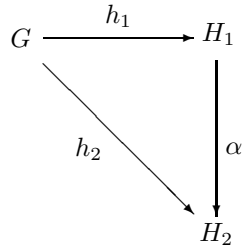
\begin{center}
\unitlength=1.00mm
\linethickness{0.4pt}
\begin{picture}(55,32)(25,6)    
\put(30.00,35.00){\makebox(0,0)[cc]{$G$}}
\put(33.00,35.00){\vector(1,0){20.00}}
\put(56.33,35.60){\makebox(0,0)[cc]{$H_1$}}
\put(33.00,31.67){\vector(1,-1){20.00}}
\put(56.00,9.00){\makebox(0,0)[cc]{$H_2$}}
\put(43.00,38.00){\makebox(0,0)[cc]{$h_1$}}
\put(58.33,20.33){\makebox(0,0)[cc]{$\alpha$}}
\put(39.00,20.33){\makebox(0,0)[cc]{$h_2$}}
\put(56.00,32.00){\vector(0,-1){20.67}}
\end{picture}\hspace{-3em}
\end{center}

\noindent
commutes.

 The category of $G$-algebras can be presented as a variety if one specifies
 the algebra $G$ by generators and relations. Generators are nullary operations
 added to the collection of primitive operations $\Omega$, while defining
 relations are added to the collection of identities specifying $\Theta$; cf.
 \cite{BPl3} and \S 2 here. This way we obtain a new variety which, generally,
 depends on the particular presentation of $G$ by generators and relations.
 However, all these varieties are equivalent as categories, and all the
 categories are, in turn, equivalent to the category of $G$-algebras. In what
 follows, the category $\Theta'$ is regarded as a variety.

 If $X$ is a set, possibly, many-sorted, and $W(X)=W$ is the free algebra in
 $\Theta$ over $X$, then the free over the same $X$ algebra $W'$ in $\Theta'$
 can be presented as the free (in $\Theta$) product $G*W$ with the homomorphism
 $h_0\colon  G\to G*W$ determined by the corresponding projection.

 Where $(H,h)$ is a $G$-algebra, every mapping $\mu\colon X\to H$ induces a
 homomorphism $\mu\colon W\to H$. Together with $h\colon G\to H$, it gives
 $\mu\colon G*W\to H$. This latter $\mu$ is a $G$-algebra homomorphism;
 commutativity of the diagram

\begin{center}
\unitlength=1.00mm
\linethickness{0.4pt}
\begin{picture}(55,32)(25,6)
\put(30.00,35.00){\makebox(0,0)[cc]{$G$}}
\put(33.00,35.00){\vector(1,0){20.00}}
\put(61.33,35.60){\makebox(0,0)[cc]{$G*W$}}
\put(33.00,31.67){\vector(1,-1){20.00}}
\put(56.00,9.00){\makebox(0,0)[cc]{$H$}}
\put(43.00,38.00){\makebox(0,0)[cc]{$h_0$}}
\put(62.33,20.33){\makebox(0,0)[cc]{$\mu$}}
\put(39.50,20.33){\makebox(0,0)[cc]{$h$}}
\put(58.00,32.00){\vector(0,-1){20.67}}
\end{picture}\hspace{-3em}
\end{center}

\noindent
follows from the definition of a free product.

 We still note that the algebra $G$ can be regarded as a $G$-algebra if we
 proceed from the identity homomorphism $G\to G$.

 Now a generalized equation has the form $w={w'}$, where $w$ and ${w'}$
are elements
 of $W'$ of the same sort. The coefficients of the equation are also from $G$.
 Such equations are resolved in $G$-algebras $H$, in particular, in the
 $G$-algebra $G$. Clearly, all the constructions considered above are applicable
 to these equations. The initial variety here is $\Theta'$.

 Of special interest is the situation when the homomorphism $h\colon G\to H$ in
 a $G$-algebra $H$ is injective. We shall prove, in this context, the following
 well-known result.

 Let us agree to say that a $G$-algebra $(H,h)$ is faithful if $h\colon G\to H$
 is an injection.

 $\bf 5.1.$ Proposition.
Suppose that $T$ is a congruence on the $G$-algebra $G*W$, and let
 $\alpha\colon G*W\to G*W/T$ be the natural homomorphism. The system of
 equations $T$ has a solution in a faithful $G$-algebra $H$ if and only if the
 homomorphism $h_0\alpha\colon G\to (G*W)/T$ is an injection.

\begin{proof}
 Assume that $h=h_0\alpha$ is an injection. Then we take $H=(G*W)/T$ and
 consider the $G$-algebra $(H,h)$, which is faithful. The point
 $\mu=\alpha\colon G*W\to H$ is a solution of the system of equations $T$.

 Now suppose that $\mu\colon G*W\to H$ is a solution of the system $T$ in a
 faithful $G$-algebra $H$ with an injection $h\colon G\to H$. By the definition
 of a homomorphism in $\Theta'$, we have a commutative diagram

\begin{center}
\unitlength=0.80mm
\linethickness{0.4pt}
\begin{picture}(50,32)(20,51)
\put(29.67,79.67){\vector(1,0){30.00}}
\put(21.33,80.00){\makebox(0,0)[cc]{$G*W$}}
\put(65,80){\makebox(0,0)[cc]{$H$}}
\put(44,51.50){\makebox(0,0)[cc]{$G$}}
\put(44.67,82.67){\makebox(0,0)[cc]{$\mu$}}
\put(48.00,55.00){\vector(3,4){15.00}}
\put(59.00,63.67){\makebox(0,0)[cc]{$h$}}
\put(41.00,55.00){\vector(-3,4){15.00}}
\put(29.67,63.33){\makebox(0,0)[cc]{$h_0$}}
\end{picture}
\end{center}

\noindent
 Since $\mu$ is a solution of $T$, the inclusion $T\subset\Ker\mu$ holds, and
 this gives one more commutative diagram

\begin{center}
\unitlength=0.80mm
\linethickness{0.4pt}
\begin{picture}(50,32)(20,48)
\put(29.67,79.67){\vector(1,0){30.00}}
\put(21.33,80.00){\makebox(0,0)[cc]{$G*W$}}
\put(65.55,80.00){\makebox(0,0)[cc]{$H$}}
\put(45.00,50.00){\makebox(0,0)[cc]{$(G*W)/T$}}
\put(44.67,82.67){\makebox(0,0)[cc]{$\mu$}}
\put(24.33,76.00){\vector(3,-4){16.00}}
\put(28.67,62.33){\makebox(0,0)[cc]{$\alpha$}}
\put(46.00,54.67){\vector(3,4){16.00}}  
\put(59.00,62.33){\makebox(0,0)[cc]{$\beta$}}
\end{picture}
\end{center}

\noindent
 Comparing the diagrams, we observe that $h_0\alpha\beta=h$. Therefore, since
 $h$ is an injection, so is $h_0\alpha$.

 As to the variety $\Theta'$, we note that if $G$ is not a one-element algebra,
 then no faithful $G$-algebra has a one-element subalgebra. For this reason,
 some system $T$ my fail to have a solution in any such an algebra at all, and
 the
 corresponding $T'$ may be empty. If $\Theta$ is the variety of all groups, then
 here all $T'$ are nonempty, but the situation changes for $G$-groups.

 Along with generalized equations, generalized identities can be considered. The
 literature on generalized equations and identities is quite extensive
 \cite{RS,Neu1,Neu2,Le,Ly,Pri}.

 We still shall make some remarks on the closure of a point. We shall make it
 apparent that if equations admit solutions in a $G$-algebra $G$, then every
 point $\mu\colon G*W\to G$ coincides with its closure.

 We proceed from the projections $h_0\colon G\to G*W$ and $h_1\colon W\to G*W$.
 For every point $\nu\colon G*W\to G$, we also have $h_0\nu =\varepsilon\colon
 G\to G$.

 For each $x\in X$, we take $x^{h_1}$ and $x^{h_1\mu h_0}$. These elements both
 belong to $G*W$ so that the equation $x^{h_1}\equiv x^{h_1\mu h_0}$ makes
 sense. We denote the system of such equations for all $x\in X$ by $T$. The
 equality $x^{h_1\mu h_0\nu}=x^{h_1\mu}$ holds for every $\nu$;
 so if $\nu=\mu$,
 then $x^{h_1\mu}=x^{h_1\mu h_0\mu}$. Therefore, $T\subset\Ker\mu$. If
 $\Ker\mu\subset\Ker\nu$, then it follows that $x^{h_1\nu}=x^{h_1\mu
 h_0\nu}=x^{h_1\mu}$. This means that $\nu$ and $\mu$ agree on $W$. Moreover,
 $g^{h_0\mu}=g=g^{h_0\nu}$ for every $g\in G$. Consequently, $\mu$ and $\nu$
 agree on $G$. But then $\mu=\nu$. Thus, the closure of the point $\mu$ only
 consists on $\mu$ itself.

We return to the situation of $G$-algebras over given $\Theta$. Let
$T$ be a congruence in a free algebra $G \ast W$. The solution is
considered in $G$-algebra $G$.
\end{proof}

 $\bf 5.2.$ Proposition.
A congruence $T''_G$ is the intersection of all congruences $\tau$ in
$G\ast W$, containing $T$, such that $G\ast W / \tau$ and $G$ are
isomorphic $G$-algebras.

The Proposition follows from 2.4 and from that $G$-algebra $G$ has no
proper subalgebras.

\paragraph{6. Algebra and topology in connection with varieties.}
\stepcounter{paragraph}
We assume in this subsection that $\Theta$ and $X$ are fixed, take $W=W(X)$,
 and consider algebraic varieties in the "affine space" $\HOM$, $G\in\Theta$.

 Suppose that $T$ is a binary relation on $W(X)$ and $A=T'$ is the corresponding
 algebraic variety. We couple with $A$ the algebra $W/A'=W/T''$ in $\Theta$.
 By the
 definition, $A'=T''$ is the intersection of all $\Ker\mu$, $\mu\in A$. So
 $W/A'$ is approximated by algebras $W/\Ker\mu$ or, what is the same, by
 algebras $\im\mu$, $\mu\in A$, which are subalgebras of $G$.

 We now shall look at the algebra $W/A'$ from another point of view.

 We shall deal with mappings (functions) $\alpha\colon A\to G$, where $A$ is
 an algebraic
 variety in $\HOM$, $G\in\Theta$. Such a function is said to be regular if there
 is $w\in W$ satisfying $\alpha(\mu)=w^{\mu}$ for every point $\mu$. If
 many-sorted algebras are considered, then the element $w$ has a definite sort,
 and the question is on a regular function of this sort.

 The function $\alpha$ can also be given by another element ${w'}$ of the same
 sort. Then $w^{\mu}={w'}^{\mu}$ for every $\mu\in A$. This means that $(w,{w'})\in
 A'$.

 For each $\mu\in A$, $A'\subset\Ker\mu$; consequently, we also have a
 homomorphism $\mu\colon W/A'\to G$. Given $w\in W$, we denote by $\bar w$ the
 corresponding element of $W/A'$. Then $\alpha(\mu)=\bar w^{\mu}$. Therefore,
 elements of $W/A'$ can be regarded as regular functions of kind $A\to G$.

 The algebra $W/A'$ itself is termed the algebra of regular functions on the
 variety $A$ with values in $G$. All this is in agreement with the presentation
 of $W/A'$ as a subdirect product of the algebras $W/\Ker\mu$, $\mu\in A$. The
 algebra $W/A'$ is also called the co-ordinate algebra of the variety $\Theta$.

 We call, furthermore, an algebra $H\in\Theta$ $G$-exact, $G\in\Theta$, if
 $(G\defis\Ker)(H)$ is the trivial congruence on $H$. Let $\mu_0\colon W(X)\to
 H$ be a surjective homomorphism for such $H$ with $\Ker \mu_0 = T$. Then,
 obviously, the algebra $H$
 can be presented as the co-ordinate algebra of the variety $A=T'$ in $\Hom
 (W(X),G)$.

 On the other hand, every algebra $W/A'$ is $G$-exact.

 The algebra $W/A'$ is an important invariant of the variety $A$. Below, we
 shall introduce, for $G$ fixed, the notion of isomorphism of two varieties. It
 will be proved that varieties $A$ and $B$ are isomorphic if and only if
 isomorphic are respective algebras $W/A'$ and $W/B'$. This is a
 generalization of a classical theorem.

 Varieties can be classified from the viewpoint of properties of their algebras,
 e.g. according to identities of the algebras.

 In particular, varieties $A$ and $B$ could be called similar if
 $$\Var(W/A')=\Var (W/B').$$

 Now assume that $A$ is a verbal variety specified by a verbal congruence $T$.
 Then $A'=T''$ is also a verbal congruence, and $W/A'$ is a free over $X$
 algebra of the corresponding subvariety of $\Theta$.

 We now return to the question on equivalence of two algebras $G_1$ and $G_2$
 from $\Theta$. If $T$ is a congruence on $W$, then we also have $A=T'_{G_1}$
 and $B=T'_{G_2}$. If $G_1$ and $G_2$ are equivalent, then $A'$ and $B'$
 coincide, and so do also the algebras $W/A'$ and $W/B'$. This is one more
 argument in favour of the notion of equivalence we discuss.

 Regretfully, we cannot speak of isomorphism of the varieties $A$ and $B$, for
 they are related to distinct $G_1$ and $G_2$. Nevertheless, there must be
 something
 in common for $A$ and $B$; at least, in the situation of the classical
 geometry. Cf. also \S 4.

 Now we shall comment on the topology on $\HOM$ connected with algebraic
 varieties. We already have noticed that the sum of two algebraic varieties
 need not
 be an algebraic variety. Because of this, in order to obtain a topology on
 $\HOM$ we regard algebraic varieties and finite unions of them to be the
 closed sets.

 These sets are described by systems of pseudoequations. For more detail, see
 \S 5.\!\! The corresponding topology is considered as a Zariski topology on $\HOM$.
 It is not clear to us what it can offer in the general situation under
 consideration. However, one useful consideration concerning the closure of a
 point can be made. We already discussed a particular situation of this sort;
 now the overall picture will be sketched.

 Let us take a point $\mu\in\HOM$. The closure of $\{\mu\}$ is $\{\mu\}''$. As
 $\{\mu\}'=T$ is $\Ker\mu$, we conclude that
    $\{\mu\}''=\{\nu,\  T=\Ker\mu\subset\Ker\nu\}$.
   The set $\HOM$ can be equipped with a
 pseudo-ordering relation by setting $\mu\le\nu$ if $\Ker\mu\subset\Ker\nu$.
 Then the closure of the point $\mu$ is the set of the points $\nu$ with
 $\mu\le\nu$. Points $\mu$ and $\nu$ are equivalent if $\mu\le\nu$ and
 $\nu\le\mu$, i.e. if $\Ker\mu=\Ker\nu$.

 Clearly, equivalence of $\mu$ and $\nu$ also means the closures of these points
 coincide. All this is well-known in the classical situation.

 As we know, two semigroups -- $\End W$ and $\End G$ -- are acting on $\HOM$.
 Let us see how the actions conform with algebraic varieties.

 First of all, recall that an action of a semigroup is defined by way of
 multiplying morphisms. If $s\in S=\End W$ and $\mu\in\HOM$, then $\mu s$ is given
 by the rule $(\mu s)(x)=\mu (s(x))$. If $\sigma\in\End G$, then, for
 $\sigma\mu$, $(\sigma\mu)(x)=\sigma(\mu(x))$. We here apply morphisms from the
 left. If $A$ is a subset of $\HOM$, then
$$
\mu\in sA\Longleftrightarrow \mu s\in A, \quad\mu\in
    A\sigma\Longleftrightarrow\sigma\mu\in A.
$$
 It is easily seen that the algebraic variety $T'=A$ for every $T$ on $W$ is
 always invariant under the action of $\End G$: if $\mu\in A$, then
 $\sigma\mu\in  A$.
 This way, on each $A$ the action structure of $\End G$ is defined. In
 particular, if $\sigma\in\Aut G$, then the points $\mu$ and $\sigma\mu$ are
 equivalent: they determine the same closure.

 If $T$ is a binary relation on $W$, then we define $sT$, $s\in\End W$, to be the
 new binary relation determined by the rule

    \begin{quote}
$w(sT){w'}$ \ if there are $w_1$ and ${w'}_1$ such that
\ $w_1^s=w$, ${w'}^s_1={w'}$ and $w_1T{w'}_1$.
    \end{quote}

 $\bf 6.1.$ Proposition.
{\rm 1.} For each $s\in\End W$ and $T$,
$$(sT)'=sT'.$$
{\rm 2.} For each $A\subset\HOM$ and $s\in\Aut W$,
$$(sA)'=sA'.$$

These equalities will be verified in \S 5 in a more general context.

We note two consequences of the proposition.

\smallskip
\noindent {\bf 1.} If $A=T'$ is an algebraic variety, then so is the set $sA$ for every
 $s\in\End W$. In other words, the system $M_G$ of all algebraic varieties in
 $\HOM$ is invariant with respect to the action of the semigroup $\End W$.

\smallskip
\noindent {\bf 2.} If $T=A'$ is a closed congruence on $W$, then so is the congruence $sT$ for
 every automorphism $s$ of $W$: if $T''=T$, then $(sT)''=sT$. The system
 $\Cl_G (W)$ is invariant with respect to action of the group $\Aut W$.

\smallskip
We shall see below that if $s$ is an automorphism, then it determines an
isomorphism between the varieties $A$ and $sA$.

The next proposition, which reveals connections with closures, is also related
to Proposition 6.1.

 $\bf 6.2.$ Proposition.
Assume that $s\in\Aut W$. Then, for all $T$ and $A$,
\begin{enumerate}
 \item[\rm 1.] $(sT)''=sT''$.
 \item[\rm 2.] $(sA)''=sA''$.
\end{enumerate}

\begin{proof}
 Of course, $(sT)'=sT'$. We apply ${}'$ once more:
 $(sT)''=(sT')'=sT''$. Likewise in the second case: $(sA)'=sA'$ and, further,
 $(sA)''=(sA')'=sA''$.
\end{proof}

 Finally, we comment on the lattice of varieties for $G$ fixed. In the classical
 situation all varieties make up a lattice, which is a sublattice of the
 distributive lattice of all subsets of $\HOM$.
 The lattice of varieties appears also in the general case. Then $A\cdot B=A\cap
 B$ and $A+B=(A\cup B)''=(A'\cap B')'$. What can we say about this lattice? Are
 there any connections with congruence lattice of $W$? What about two algebras
 $G_1$ and $G_2$ when the respective variety lattices are isomorphic? We have
 not examined these questions.

\paragraph*{7. Relation to the $\bf\Theta$-structure of algebras.}
\stepcounter{paragraph}
 This subsection is concerned with the subject of the preceding one. We here
 equip the set $\HOM$ with the structure of the variety $\Theta$. This can only
 be done in the case of one-sorted algebras, for the set $\HOM$ is always
 one-sorted. So let $\Theta$ be a variety of one-sorted $\Omega$-algebras.

Let $X$ be the set of variables, and assume that $G$ is an algebra from
 $\Theta$. Then $G^X$ is also an algebra in $\Theta$. The $\Theta$-structure of
 $G^X$ can be transferred to $\HOM$.

 If $\omega\in\Omega$ is an $n$-ary operation, then for homomorphisms
 $\mu_1,\ldots,\mu_n\in$\\ $\HOM$
$$x^{\mu_1\cdots\mu_n\omega}=x^{\mu_1}\cdots x^{\mu_n}\omega.$$
 But we cannot be sure that
$$w^{\mu_1\cdots\mu_n\omega}=w^{\mu_1}\cdots w^{\mu_n}\omega.$$
 for arbitrary $w\in W$. This is so only under some specific conditions, which
 are discussed below.

 Let $\omega_1$ be an $n$-ary operation and $\omega_2$--an $m$-ary operation
 from
 $\Omega$, none of them nullary, and consider a matrix $(x_{ij})$,
 $i=1,\ldots,m$, $j=1,\ldots,m$, consisting of variables. Set
$$
w_1=(x_{11}\cdots x_{1n}\omega_1)\cdots (x_{m1}\cdots
        x_{mn}\omega_1)\omega_2,
$$
$$
w_2=(x_{11}\cdots x_{m1}\omega_2)\cdots (x_{1n}\cdots
        x_{mn}\omega_2)\omega_1.
$$
 The formula $w_1\equiv w_2$ is a kind of a commutation law for the operations
 $\omega_1$ and $\omega_2$. If $\omega_1=0_{\alpha}$ is a nullary operation, and
 $\omega_2=\omega$ is arbitrary, then their commutation means that
$$0_{\alpha}\cdots 0_{\alpha}\omega=0_{\alpha}.$$
 Commutation of two nullary operations $0_{\alpha}$ and $0_{\beta}$ means that
 $0_{\alpha}=0_{\beta}$.

 The commutation law can be applied to coinciding operations, too. For
 example, in a group this law, when applied to multiplication, means that the
 group is Abelian. This need not be so in a semigroup.

 An algebra $G\in \Theta$ is said to be commutative if the commutation law
 holds in it for every two operations, including the case of equal operations.

 $\bf 7.1.$ Proposition.
If $G$ is commutative, then
$$w^{\mu_1\cdots\mu_n\omega}=w^{\mu_1}\cdots w^{\mu_n}\omega$$
 for every operation $\omega\in \Omega$, all $\mu_1,\ldots,\mu_n$, and every
 $w\in W$.

 Now assume that $G$ is commutative and that $T$ is a set of formulas of kind
 $w\equiv {w'}$.

 $\bf 7.2.$ Proposition.
The algebraic variety $A=T'$ is a subalgebra of $\HOM$.

\begin{proof}
 We have to find out whether the set $A$ is closed under the operations from
 $\Omega$. Let $\omega\in\Omega$ be an $n$-ary operation, and let
 $\mu_1,\ldots,\mu_n\in A$. We shall check that $\mu_1\cdots\mu_n\omega\in A$.
 Suppose that $w\equiv {w'}\in T$; then
$$w^{\mu_1\cdots\mu_n\omega}=w^{\mu_1}\cdots w^{\mu_n}\omega={w'}^{\mu_1}\cdots
{w'}^{\mu_n}\omega={w'}^{\mu_1\cdots\mu_n\omega},\quad \mu_1\cdots\mu_n\omega\in
A.$$
\end{proof}

 We can use all this as follows.

 Assume that $\Omega_0$ is a subset of $\Omega$. Any algebra $G\in\Theta$ can be
 regarded as an $\Omega_0$-algebra. Considered this way, it may turn out to be
 commutative, though, in general, it may be non-commutative as well. We may
 choose $\Omega_0$ in several ways. In doing so, we also can apply the above
 considerations and conclude for the respective $T$'s that $A=T'$ is an
 $\Omega_0$-closed variety.

 Let us draw some consequences of this for the classical situation. The
 operation system $\Omega$ consists here of addition, multiplication, zero, unit
 and scalars. If we take $\Omega_0$ to contain addition, scalars and the zero,
 then the corresponding algebras are vector spaces, and they are commutative.
 The related equations are of the form
$$\alpha_1 x_1+\cdots +\alpha_n x_n=0.$$
 The corresponding varieties are $\Omega_0$-algebras.

 In another case $\Omega_0$ consists of multiplication and the unit. This also
 is a commutative collection. The equations take either the form
$$x_1^{n_1}\cdots x_k^{n_k}=x_1^{m_1}\cdots x_l^{m_l}$$
 of the form $x_1^{n_1}\cdots x_k^{n_k}=1$.

 For the corresponding $T$'s, the varieties $A=T'$ are invariant
 under $\Omega_0$.

 In particular, the parabola $y^2=x$ is closed under multiplication, while the
parabola $y^2=2px$ does not posses this property. The scalars do not
commute with
 multiplication. This is true also of hyperbolas $xy=1$ and $xy=a$. The zero
 commutes with addition and multiplication, but the unit does not commute with
 addition; also the zero and the unit do not commute.

 The general theory developed here is, of course, applicable when $\Theta$ is
 the variety of $K$-modules, where $K$ is a commutative ring with unit. In this
 case, the algebraic varieties in $\HOM$ for every $G\in\Theta$ are submodules.
 However, not every submodule is a variety.

 Let us comment on the latter observation. Suppose that $W=KX$ is the free
 module over $X=\{ x_1,\ldots ,x_n\}$, that $G$ is some other module and that
 $T$ is a submodule of $KX$. According to the general theory, the corresponding
 algebraic variety $A=T'$ is of the form $A=\mu_0\Hom (KX/T,G)$, where $\mu_0$
 is the natural homomorphism. We restrict below the discussion to the simple
 case
 when $K$ is a field, and assume that the space $KX/T$ is $m$-dimensional,
 $m<n$. Under a suitable enumeration of elements in $X$, $KX/T$ admits a basis
 consisting of elements $x_1^{\mu_0},\ldots,x_m^{\mu_0}$. With the homomorphism
 $\mu_0$,  the $(n\times m)$-matrix of the kind
$$\mu_0=\left(\begin{array}{ccc}
  1&&\\
  & \ddots& \raisebox{14pt}{\huge $0$}\\
  \raisebox{5pt}{\huge $0$}&&1 \\
  \cdots& \cdots& \cdots \\
  \cdots& \cdots& \cdots
       \end{array}\right)$$
is related.

 Now assume that $G$ is $k$-dimensional. Then elements of $\Hom (KX/T,G)$ are
 presented by matrices
$$\nu=\left(\begin{array}{ccc}
       \alpha_{11}& \cdots& \alpha_{1k}\\
       \vdots& \ddots& \vdots\\
       \alpha_{m1}& \cdots& \alpha_{mk}\end{array}\right)$$
 Elements of the variety $A$ are composed of $(n\times k)$-matrices of the kind
$$\mu_0\nu=\left(\begin{array}{ccc}
       \alpha_{11}& \cdots& \alpha_{1k}\\
       \vdots& \ddots& \vdots\\
       \alpha_{m1}& \cdots& \alpha_{mk}\\
       \cdots& \cdots& \cdots\\
       \cdots& \cdots& \cdots\\
       \cdots& \cdots& \cdots\\
       \end{array}\right)$$
The upper $(m\times k)$-part of such a matrix is quite arbitrary, for $\nu$ is
arbitrary. The lower part depends of the fixed matrix $\mu_0$.
Now it is clear that there are subspaces in $\Hom (KX,G)$ that are not
algebraic varieties.

\paragraph*{8. Additional remarks.}
\stepcounter{paragraph}
 We begin with some notes concerning the case $\Theta$ is the variety of groups.
 First of all, we make a simple observation. The set $X$ is assumed to be fixed.

 We denote by $F=F(X)$ the free group over $X$. Suppose that $G$ is a
 torsion-free group and that $T$ is a normal subgroup of $F$. The kernel of every
 homomorphism $\nu\colon F/T\to G$ contains elements of finite order
from $F/T$.
 For this reason, the kernel $(G\defis\Ker)(F/T)$ contains all elements of
 finite order from $F/T$. This also means that all $\varphi\in F$ with
 $\varphi^n\in T$ for some $n$ belong to $T''$.

 If $F/T$ is a nilpotent group, then all such $\varphi$'s make up a normal
 subgroup of $F$, which is naturally denoted by $\sqrt{T}$. Now
 $T'' \supset \sqrt{T}$. In a number of cases even the equality
 $T''=\sqrt{T}$ holds. This is  something like the Hilbert theorem.

 We shall further discuss geometric equivalence of groups.

 $\bf 8.1.$ Proposition.
Suppose that $G_1$ and $G_2$ are equivalent groups. If $G_1$ is
torsion-free, then so is $G_2$.

\begin{proof}
 Assume that $G_1$ is torsion-free and that $G_2$ has a cyclic subgroup $H$ of
 order $n$. For $T$, take the verbal subgroup of $F$ over the variety of groups
 of exponent $n$. Then $(G_1\defis\Ker)(F/T)=F/T$. Now let $\nu\colon F/T\to H$ be
 a non-trivial homomorphism. It also is an element of $\Hom (F/T,G_2)$. Since
 the kernel $\Ker\nu$ differs from $F/T$, we infer that $(G_2\defis\Ker)(F/T)\ne
 F/T$. Therefore, $G_1$ and $G_2$ are not equivalent, and this contradicts the
 supposition.
\end{proof}

 If both $G_1$ and $G_2$ are periodic and if they are equivalent, then they must
 have the same exponent.

 The next question seems to be simple. Suppose that $G_1$ and $G_2$ are
 equivalent and $G_1$ is periodic. Is also $G_2$ periodic?

 We now shall consider separately the case when $X$ only consists of one element
 $x$. We shall show that, under this assumption, any two torsion-free groups are
 equivalent.

 Let $G_1$ and $G_2$ be such groups, $F=F(X)$ be an infinite cyclic group, $T$
 be a subgroup of $F$. If $T$ only consists of the unit, then
 $(G_1\defis\Ker)(F/T)=(G_2\defis\Ker)(F/T)$, and both these kernels only
 consist of the unit. If $T$ still contains something else, then $F/T$ is
 finite, and
 $(G_1\defis\Ker)(F/T)=(G_2\defis\Ker)(F/T) = F/T$.

 This conclusion does not remain valid if $X=\{ x,y\}$. If $G_1$ is a
 commutative, and $G_2$ is a non-commutative group, both torsion-free, then they
 are not equivalent.

 Already in the classical algebraic geometry it can be
 proved that the equality $\Var (G_1)=\Var (G_2)$ does not imply equivalence of
 $G_1$ and $G_2$. It is easily seen that the some holds for groups. Let us
 demonstrate this.

 Assume we are  given a surjective group homomorphism $\delta\colon G\to H$,
 where $G$
 is torsion-free and $H$ is periodic. Let $G_1=G$ and $G_2=G\times H$. Then
 $G_1$ and $G_2$ are not equivalent, but $\Var (G_1)=\Var (G_2)$.

 The following problem admits a simple solution. Find all groups $G$ for which
 all invariant subgroups $T$ of $F(X)$ are closed.

 Now we pass to the general situation. It is not difficult to observe that $T$
 is a congruence on $W=W(X)$, then $T=T''$ for some $G\in\Theta$. We
 can take the algebra $W/T$ for $G$. In this connection, we note one more
 problem which  is rather ambiguous.

 Let $T$ and $T_1$ are two congruences on $W(X)$ with $T\subset T_1$. What can
 be said concerning existence of $G\in\Theta$ such that $T''_G=T_1$?

 Such a group $G$ does not exist, for instance, if the congruence $T$ is verbal
 and $T_1$ is not. If both congruences are fully characteristic, then the
 problem can be solved in a simple manner. Indeed, the following holds:
$$(W/T_1\defis\Ker)(W/T)=T_1/T.$$
 Now if $G=W/T_1$, then $T''_G=T_1$.

 We now make a remark that also concerns with an arbitrary $\Theta$. Assume that
 $G$ and $H$ are algebras from $\Theta$. We treat $\Hom (G,H)$ as the set of
 representations of $G$ into $H$. We shall couple with it a variety of
 representations, which will be determined up to isomorphism of algebraic
 varieties.

 Let the algebra $G$ be specified by generators and defining relations. Let,
 furthermore, $X$ be the set of its generators, and $T$ be the congruence on
 $W(X)$  generated by the relations. Then there is a surjective homomorphism
 $\mu_0\colon W\to G$ with the kernel $T=\Ker\mu_0$. The corresponding algebraic
 variety $A=T'$ is specified in $\Hom (W,H)$. Then $A=\mu_0\Hom (G,H)$ is the
 variety of representations of $G$ in $H$ we are interested in. As we shall see,
 passage to another system of generators and relations leads to isomorphism of
 algebraic varieties.

 Turning back to groups, let us consider two groups $G$ and $H=\Aut (V)$ where
 $V$ is a module over some $K$. Then the question is of the variety of linear
 representations of the given $G$ in a linear group $H$. It is an algebraic
 variety in $\Hom (W,H)$; cf. \cite{PlRap}. One can consider various
 subvarieties of it and relate them with classification problem for
 representations. Also, the problem of geometrical equivalence of the groups
 $\Aut V_1$ and $\Aut V_2$ naturally arises here, $V_1$ and $V_2$ being various
 modules (over the same  $K$, or not).

 One also can speak about geometric equivalence of two representations on the
 basis of the variety $\Theta$ of all representations over a given ring $K$; see
 \cite{PlV}. In particular, two irreducible and faithful representations of
 finite groups over the same field are geometrically equivalent if and only if
 they are isomorphic. Algebraic varieties related to linear representations
 motivate various interesting ideas. This is a separate subject.

 Now again we shall make some general observations. The variety $\Theta$ is
 arbitrary, and fixed are the free algebra $W=W(X)$ and $G\in\Theta$. For every
 $\mu\colon W\to G$, the collection of the elements $x^{\mu}$ with $x\in X$ is a
 generating set of $\im\mu$. If, furthermore, $T$ is a congruence on $W$ and
 $\mu\in A=T'$, then $T\subset\Ker\mu$; this means that $T$ is induced in the
 system of defining relations of the algebra $\im\mu$. {\it Therefore, we have
 information about the generators and relations of algebras of the kind
 \/ $\im\mu$
 for all $\mu\in A=T'$}.

 For example, if the question is of groups and if $T$ contains all commutators
 $[x,y]$ for $x,y\in X$, then all subgroups of kind $\im\mu$, $\mu\in A=T'$, of
 any $G$ are commutative.

 In conclusion of the section, we note that the theory we deal with here was
 stimulated, in  considerable extent, by investigations of equations in groups.
 These investigations, in they turn, are connected with geometrical
algebra; see
 ~\cite{Gro1,Gro2,RS}. See also \cite{Paulin} as a survey of works on
geometrical algebra, in particular, of the works of E.~Rips and Z.~Sela.

 The geometric approach clears some ways for seeking solutions. Generally the
 aims of algebraic geometry are wider. We have in mind both introducing
 geometric  concepts in algebraic structures and algebraic iterpretation  of
 the arising geometric structures. With respect to this, geometric algebra and
 algebraic geometry are close to each other; however, they are oriented to
 different geometric structures. But speaking generally,
 geometry and algebra in
 either field are heavily intertwined. Our interests are focused chiefly on
 algebra. It is difficult to perceive that general algebraic varieties could be
 well-connected with substantial geometry.

\section*{\S4. Varying the variables set,
the base variety and the base algebra}
\setcounter{paragraph}{0}
\paragraph{1. Changing $\bf X$.}
\stepcounter{paragraph}
We count the variety fixed, and change the set of variables, $X$. So, we deal
with $W(X)$ and $W(Y)$.

 $\bf 1.1.$ Theorem.
If $X \subset Y$ and algebras $G_1$ and $G_2$ are $Y$-equivalent, then they are
also $X$-equivalent.

\begin{proof}
We treat $W(X)$ as a subalgebra of $W(Y)$, and take some $G \in
\Theta$. Then every $\mu\colon W(Y) \to G$ induces $\nu\colon W(X) \to G$.
On the other hand, there are several mappings $\mu\colon W(Y) \to G$
inducing a given $\nu\colon W(X) \to G$. We shall write $\nu = \mu^\alpha$.
\medskip

We shall show that if $A$ is an algebraic variety in $\HOM$, then its full
inverse $\alpha$-image $B$ is an algebraic variety in $\Hom(W(Y),G)$. Suppose
that $A = T'$, where $T$ is a binary relation on $W(X)$. It can be considered as
a relation on $W(Y)$ as well. We shall write $T = T_X$ and $T = T_y$,
respectively, in this connection. Let us check that $T'_y=B$.

First observe that if $\nu = \mu^\alpha$, then $\Ker\nu = \Ker\mu \cap W(X)$.

Now let $\mu \in B$. Then $\mu^\alpha = \nu \in A$ and
    \(T = T_X \subset \Ker\nu = W(X) \cap \Ker\mu\).
But then $T = T_Y \subset \Ker\mu$ and $\mu \in T'_Y$.
Conversely, let $\mu \in T'_Y$. Then
    \(T = T_Y \subset \Ker\mu \cap W(X) = \Ker\nu\),
where \(\nu = \mu^\alpha\), and, further,
    \(\nu \in T' = T'_X = A\) and $\mu \in B$.

Now we make some remarks on congruences.

  If $T$ is a congruence on $W(Y)$, then we have a congruence $W(X) \cap T$ on
$W(X)$. Being a binary relation, it generates a congruence on $W(Y)$. The latter
one is included in $T$ but does not, generally, coincide with $T$.

Assume now that $A$ is a variety in $\HOM$, $A' = T_X$ and $B$ is the full
inverse image of $A$, and let $T_Y = B'$. We shall verify that
    \( T_Y \cap W(X) = T_X \).

Clearly, \(B' = T_Y = \bigcap_{\mu \in B}\Ker\mu\). Furthermore,
\[
T_Y \cap W(X) = (\bigcap_\mu \Ker\mu) \cap W(X) = \bigcap_\mu(\Ker\mu \cap W(X)) =
    \bigcap_{\nu \in A}\Ker\nu = T_X.
\]

Now suppose that $T = T_X$ is a closed congruence on $W_X$. For it, we shall
construct a closed congruence $T_Y$ in $W(Y)$ so, that $T_Y \cap W(X) = T_X$.
Take $A = T'_X$ and let $B$ be the full preimage of $A$. Then take $B' = T_Y$.
It is a closed congruence, and $A' = T_X$. Now indeed,  $T_Y \cap W(X) = T_X$.

Now, we return to the theorem. Assume $G_1$ and $G_2$ are $Y$-equivalent. This
means that the congruence $T_Y$ is $G_1$-closed if and only if it is
$G_2$-closed. We shall check the same for $X$. Take $T_X$ in $W(X)$ and
assume that this
congruence is $G_1$-closed. Suppose that $T_Y$ is $G_1$-closed congruence on
$W(Y)$ such that $T_Y \cap W(X) = T_X$. The congruence $T_Y$ is $G_2$-closed. We
have to prove that $T_X$ is $G_2$-closed, too. To proceed, some additional
remarks are needed.

\smallskip
\noindent Let again $G$ be any algebra, and let $B$ be an algebraic variety in
$\Hom(W(Y),G)$ determined by some closed $T_Y = B'$. We denote by $A$ the
$\alpha$-image of $B$ in $\HOM$, and check that then $A' = B'\cap W(X)$.
Suppose that $(w,{w'}) \in A'$. Both $w$ and ${w'}$ are elements of $W(X)$, and
$(w,{w'}) \in \Ker \nu$ for every $\nu \in A$. Now if $\mu \in B$, then
$\mu^\alpha = \nu \in A$. This means that $w^\nu = w^{\prime\nu}$,
$w^\mu = w^{\prime\mu}$ and, therefore, $(w,{w'}) \in \Ker \mu$. Since this holds
for every $\mu \in B$, we conclude that $(w,{w'}) \in \bigcap_{\mu \in B} \Ker \mu
= B'$. Therefore, $(w,{w'}) \in B' \cap W(X)$.

If, conversely, $(w,{w'}) \in B' \cap W(X)$, then $w^\mu = w^{\prime\mu}$ for
every $\mu \in B$. Take any $\nu \in A$ of the kind $\nu = \mu^\alpha$.
Since $w$ and ${w'}$ belong to $W(X)$, we obtain that
$w^\nu = w^\mu = w^{\prime\mu} = w^{\prime\nu}$ and
$(w,{w'}) \in \Ker\nu$. This holds for every $\nu \in
A$, thereby $(w,{w'}) \in A'$.

In particular, if $G = G_2$, then $T_Y \cap W(X) = T_X = A'$ for some suitable
$A$. This means that the congruence $T_X$ is $G_2$-closed. As to $A$, this set
need not be an algebraic variety.

We have demonstrated that if a congruence $T_X$ is $G_1$-closed, then it is
$G_2$-closed, and the converse also holds. So, $G_1$ and $G_2$ are
$X$-equivalent. The proof of the theorem is completed.
\end{proof}

As we saw, the converse is not true: $X$-equivalence does not imply
$Y$-equi\-va\-len\-ce.

 $\bf 1.2.$ Problem.
Given finite $X$ and $Y$ with $X \subset Y$, find algebras
$G_1$ and $G_2$ that are $X$-equivalent but not $Y$-equivalent. Is this always
possible?

 $\bf 1.3.$ Problem.
Is it true or not that $X$-equivalence of $G_1$ and $G_2$
for every finite $X$ implies their $Y$-equivalence for $Y$ enumerable?

Let us add a remark which will be used later. {\it Let $X \subset Y$, and
$T=T_X$ be a binary relation in $W(X)$. If we consider $T$ as $T_Y$ in
$W(Y)$, then $T''_Y \cap W(X) = T''_X.$}

Indeed, if $T'_X=A$, then $T'_Y=B$ is a full coimage
of $A$. By the definition, $T''_X=A'$, $T''_Y=B'$, and, as we have
seen, $A'=B'\cap W(X)$.
\medskip

Now, we pass to another important subject. We count the variety $\Theta$ fixed,
and take $X$ and $Y$ either distinct or coinciding. We then have the algebras
$W(X)$ and $W(Y)$, respectively, and suppose the algebra $G$ to be given.  We
are going to co-ordinate the varieties in $\Hom(W(Y),G)$ and in $\HOM$.

We consider the set $\Hom(W(Y),W(X))$,\! which becomes the semigroup $\End\! W$
when $X$ and $Y$ coincide.

For every $s\colon W(Y) \to W(X)$ and every $\nu\colon W(X) \to G$, we have
$\mu = \nu s\colon W(Y) \to G$. This gives us the mapping
$$
\tilde s\colon \HOM \to \Hom(W(Y), G).
$$
If, furthermore, $A = T'_1$ is an algebraic variety in $\HOM$ and $B = T'_2$ is
an algebraic variety in $\Hom(W(Y),G)$, and if $T_1$ and $T_2$ are congruences
on $W(X)$ and $W(Y)$, respectively, then $\tilde s$ determines a morphism
$s\colon A \to B$ if $\nu s \in B$ for every $\nu \in A$.

 $\bf 1.4.$ Proposition. The element $s \in
\Hom(W(Y),W(X))$ determines a morphism $s\colon A \to B$ if and
only if
    $w\, T''_2\, {w'}$ implies $s(w)\,T''_1\,s({w'})$.

\begin{proof}
Assume that $\nu s \in B$ for every $\nu \in A$ and that $w\, T''_2\, {w'}$.
We need to prove
that $s(w) T''_1 s({w'})$. We have $T''_1 = A'$ and $T''_2 = B'$. Moreover,
    $A' = \bigcap_{\nu \in A} \Ker \nu$.
We shall check that, for all $\nu \in A$,  $(s(w), s({w'}))
\in \Ker\nu$ or, in other words,
$\nu s(w) = \nu s({w'})$. By the definition, $T''_2 = B' = \bigcap_{\mu \in B}
\Ker\mu$. Hence, $wT''_2{w'}$ means that $w^\mu = w^{\prime\mu}$. In particular,
this is true of $\mu = \nu s$, and then $\nu s(w) = \nu s({w'})$.

To prove the converse, we assume that
    $w\, T''_2\, {w'}$ implies $s(w)\, T''_1\, s({w'})$ and
that $\nu \in A$ is given. We shall check that $\nu s \in B$ or,
equivalentially,
$\nu s(w) = \nu s({w'})$ whenever $w\, T''_2\, {w'}$. Suppose the latter condition
is fulfilled. Then also $s(w)\, T''_1\, s({w'})$. Now, if $\nu \in A$, then
$\nu s(w) = \nu s({w'})$, $\nu s \in B$.
\end{proof}

This proposition has the following application.

 $\bf 1.5.$ Proposition.
To every morphism $s\colon A \to B$, there is an algebra homomorphism
    $\sigma\colon W(Y)/B' \to W(X)/A'$.
The converse also holds: every algebra homomorphism $\sigma$ induces a mormhism
$s$ between the respective algebraic varieties.

\begin{proof}
Assume we are given a morphism $s\colon A \to B$, $s \in \Hom(W(Y), W(X))$.
The homomorphism
    $s\colon\! W(Y)\! \to\! W(X)$
and the natural homomorphism
    $\sigma_0\colon\! W(X)$\\ $\to W(X)/T''_1$
give    $s\sigma_0\colon W(Y) \to W(X)/T''_1$.
By Proposition 1.4, the congruence $T''_2$ is included in the kernel of the
homomorphism $s\sigma_0$. Because of this, also the homomorphism
    $\sigma\colon W(Y)/T''_2 \to W(X)/T''_1$
is defined. It remains to observe that $T''_2 = B'$ and $T''_1 = A'$.
We proceeded here from $A = T'_1$ and $B = T'_2$.

We pass to the final part of the proposition.
Assume we are given a homomorphism
    $\sigma\colon W(Y)/T''_2 \to W(X)/T^{\prime\prime}_1$.
There is a related commuting diagram

\begin{center}
\unitlength=.9mm
    \begin{picture}(40,50)(0,-5)
    \put(10,5){\vector(1,0){20}}
    \put(5,30){\vector(0,-1){20}}
    \put(35,30){\vector(0,-1){20}}
    \put(10,35){\vector(1,0){20}}
    \put(7,5){\makebox(0,0)[rc]{$W(Y)/T''_2$}}
    \put(33,5){\makebox(0,0)[lc]{$W(X)/T''_1$}}
    \put(7,35){\makebox(0,0)[rc]{$W(Y)$}}
    \put(33,35){\makebox(0,0)[lc]{$W(X)$}}
    \put(20,0){\makebox(0,0)[cb]{$\sigma$}}
    \put(-1,20){\makebox(0,0)[lc]{$\sigma_1$}}
    \put(42,20){\makebox(0,0)[rc]{$\sigma_0$}}
    \put(20,40){\makebox(0,0)[ct]{$s$}}
    \end{picture}
\end{center}

\noindent
where $\sigma_1$ and $\sigma_0$ are the natural homomorphisms, and $s$ also is
specified in a natural way.

Now assume that $w T''_2 {w'}$. This means that
    \( w^{\sigma_1} = {w'}^{\sigma_1} \).
But then
    \( w^{\sigma_1 \sigma} = {w'}^{\sigma_1 \sigma} \)
and     \( w^{s\sigma_0} = {w'}^{s\sigma_0} \),
whence  \( s(w) T''_1 s({w'}) \).
Now it follows from Proposition 1.4 that we have a morphism
    \( s \colon A \to B \).
\end{proof}

We now pass to examples.

Assume we are given
    \( s\colon\! W(Y)\! \to\! W(X) \),
with the corresponding
    \( \tilde s \colon\! \HOM$\\ $\to \Hom(W(Y),G) \),
$G \in \Theta$. For every subset $B$ of $\Hom(W(Y),G)$, we define $sB = A$, a
subset of $\HOM$, by the rule:
    \( \mu \in A = sB \Leftrightarrow \mu s \in B \).
Moreover, where $T$ is a binary relation on $W(Y)$, we define the binary
relation $sT$ on $W(X)$, as above, by the rule: $w\, sT\, {w'}$ if there are $w_1$
and ${w'}_1$ in $W(Y)$ such  that
    \( w_1^s = w \), \( w_1^{\prime s} = {w'}\)
and     \( w_1T{w'}_1 \).
Again, $(sT)' = sT'$.

Let us prove this. Assume that $\mu \in (sT)'$ or, what amounts to the same,
$sT \subset \Ker \mu$. We shall check that $\mu \in sT'$, i.e. $\mu s \in T'$,
    \( T \subset \Ker \mu s \).
Let $w_1 T {w'}_1$. We have to see whether
    \( (s(w_1) )^{\mu} = s({w'}_1)^\mu \).
Take $s(w_1) = w$, $s({w'}_1) = {w'}$. Then $w\, sT\, {w'}$. Since $sT \subset \Ker
\mu$, $w^\mu = {w'}^\mu$. Thus,
    \( s(w_1)^\mu = s({w'}_1)^\mu \).
We obtain that $\mu \in sT'$.

Conversely, assume that $\mu \in sT'$. We shall check that $\mu \in (sT)'$. Let
$w\, sT\, {w'}$. We have to make sure that $w^\mu = {w'}^\mu$. By the condition,
$w_1, {w'}_1 \in W(Y)$ with
$s(w_1) = w$,  $s({w'}_1) = {w'}$  and  $w_1 T {w'}_1$. Since $\mu \in sT'$,
$\mu s \in T'$  and
    \( s(w_1)^\mu = s({w'}_1)^\mu \).
This gives us $w^\mu = w^{\prime \mu}$, $\mu \in (sT)'$.

In particular, if $B = T'$ is an algebraic variety in $\Hom(W(Y), G)$, then $A =
sB = sT' =(sT)'$ is an algebraic variety in $\HOM$.

Moreover, if $\mu \in A$, then $\mu s \in B$, and $s\colon  A \to B$
is a morphism. In general, $s\colon  A \to B$ is a morphism if
$A \subset sB$.

We further consider a particular situation.
Let $Y$ be a subset of $X$, and take for $s$ the corresponding injection
$s\colon W(Y) \to W(X)$. Then
\[
\widetilde{s}\colon \HOM \to \Hom(W(Y),G)
\]
is the projection which we have already used. If $B$ is an algebraic variety in
$\Hom(W(Y),G)$, then $sB = A$ is the corresponding full preimage. We have the
morphism $s\colon A \to B$. Cf. the proof of Theorem 1.1.

 $\bf 1.6.$ Definition.
A morphism $s\colon A \to B$, $s \in \Hom(W(Y),W(X))$ is an algebraic variety
isomorphism if it has the inverse morphism $s'\colon B \to A$.

Here
    \( s' \in \Hom(W(X),W(Y)) \),
and, for every $\nu \in A$ and $\mu \in B$,
    \( \nu ss' = \nu \) and \( \mu s's = \mu \).
If there is such $s$, then the varieties $A$ and $B$ are isomorphic.

 $\bf 1.7.$ Theorem.
Varieties $A$ and $B$ are isomorphic if and only if isomorphic are the
respective algebras $W(X)/A'$ and $W(Y)/B'$.

begin{proof}
Assume that $A$ and $B$ are isomorphic and that $s$ and $s'$ are the
respective morphisms. We also have
\( s\colon W(Y) \to W(X) \)  and \( s'\colon W(X) \to W(Y) \),
and, simultaneously, the homomorphisms
    \( \sigma\colon W(Y)/T''_2 \to W(X)/T''_1 \)
and
    \( \sigma'\colon W(X)/T''_1 \to W(Y)/T''_2 \).
Let us check that they are inverse to each other.

We need to see whether $s's(w)\, T''_2\, w$ for every $w \in W(Y)$ and
$ss'(w_1)T''_1\, w_1$ for every $w_1 \in W(X)$.

The condition $s's(w)\, T''_2\, w$ means that $(s's(w))^\mu = w^\mu$ for every
$\mu \in B$. In other notation this means that $\mu s's(w) = \mu(w)$. Since $\mu
s's = \mu$, the equality holds. Likewise, the other condition is also fulfilled.
Therefore, $\sigma$ is an algebra isomorphism.

Now assume that
    \( \sigma \colon W(Y)/T''_2 \to W(X)/T''_1 \)
is an algebra isomorphism and
    \( \sigma' \colon W(X)/T''_1 \to W(Y)/T''_2 \)
is the inverse isomorphism. Let us consider the commuting diagrams

\centerline
{\unitlength = .9mm
    \begin{picture}(40,50)(0,-5)
    \put(10,5){\vector(1,0){20}}
    \put(5,30){\vector(0,-1){20}}
    \put(35,30){\vector(0,-1){20}}
    \put(10,35){\vector(1,0){20}}
    \put(7,5){\makebox(0,0)[rc]{$W(Y)/T''_2$}}
    \put(33,5){\makebox(0,0)[lc]{$W(X)/T''_1$}}
    \put(7,35){\makebox(0,0)[rc]{$W(Y)$}}
    \put(33,35){\makebox(0,0)[lc]{$W(X)$}}
    \put(20,0){\makebox(0,0)[cb]{$\sigma$}}
    \put(-1,20){\makebox(0,0)[lc]{$\sigma_1$}}
    \put(42,20){\makebox(0,0)[rc]{$\sigma_0$}}
    \put(20,40){\makebox(0,0)[ct]{$s$}}
    \end{picture}
\hspace{2.6cm}
    \begin{picture}(40,50)(0,-5)
    \put(10,5){\vector(1,0){20}}
    \put(5,30){\vector(0,-1){20}}
    \put(35,30){\vector(0,-1){20}}
    \put(10,35){\vector(1,0){20}}
    \put(7,5){\makebox(0,0)[rc]{$W(X)/T''_1$}}
    \put(33,5){\makebox(0,0)[lc]{$W(Y)/T''_2$}}
    \put(7,35){\makebox(0,0)[rc]{$W(X)$}}
    \put(33,35){\makebox(0,0)[lc]{$W(Y)$}}
    \put(20,0){\makebox(0,0)[cb]{$\sigma'$}}
    \put(-1,20){\makebox(0,0)[lc]{$\sigma_0$}}
    \put(42,20){\makebox(0,0)[rc]{$\sigma_1$}}
    \put(20,40){\makebox(0,0)[ct]{$s'$}}
    \end{picture}
}

We shall prove that the morphisms $s\colon A \to B$ and $s'\colon B \to A$ are
mutually inverse--i.e. that $\nu ss' = \nu$ and $\mu s's = \mu$ for all
$\nu \in A$ and $\mu \in B$. Take any $w \in W(X)$, and check that
    \( \nu ss'(w) = \nu(w) \).
Clearly,
    \( \sigma_1s'(w) = \sigma'\sigma_0(w) \).
Apply $\sigma$; then
\( \sigma\sigma_1s'(w) = \sigma_0(w) = \sigma_0ss'(w) \).
This gives $ss'(w)\,T''_1\,w$. But then $\nu ss'(w) = \nu(w)$ and, furthermore,
$\nu ss' = \nu$ for every $\nu \in A$. Likewise, $\mu s's = \mu$.

 $\bf 1.8.$ Proposition.
If $s\colon A \to B$ is an isomorphism, then it is a bijection between $A$ and
$B$.

\begin{proof}
Assume that
    \( \nu_1s = \nu_2s \in B \)
for $\nu_1, \nu_2 \in A$. We apply the inverse $s'$; then
$\nu_1ss' = \nu_2ss'$ and $\nu_1 = \nu_2$. Now assume that $\mu \in B$. Then
$\mu = \mu s's$, $\mu s' \in A$, $\mu s' = \nu$ and $\nu s = \mu$.
\end{proof}

If $s\colon A \to B$ is a morphism, then $A \subset sB$. The converse also holds.
$A \subset sB$ means that $s\colon A \to B$ is a morphism. It may seem that, if
$s$ is an isomorphism, then $A =sB$; however, it is not the case. If
$\nu \in sB$,
then $\nu s \in B$ and $\nu ss' \in A$. But we cannot claim that
$\nu ss' = \nu$, for we do not know whether, and cannot conclude that,
 $\nu \in A$.

Theorem 1.7 holds also in the case $X = Y$, and then $s$ and $s'$ are elements
of $\End W$. If, in particular, $s \in \Aut W$ and $s^{-1} = s'$, then, for any
$B$, we take $A = sB$, and $s\colon A \to B$ is a variety isomorphism with the
inverse $s' = s^{-1}$. Here, $W = W(X) = W(Y)$, and some algebra $G \in \Theta$
is presumed to be given. We cannot claim that every isomorphism, for $X$ fixed,
is determined by an automorphism.

It was already noticed, and Theorem 1.7 confirms
this, that varieties can be distinguished on the level of the
respective algebras, and using properties of the algebras. The corresponding
theorem is well-known in  the classical situation; here, varieties are
distinguished geometrically, too.

We now give an application of the theorem proved.

\smallskip
\noindent Assume that $G$ and $H$ are two algebras from $\Theta$, and consider the set
$\Hom(G,H)$. Suppose  $G$ is specified by generators and relations. Let $X$
be the set of generators and $W(X)$--the corresponding free algebra. We have
the canonic surjective homomorphism
    \( \mu_0\colon W(X) \to G \).
Its kernel, $\Ker \mu_0 = T$, is regarded as a defining relation. We now pass
to      \( \Hom(W(X),H) \)        and consider here the algebraic variety $A$
determined by the congruence $T$, $A = T'$. As we know,
    \( A = \mu_0 \Hom(G,H) \).
We take, furthermore,
    \( \tau = (H\defis\Ker)(G) \),
and consider the composition homomorphism
\[
W(X) \stackrel{\mu_0}{\longrightarrow} G
        \stackrel{\mu_1}{\longrightarrow} G/\tau,
\]
where $\mu_1$ is the natural homomorphism. Then
    \( T'' = \Ker \mu_0\mu_1 = \mu^{-1}_0(\tau) \),
and $W(X)/T''$ is an algebra isomorphic to $G/\tau$.

Now suppose that the algebra $G$ is specified by generators and
relations in two different ways. Let $Y$ be the new system of
generators, let $W(Y)$ be the free algebra corresponding to it,
and let
    \( \mu_1 \colon W(Y) \to G \)
be the canonical homomorphism with the kernel $T_1 = \Ker \mu_1$.
We pass to $\Hom(W(Y), H)$ and take here $B = T'_1$. Just as above, $T''_1$ is
the kernel of the composition homomorphism
    \( W(Y) \stackrel{\mu_1}{\longrightarrow} G
        \stackrel{\mu_1'}{\longrightarrow} G/\tau \),
and the algebra $W(Y)/T''_1$ is isomorphic to $G/\tau$. The algebras $W(X)/T''$
and $W(Y)/T''_1$ are isomorphic; hence, so are the varieties $A$ and $B$.

We have proved

 $\bf 1.9.$ Proposition.
Suppose that, in $\Theta$, given are algebras $G$ and $H$, and that $G$ is
specified by generators and relations in different ways. Then the respective
algebraic varieties related with $H$ are isomorphic.

It also follows from the above that if $T_1$ is a congruence on $W(X)$ and
$T_2$ is a
congruence on $W(Y)$, and if the algebras $W(X)/T_1$ and $W(Y)/T_2$ are
isomorphic, then, for a given $G$, the algebras $W(X)/T''_1$ and $W(Y)/T''_2$
are isomorphic. Isomorphic are also the varieties $T'_1$ and $T'_2$.

\paragraph{2. Changing $\bf\Theta$.}
\stepcounter{paragraph}
Let us consider the situation when some subvariety $\Theta_0$ containing an
algebra $G$ is picked out of the given $\Theta$. We are interested in
connections between geometries for $G$ relatively to $\Theta$ and $\Theta_0$. We
count the set of variables $X$ fixed, and denote the corresponding free algebras
over $X$ by $W$ and $W_0$, respectively. To the variety $\Theta_0$ there is the
verbal congruence $T_0$ on $W$, and we have the natural epimorphism
$\mu_0 \colon W \to W_0$ with the kernel $T_0$.

Assume, furthermore, that $A_0$ is a subset of $\Hom(W_0,G)$, and take
    \( A = \mu_0A_0 = \{\mu = \mu_0\nu,\ \nu \in A_0\} \).
Of course, the passage $\nu \mapsto \mu$ determines a bijection $A_0 \to A$.

Every element $\mu \in \HOM$ is uniquely presented in the form $\mu = \mu_0\nu$,
$\nu \in \Hom(W_0,G)$; therefore, if $A$ is a subset of $\HOM$, then
$A = \mu_0A_0$, $A_0 \subset \Hom(W_0,G)$.

We aim to demonstrate that algebraic varieties here are linked with algebraic
varieties.

Every algebraic variety $A$ can be presented as $A = T'$, where $T$ is a
congruence on $W$ including the congruence $T_0$. $T$ leads to the congruence
$T/T_0$ on $W_0$. The passage $T \mapsto T/T_0$ is a bijection between
the congruences on $W_0$ and those on $W$ including $T_0$. Moreover, $wT{w'}$
holds if and only if so does
    \( w^{\mu_0}\, (T/T_0)\, w^{\prime\mu_0} \).

 $\bf 2.1.$ Proposition.
The following relationship always holds:
\[
T' = \mu_0(T/T_0)'.
\]

\begin{proof}
Let $\mu \in A = T'$, $\mu = \mu_0\nu$. We have to check that
    \( \nu \in A_0 = (T/T_0)' \).
Assume that
    \( w^{\mu_0}\, (T/T_0)\, w^{\prime\mu_0} \).
Then  $wT{w'}$ and $w^\mu = {w'}^\mu$ as well. This gives
    \( w^{\mu_0\nu} = w^{\prime\mu_0\nu} \),
    \( (w^{\mu_0})^\nu = (w^{\prime\mu_0})^\nu \)
and $\nu \in A_0$.

Conversely, if $\nu \in A_0$ and $wT{w'}$ holds, then also
    \( w^{\mu_0}\,(T/T_0)\, {w'}^{\mu_0} \)
and     \( w^{\mu_0\nu} = {w'}^{\mu_0\nu} \),
$w^\mu = {w'}^\mu$, $T \subset \Ker \mu$, $\mu \in A$. The proposition is proved.
\end{proof}

In particular, if $A = \mu_0A_0$ and $A$ is the algebraic variety determined by
the congruence $T \supset T_0$, then $A_0$ is the algebraic variety determined
by $T/T_0$.

If $A_0$  is an algebraic variety, then $A_0 = (T/T_0)'$ for some $T$, and
    \( T' = A = \mu_0A_0 = \mu_0(T/T_0)' \).
Consequently, $A = \mu_0A_0$ is also an algebraic variety.

 $\bf 2.2.$ Proposition.
Suppose that $A$ is a subset of $\HOM$, $A = \mu_0A_0$. Then
\[      A'_0 = A'/T_0   .\]

\begin{proof}
Let $T = A'$. By the definition,
    \( T = \bigcap_{\mu \in A} \Ker \mu \).
Since always $T_0 \subset \Ker\mu$, we conclude that $T_0 \subset T$; so the
quotient $T/T_0$ makes sense. We need to check that $A'_0 = T/T_0$. Assume that
    \( w^{\mu_0}\, (T/T_0)\, {w'}^{\mu_0} \).
Then also $wT{w'}$, and for every $\mu \in A$ such that $\mu = \mu_0\nu$ with
$\nu \in A_0$, we have
$w^\mu = w^{\prime\mu}$, $w^{\mu_0\nu} = w^{\prime\mu_0\nu}$
and     \( (w^{\mu_0}, w^{\prime\mu_0}) \in Ker \nu \).
This holds for any $\nu \in A_0$; therefore,
    \( (w^{\mu_0}, w^{\prime\mu_0}) \in A'_0 \).

If, conversely,
    \( (w^{\mu_0}, w^{\prime\mu_0}) \in A'_0 \)
holds, then
    \( w^{\mu_0\nu} = w^{\prime\mu_0\nu}      \)
for every $\nu \in A_1$. But then, for every $\mu \in A$,
    $w^\mu = w^{\prime\mu}$, $(w,{w'}) \in T$. From here,
    \( (w^{\mu_0}, w^{\prime\mu_0}) \in T/T_0 \).
\end{proof}

We now want to link $A''$ with $A''_0$, and $T''$ with $(T/T_0)''$.

 $\bf 2.3.$ Proposition.
For $T \supset T_0$,
\[
T''/T_0 = (T/T_0)''.
\]

\begin{proof}
Let $A = T'$, and assume that $A = \mu_0A_0$. Then $A_0 = (T/T_0)'$,
$T'' = A'$ and $A'/T_0 = A'_0$. Therefrom, $T''/T_0 = (T/T_0)''$.
\end{proof}

 $\bf 2.4.$ Proposition.
If $A = \mu_0A_0$, then
\[
A'' = \mu_0A_0''.
\]

\begin{proof}
Clearly, $A'_0 = A'/T_0$. Let $A' = T$; then
\( A'' = T' = \mu_0(T/T_0)' = \mu_0(A'/T_0)' = \mu_0(A'_0)' = \mu_0A''_0 \).
\end{proof}

We now want to prove that the property of two algebras to be equivalent do not
depend geometrically from the variety which they belong to.

 $\bf 2.5.$ Proposition.
Suppose that $G_1$ and $G_2$ are two algebras in $\Theta$ and that they both
belong to a subvariety $\Theta_0$. For given $X$, the algebras are equivalent in
$\Theta$ if and only if they are equivalent in $\Theta_0$.

\begin{proof} We proceed from the free algebras $W = W(X)$ and $W_0 = W_0(X)$. Assume
that algebras $G_1$ and $G_2$ are equivalent in $\Theta_0$. We need to check
that they
are equivalent in $\Theta$. It suffices to consider a congruence $T$ on $W$
that includes $T_0$. Then
    \( (T/T_0)''_{G_1} = (T/T_0)''_{G_2} \).
It follows that
    \( T''_{G_1}/T_0 = T''_{G_2}/T_0 \).
Hence, $T''_{G_1} = T''_{G_2}$, and $G_1$ and $G_2$ are equivalent
in $\Theta$.

Now assume that $G_1$ and $G_2$ are equivalent in $\Theta$. We shall check
their equivalence in $\Theta_0$.

We present any congruence on $W_0$ as $T/T_0$, where $T$ is a congruence on $W$
including $T_0$. We have to prove that
\[
(T/T_0)''_{G_1} = (T/T_0)''_{G_2}.
\]
By the assumption, $T''_{G_1} = T''_{G_2}$, and then
    \( T''_{G_1}/T_0  = T''_{G_2}/T_0 \).
The needed equality now follows; hence, $G_1$ and $G_2$ are equivalent in
$\Theta$.
\end{proof}

We could make use of this observation, for example, as follows. Let $\Theta$ be
the variety of groups, let $G_1$ and $G_2$ be commutative groups, and suppose we
want to know if they are equivalent. It is sufficient to proceed from the
variety $\Theta_0$ of commutative groups.

Under the assumption that $X$ is infinite, the test on equivalence runs, in
general, as follows. Given $G_1$ and $G_2$, we find $\Var G_1$ and
$\Var G_2$. If the varieties are distinct, then $G_1$ and $G_2$ are not
equivalent. Otherwise, let $\Theta_0 = \Var_{G_1} = \Var_{G_2}$. The further
checking is fulfilled in $\Theta_0$.

Our concern is, furthermore, with the following problem. Given are $\Theta$ and
$G \in \Theta$. We look for conditions under which an algebraic variety
$A = T'$ can be
specified by a finite $T$ or, alternatively, when are all $A$ for a given $G$
finitely based. The answer depends, generally, on $\Theta$. For this reason, it
would be better to put the question still another way. Suppose that
$\Theta_0 = \Var G$. Then every $A$ can be presented as $A = \mu_0A_0$, where
$A_0$ is a variety in the corresponding $\Hom(W_0,G)$. $A$ being finitely based
means here that this $A_0$ is finitely based.

There are several  problems that deserve attention.

\paragraph{3. Changing the algebra $\bf G$.}
\stepcounter{paragraph}
We are now interested in the following problem: the congruence $T$ on $W(X) = W$
being fixed, what are the connections between algebraic varieties for the algebra
$G$, its subalgebras, and its homomorphic images.

It was already noted that if $\Hom(W,H)$ is considered to be a subset of $\HOM$
whenever $H$ is a subalgebra of $G$, then for every $T$,
    \( T'_H = T'_G \cap \Hom(W,H) \).
Now, what are connections between $T''_H$ and $T''_G$?

 $\bf 3.1.$ Proposition.
Suppose that
    \( T'_G = T'_H \cup B \)
and that the intersection of $B$ and $T'_H$ is empty. Then
\[
    T''_G = T''_H \cap B'.
\]

\begin{proof}
Clearly,
    \( T''_H = \bigcap_{\nu \in T'_H} \Ker \nu \),
    \( T''_G = \bigcap_{\mu \in T'_G} \Ker \mu \),
and what is needed follows.
\end{proof}

In particular, $T''_G \subset T''_H$, and if $T''_H = T$, then
$T''_G=T$.

Obviously, every algebraic variety $A$ in $\Hom(W,H)$ can be presented as
    \( A = B \cap \Hom(W,H) \),
where $B$ is a variety in $\Hom(W,G)$.
If $A$ is not a variety, then
    \( A'' = (A')'_G \cap \Hom(W,H) \).

Now assume that we are given a surjective homomorphism $\delta\colon G \to H$
and a congruence $T$ on $W$. Then we also have
    \( \widetilde{\delta}\colon \Hom(W,G) \to \Hom(W,H) \),
$\widetilde{\delta}(\mu) = \mu\delta$, $\mu \in \HOM$. This mapping is
surjective.

    Take $T'_H = A$ and $T'_G =B$. Also, the mappings $\delta_*$ and
$\delta^*$ will be needed. We remind their definitions. The first mapping takes
a subset of $\Hom(W,H)$ into its $\widetilde\delta$-preimage in $\HOM$,
while the
other one acts into the opposite direction and gives the images. We shall
concern with the sets $\delta_*(A)$ and $\delta^*(B)$.

We immediately obtain the inclusions $B \subset \delta_*(A)$ and
$\delta^*(B) \subset A$, and they both are strict. A natural
question arises concerning the varieties $(\delta_*(A))''_G$ and
$(\delta^*(B))''_H$. Since $A$ is a variety, $(\delta^*(B))''_H
\subset A$. Generally, this inclusion is also strict. Probably,
there is nothing of interest that could be added to in the general
case.

The main difficulty is that there is no method which could enable one
to build up, from a given congruence $T$ and a homomorphism
$\delta\colon G \to H$,
new congruences only depending on $T$ and $\delta$ and making it possible
to compute the corresponding varieties.

\paragraph{4. The category of varieties.}
\stepcounter{paragraph}
Our aim here is to precise the definition of the category of algebraic
varieties
relatively to given $\Theta$ and $G \in \Theta$. We shall vary the set of
variables $X$ and deal with various algebras $W(X)$. The question is on
equational varieties.

Assume that $A = T'_1$ and $B = T'_2$ are varieties in $\Hom(W(X), G)$ and
$\Hom(W(Y), G)$, respectively. The morphism $\alpha\colon A \to B$ is now
interpreted as a regular mapping. This is a mapping for which there is a
homomorphism $s\colon W(Y) \to W(X)$ such that $\alpha(\mu) = \mu s$ for all
$\mu \in A$. This $s$ need not to be uniquely determined by $\alpha$.

 $\bf 4.1.$ Proposition.
The equality $\mu s = \mu s'$ holds for all $\mu \in A$ if and only if $s$ and
$s'$ induce the same homomorphism
\[
\sigma\colon W(Y)/B' \to W(X)/A'.
\]

\begin{proof}
Assume that $\mu s = \mu s' \in B$. Then we have the diagrams

\centerline
{\unitlength = .9mm
    \begin{picture}(40,50)(0,-5)
    \put(10,5){\vector(1,0){20}}
    \put(5,30){\vector(0,-1){20}}
    \put(35,30){\vector(0,-1){20}}
    \put(10,35){\vector(1,0){20}}
    \put(7,5){\makebox(0,0)[rc]{$W(Y)/T''_2$}}
    \put(33,5){\makebox(0,0)[lc]{$W(X)/T''_1$}}
    \put(7,35){\makebox(0,0)[rc]{$W(Y)$}}
    \put(33,35){\makebox(0,0)[lc]{$W(X)$}}
    \put(20,0){\makebox(0,0)[cb]{$\bar s$}}
    \put(-1,20){\makebox(0,0)[lc]{$\sigma_1$}}
    \put(42,20){\makebox(0,0)[rc]{$\sigma_0$}}
    \put(20,40){\makebox(0,0)[ct]{$s$}}
    \end{picture}
\hspace{2.6cm}
    \begin{picture}(40,50)(0,-5)
    \put(10,5){\vector(1,0){20}}
    \put(5,30){\vector(0,-1){20}}
    \put(35,30){\vector(0,-1){20}}
    \put(10,35){\vector(1,0){20}}
    \put(7,5){\makebox(0,0)[rc]{$W(X)/T''_1$}}
    \put(33,5){\makebox(0,0)[lc]{$W(Y)/T''_2$}}
    \put(7,35){\makebox(0,0)[rc]{$W(X)$}}
    \put(33,35){\makebox(0,0)[lc]{$W(Y)$}}
    \put(20,0){\makebox(0,0)[cb]{$\bar s'$}}
    \put(-1,20){\makebox(0,0)[lc]{$\sigma_0$}}
    \put(42,20){\makebox(0,0)[rc]{$\sigma_1$}}
    \put(20,40){\makebox(0,0)[ct]{$s'$}}
    \end{picture}
}
We shall prove that, under the assumption,
        \( \overline{s\mathstrut} = \overline{s'} \).
Suppose that $w$ is any element of $W(Y)$. Let us take
    \( w^{\sigma_1\bar s} = w^{s\sigma_0} \)
and     \( w^{\sigma_1\bar s'} = w^{s'\sigma_0} \),
and check that these elements coincide. This means that
        \( (w^s, w^{s'}) \in A' \),
i.e.    \( (w^s)^\mu = (w^{s'})^\mu \)
whenever $\mu \in A$. Since $\mu s = \mu s'$, the latter equality holds.
We have proved a half of the proposition.

To prove the converse, assume that
        \( \overline{s\mathstrut} = \overline{s'} \).
Then
\(      w^{\sigma_1\bar s} = w^{s\sigma_0} = w^{\sigma_1\bar{s}'} =
        w^{s'\sigma_0}  \),
and $(w^s, w^{s'}) \in A'$. Thus $(\mu s)(w) = (\mu s')(w)$ for every
$\mu \in A$. This goes for any $w \in W(Y)$, whence $\mu s = \mu s'$.
\end{proof}

We see that the homomorphism
        \( \sigma\colon W(Y)/B' \to W(X)/A' \)
corresponds to any regular mapping $\alpha\colon A \to B$ in a one-to-one
manner.

Now assume that we are given regular mappings $\alpha\colon A \to B$ and
$\beta\colon B \to C$, $\alpha(\mu) = \mu s_1$, $\mu \in A$ and
$\beta(\nu) = \nu s_2$, $\nu \in B$. Then
\(      \beta(\alpha(\mu)) = (\alpha\beta)(\mu) = \mu s_1s_2    \).
So $\alpha\beta$ is a regular mapping determined by the product $s_1s_2$. The
unit map $\varepsilon\colon A \to A$ is given by the unit of $\End W(X)$.

In this way we arrive at the category of algebraic varieties,
which we denote by
$\Ka_G$. The objects of the category are the algebraic varieties, and
morphisms of $\Ka_G$ are regular mappings. The passage from every $A$ to the
respective the
algebra $W(X) / A'$
in $\Theta$ is a contravariant functor from $\Ka_G$ to the category
$\Theta$.

The category $\Ka_G$ provides the endomorphism semigroup $\End A$ and the
automorphism group $\Aut A$ of every object $A$.
The group $\Aut A$ is anti-isomorphic to $\Aut(W/A')$, and $\End A$ is
anti-isomorphic to
$\End(W/A')$.
\smallskip

Now suppose that two algebras, $G_1$ and $G_2$, are given, and take the related
categories $K_1 = \Ka_{G_1}$ and $K_2 = \Ka_{G_2}$.

Let $A$ be a variety in $\Ka_1$ with a definite set of variables
$X$, and let $A = T'_{G_1}$. We assume that the congruence $T$ is
closed under $G_1$, $T''_{G_1} = T$, and set $F(A) = T'_{G_2}$.
Furthermore, we take a morphism $\alpha\colon A \to B$ from
$\Ka_1$ and suppose that it is produced by some $s\colon W(Y) \to
W(X)$. We also suppose that the variety $B$ is related with $Y$:
$B = T'_{1G_1}$. The congruence $T_1$ here is also assumed to be
$G_1$-closed. Since $\alpha\colon A \to B$ is a morphism, $T_1$
and $T$ are connected as follows: $w T_1 {w'}$ implies $s(w)\, T\,
s({w'})$. Now let $\mu \in F(A) = T'_{G_2}$. Then $\mu s(w) = \mu
s({w'})$ if $w T_1 {w'}$. This means that $\mu s \in T'_{1G_2} =
F(B)$. Thus, at the same time, we have the morphism
\[      s\colon F(A) \to F(B)   .\]
This way, a functor $F\colon K_1 \to K_2$ is defined.

In what follows, we deal with the situation when all the sets $X$ are finite.

 $\bf 4.2.$ Theorem.
If the algebras $G_1$ and $G_2$ are equivalent with respect to all $X$, then the
categories $K_1$ and $K_2$ are equivalent.

\begin{proof} The function $F\colon K_1 \to K_2$ is defined just as above. Given any
object $A$, we take $A' = T$, and also $A'' = A = T'_{G_1}$ and $F(A) =
T'_{G_2}$. Here, the congruence $T$ is both $G_1$-closed and $G_2$-closed.

Now if we are given a morphism $\alpha\colon A \to B$ determined by some
$s\colon W(Y) \to W(X)$, then $s$ yields also the morphism
$F(\alpha)\colon F(A) \to F(B)$.

Likewise, $F'\colon K_2 \to K_1$ is constructed.

Eventually, this way we obtain an equivalence between categories. Indeed, let
$A$ be
an object from $K_1$. Then $A = T'_{G_1} = T'$, where $T$ is a simultaneously
$G_1$-closed and $G_2$-closed congruence. Also, $F(A) = T'_{G_2}$ and,
further,
\(      F'(F(A)) = F'(T'_{G_2}) = T'_{G_1} = T' = A     \).
Now assume that the morphism $\alpha\colon A \to B$ is given by some
$s\colon W(Y) \to W(X)$. This $s$ yields the morphisms
\(      F(\alpha)\colon F(A) \to F(B)   \)
and
\(      F'F(\alpha)\colon A \to B     \).
Since $s$ gives both $\alpha$ and $F'F(\alpha)$, we conclude that
$F'F(\alpha) = \alpha$.
\end{proof}

 $\bf 4.3.$ Problem.
Is the converse true, i.e. does equivalence of the categories $K_1$ and
$K_2$ imply that the respective algebras are equivalent? Is there any
necessity for introducing
a new notion of equivalence for algebras via  equivalence of these categories?

\setcounter{paragraph}{0}
\section*{\S5. Algebraic logic and algebraic varieties}
\paragraph*{1. Basic concepts.}
\stepcounter{paragraph}
We will generalize here the notion of an algebraic variety. Need for such a
 generalization already appeared when the sum of two varieties was dealth with.
 The sum cannot be given by means of equational logic, and we are going to
 generalize the very
 notion of an equation, and that of a solution as well.

 At given variety $\Theta$, the corresponding Halmos algebra $U$ is considered
 instead of the free algebra $W=W(X)$; naturally, we regard the set of
 variables, $X$, to be infinite. Of course, $U$ is an algebra with equalities.
 As to the collection $\Phi$ of relation symbols, it may be either empty or
 nonempty. If $\Phi$ is empty, we speak merely of algebras, while in the case
 $\Phi$
 is nonempty we are dealing with models. A model is of the form $(G,\Phi,f)$
 where $G \in \Theta$ and $f$ is a state realizing $\Phi$ in $G$. Every such a
 model determines a homomorphism $\hat f\colon U \to V_G$. If $u\in U$, we
 also write $\hat f(u)=f*u$ and consider this element of $V$ as a subset of
 $\HOM$,
 i.e. of the same "affine space". If $\Phi$ is empty, we write $f=f_G$.

 We now consider a formula $u$ as an equation; a formula of kind $w\equiv {w'}$ is an
 equation of a special form. The point $\mu\in\HOM$ is a solution of the
 "equation" $u$ in a model $(G, \Phi,f)$ if $\mu\in f*u$. This definition
 conforms well with what we saw in the case of equations of kind $w \equiv {w'}$,
 and such a generalization is useful also in classical geometry.

 We now can consider $f*u$ as the algebraic variety related with the model given
 and determined by the formula $u$. In a database, $f*u$ is the answer to the
 query $u$ at the state $f$. Therefore, the answer to a query can be treated as
 an algebraic variety.

 In general, such a generalized algebraic variety in $\HOM$ is given by a
 collection $T$ of formulas from $U$. In a database context $T$ can be
 thought
 to be a system on queries. The common reply to the system is the corresponding
 variety. Again, a Galois correspondence can be established between such collections
 $T$ and subsets of $\HOM$, and it is in agreement with what was said earlier.
 The connection is set up as follows.

 If $T$ is a set of formulas (elements of $U$) and $(G, \Phi,f)$ is a model,
 $G\in\Theta$, then we let
 $$ T'= \bigcap_{u\in T} (f*u).$$
 Here, all the $f*u$ are subsets of $\HOM$, and $T'$ is a subset of $\HOM$.

 If, on the other hand, $A$ is a subset of $\HOM$, then
 $$A'=\{ u,\ A\subset f*u \} .$$

 As mentioned above, this definition agrees with what appears in the case of
 equational logic. Let us explain this point.

 First of all, we remind that if $w\equiv {w'}$ is an equality, then
    $\mu\in f_G  *(w\equiv {w'} )\Leftrightarrow w^{\mu}={w'}^{\mu}$.
Hence,
 $$
 \mu\in f_G *(w\equiv {w'})\Leftrightarrow (w,{w'})\in \Ker \mu.
 $$
 Now assume that all the formulas $u\in T$ are equalities.
 We shall check that
 $$
 {\bigcap _{w\equiv {w'}\in T} f*(w\equiv {w'})} \,=\,
    \{ \mu,\ T\subset \Ker \mu \}.
$$
 If $\mu \in\bigcap f*(w \equiv {w'})$, then $(w,{w'})\in\Ker\mu$ for every $w\equiv
 {w'}\in T$,  and $T\subset\Ker\mu$.
 If, conversely, $T\subset\Ker\mu$, then $w\equiv
 {w'}\in T$ implies $w^{\mu}={w'}^{\mu}$, and $\mu\in f*(w\equiv {w'})$. This holds
 for  all $w\equiv {w'}\in T$.

 Now assume that, for a given $w\equiv {w'}$, $A$ is a subset of $f*(w\equiv {w'})$.
 This means that $w^\mu ={w'}^{\mu}$ for every $\mu\in A$, i.e.
 $(w,{w'})\in\bigcap_{\mu\in A}\Ker\mu$. If, on the other hand,
 $(w,{w'})\in\bigcap_{\mu\in A}\Ker\mu$, then, for every $\mu\in A$,
 $w^\mu={w'}^{\mu}$, $\mu\in f*(w\equiv {w'})$, and furthermore, $A\subset
 f*(w\equiv {w'})$.

 It also is easily seen that the passages $T\to A=T'$ and $A\to T=A'$ really
 give a Galois correspondence. Obviously, $T_1 \subset T_2$ implies $T'_2 \subset
 T'_1$ and $A_1 \subset A_2$ implies $A'_2 \subset A'_1$. Let us verify that
 $A\subset A''$ and $T\subset T''$.

 We know that $A''=(A')'=\bigcap_{u\in A'} f*u$. Let $\mu\in A$, and let $u\in
 A'$, $A\subset f*u$. Then $\mu\in f*u$, and this is the case for all $u\in A$
 and $\mu\in A''$.

 Furthermore, $T''=(T')'=\{ u,\ T'\subset f*u\} $. Here, $T'=\bigcap_{u\in T}
 f*u$, and then $T'\subset f*u$ for every $u\in T$. Hence, $u\in T$ implies
 $u\in T''$.

 We now call a subset $A\subset \HOM$ an algebraic variety for the model
 $(G,\Phi ,f)$ if $A=T'$ for some $T\subset U$. $A$ is an algebraic variety if
 an only if $A=A''$.

 It follows immediately from the definitions that, for every $A$, the set $T=A'$
 is a Boolean filter in $U$. A closed Boolean filter $T$ is one with $T=T''$.

 Now we shall consider the most elementary connections arising on the higher
 level we have
 reached. First of all, we notice that, for every $A\subset \HOM$, we formerly
 had a
 congruence $A'=T$ of the free algebra $W$ that depended on the algebra $G$.
 Here we have a Boolean filter $A'=T$ of the Halmos algebra $U$ that depends,
 in  general, on the model $(G,\Phi,f)$. Clearly, the congruence $T$ can also
 be presented in the Halmos algebra $U$.

 Let as consider separately the case when $A$ is all the set $\HOM$. Then the
 congruence $T=A'$ is the equational theory $T(G)$ of the algebra $G$. If,
 furthermore, $T=A'=\{ u,\ A\subset f*u\} $, then for $A= \HOM$ we obtain that
 $f*u=\HOM$, and $T=\Ker\hat f$, $T=T(G,\Phi,f)$ is the elementary theory of the
 model under consideration. Here $T$ is a filter of the Halmos algebra $U$.
 A question naturally arises, for  which $A$ the respective $T=A'$ is a filter
 of $U$. The condition looks as
    $ A\subset f*u\Rightarrow  A\subset f*\forall (X)u$,
and the right hand inclusion is equivalent to
 $A\subset \forall (X)(f*u)$.

 $\bf 1.1.$ Proposition.
 If for some set $A$ the corresponding $T=A'$ is a filter of the Halmos algebra
 $U$,  then  either $T$ is an improper filter or $T=T(G,\Phi,f)$.

\begin{proof}
  Assume first that $A$ is empty. Then $A\subset f*u$ for all $u$, and $T=A'=U$
  is an improper filter, which contains no model.

  Now assume that $A$ is nonempty, $A\subset f*u$, and $T=A'$ is a filter. Then
  $\forall (X)u\in T$ and $A\subset \forall (X)(f*u)$. If $f*u$ is a proper
  subset of $\HOM$, then the set $\forall (X)(f*u)$ is empty. This contradicts
  the assumption. Consequently, $u\in\Ker\hat f$ for every $u\in T$,
  and $T$ is included in the elementary theory of the model $(G,\Phi,f)$.
  On the other hand, for every $A$,
  the elementary theory of the model is included in $A'=T$. Thus, if $A$ is
  nonempty and $A'=T$ is a filter, then $T=T(G,\Phi,f)$.
  \end{proof}

  We already know that it is the case if $A=\HOM$. So, one thing more that we
  have to comprehend is wether it is possible when $A$ is a proper subset of
  $\HOM$.

  Now assume that $T$ is a filter of the Halmos algebra $U$.

 $\bf 1.2.$ Proposition.
  Either $T'$ is empty or $T'=\HOM$.

\begin{proof}
  Clearly, $T'=\bigcap_{u\in T}f*u$. Hence, if $f*u$ is empty for some
  $u$, then so is $T'$. If $f*u$ is a proper subset of $\HOM$, then $\forall
  (X)(f*u)$ is empty, and $f*\forall (X)u$ is empty as well. But $\forall
  (X)u\in T$, hence, $T'$ is empty. Consequently, if $T'$ is nonempty, then
  $f*u=\HOM$ for all $u\in T$, and $T'=\HOM$.
\end{proof}

Now let $A$ be nonempty subset of $\HOM$ such that $T=A'$ is a filter. Then
 $A''=\HOM$. In particular, if $A$ is a proper algebraic variety, then $T=A'$ is
 not a filter of the Halmos algebra $U$.

 Also, the following notes are self-evident.

 If $T$ only contains the zero, then the set $T'$ is empty. If the unit is the
 single element of $T$, then $T'=\HOM$. For $T$ empty, again $T'=\HOM$,
 and if $T=U$, then $T'$ is empty. The case when $T$ consists of equalities was
 considered earlier.

 Clearly, the intersection of varieties is a variety, and for every $A$, $A''$ is
 the intersection off all varieties including $A$.

 $\bf 1.3.$ Proposition.
 The sum of a finite number of the varieties is also a variety.

\begin{proof}
  It suffices to deal with two summands. Suppose that $A=T'_1$ and $B=T'_2$,
  and let $T=T_1\vee T_2$ be the set of formulas $u\vee v$ with $u\in T_1$ and
  $v\in T_2$. Then
$$
T'=\bigcap_{\stackrel{{\mbox{\scriptsize $u\! \in\! T_1$}}}{v\in T_2}}
    f*(u\vee v) =
   \bigcap_{\stackrel{{\mbox{\scriptsize $u\! \in\!  T_1$}}}{v\in T_2}}
    ((f*u)\cup (f*v))=
$$
$$
(\bigcap_{u\in T_1}(f*u))\cup (\bigcap_{v\in
 T_2}(f*v))=T'_1 \cup T'_2 =A\cup B.$$
  So, $A\cup B=T'$ is an algebraic variety.
\end{proof}

 We note in addition that if $A=T'$ and $T$ is finite, say, $T=\{
 u_1,\ldots,u_n\} $, then $A$ can be given by the formula
 $u=u_1\land\cdots\land
 u_n$. In this case, the set $\neg A$ is also a variety, and is determined by
 the formula $\neg u$. In general, the complement of a variety  need not  be a
 variety itself. This, particularly, means that the system of all
varieties for a
 given model $(G,\Phi,f)$ may be not a subalgebra of the Boolean algebra
 $\MM_{G}$  of all subsets of $\HOM$.

 On the other hand, the following proposition holds.

 $\bf 1.4.$ Proposition.
 If $A=T'$ is an algebraic variety, then for every $s\in \End W$, the set $sA$
 is also an algebraic variety.

\begin{proof}
 Assume that $sT=\{ su,\ u\in T\}$. We shall verify that
 $$(sT)'=sT'.$$
  Assume $\mu\in sT'$. Then
    $\mu s\in T'=\bigcap_{u\in T} (f*u)$.
For every $u\in T$, $\mu s\in f*u$ and $\mu\in s(f*u)=f*su$. Here, $su$ is
any element of  $sT$, and $\mu\in (sT)'$. Now assume that $\mu\in (sT)'$,
  $\mu\in\bigcap_{u\in T}(f*su)$. Then for any $u\in T$, $\mu\in
  f*su=s(f*u)$ and $\mu s\in f*u$. Since it is so for all $u\in T$, we conclude
  that $\mu s\in T'$ and $\mu\in sT'$.
\end{proof}

 In particular, $sA$ is a variety determined by $sT$.

 Assume now that $A = T'$ and that $Y$ is a subset of $X$.
 Let $\exists (Y)T$ stand for
 $\{\exists(Y)u,\ u \in T\}$. We shall verify the inclusion
  $$\exists (Y)T'\subset (\exists (Y)T)'.$$

 Given $\mu\in\exists (Y)T'$, we choose $\nu\in T'$ so that $\mu (x)=\nu (x)$
 for all $x\notin Y$. For every $u\in T$, we have $\nu\in f*u$, and then
 $\mu\in\exists (Y)(f*u)=f*\exists (Y)u$. Since $\exists (Y)u$ is an element of
 $\exists (Y)T$, we conclude that $\mu\in (\exists (Y)T)'$. The needed inclusion
 is proved.

 The converse inclusion is, in general, not true. However, if the set $T$ is
 finite and $A=T'$, then $\exists (Y)A$ is a variety for every $Y\subset X$.
 Indeed, if $T$ is finite, then $A=f*u$ where $u$ is a conjunction of elements of
 $T$. Then $\exists (Y)A=f*\exists (Y)u$.

 Apparently, it is not, in general, the case that for any variety $A$ the set
 $\exists (Y)A$ is again a variety.

 $\bf 1.5.$ Proposition.
 Suppose that $A$ is a subset of $\HOM$ and that $s\in \Aut W$. Then
$$(sA)'=sA'.$$

\begin{proof}
  Assume that $u\in (sA)'$. Then $sA\subset f*u$, $A\subset
  s^{-1}(f*u)=f*s^{-1}u$, $s^{-1}u\in A'$. Hence, $u\in sA'$. Conversely, assume
  that $u\in sA'$, i.e. $u=sv$ with $v\in A'$. Then $A\subset f*v$, $sA\subset
  s(f*v)=f*sv=f*u$, whence $u\in (sA)'$.
\end{proof}

We note here that the inclusion $sA'\subset sA$ holds for arbitrary $s\in\End
 W$.

 It follows from the proposition that if $T=A'$ is a closed collection of
 formulas, then so is $sT$ for every $s\in\Aut W$.

 Moreover, the two preceding propositions imply that, for every $s\in\Aut
 W$ and all $A$ and $T$,
$$(sA)''=sA'',\quad (sT)''=sT''.$$

 Let again $Y\subset X$, and let $A$ be arbitrary. Then we have the inclusion
$$\exists (Y)A'\subset (\exists (Y)A)'.$$
 Indeed, suppose that $u\in\exists (Y)A'$, $u=\exists (Y)v$, $v\in A'$, $A\subset
 f*v$. Then $\exists (Y)A\subset\exists (Y)(f*v)=f*\exists (Y)v=f*u$. This gives
 $u\in (\exists (Y)A)'$.

 The converse inclusion does not hold.

 $Y$ may be the whole set $X$. All formulas in $\exists (X)A'$ are closed. If $A$
 is nonempty, then $\exists (X)A=\HOM$ and $(\exists (Y)A)'$ is a filter of $U$.
 It does not consist of closed formulas only.

 Let us present some remarks concerning the action of semigroup $\End G$. A more
 general situation involving the endomorphism semigroup of the model
 $(G,\Phi,f)$ could be considered. We, however, shall restrict ourselves to
 $\End G$, and the case in point will then be empty $\Phi$, i.e. we shall only
 deal  with algebras. Let us consider two instances.

\smallskip
\noindent{\bf 1.} Suppose $\delta\in\Aut G$. There is an
automorphism $\delta_*$ of the Halmos
 algebra $\MM_G$ which corresponds to this element. For every formula $u\in U$,
$$\delta_* (f_G *u)=f_G *u.$$
 It follows here from that $\delta_* (A)=A$ for any algebraic variety $A$. This
 also means that $\mu\delta\in A$ whenever $\mu\in A$, i.e. every algebraic
 variety is invariant under the group $\Aut G$.

 It is not, in general, the case with arbitrary endomorphism, and the argument
 does  not work also when $\Phi$ is nonempty.

\smallskip
\noindent{\bf 2.} Now suppose that $\delta\in\End G$ and that an universal formula $u$ is of
 kind $$w_1\equiv w_1'\vee\cdots\vee w_n\equiv w_n',$$
i.e.  is a pseudoindentity regarded as a pseudoequation. In this case the
variety
 $f_G *u$ is invariant relatively to $\delta$. Any variety specified by a
 collection of pseudoindentities is also invariant under $\End G$.

\paragraph*{2.  Generalized varieties, topology and other subjects.}
\stepcounter{paragraph}
 We shall make some notes on topology on $\HOM$ related to generalized varieties
 we are considering here. The model $(G,\Phi,f)$ is supposed to be fixed.

 A formula $u\in U$ is said to be positive if it can be built up from the basic
 ones without using negations. All the other operations from the  Halmos algebra
 signature are admitted.

 We shall deal with varieties in $\HOM$ determined by collections $T$ of
 positive formulas. Such varieties might be termed positive. The intersection of
 any collection of positive varieties is again a positive variety, and so is the
 sum of a finite number of positive varieties. This induces a topology on $\HOM$
 with positive varieties as closed sets. This topology is more refined than the
 defined above Zariski topology.

 As to the Zariski topology, its closed sets here are exactly the algebraic
 varieties determined by collections of formulas which we consider to be
 pseudoindentities. This is a particular case of positive formulas. It is
 naturally here to proceed from an empty $\Phi$ and consider algebras
 $G\in\Theta$. In the classical case, pseudoindentities reduce to indentities.

 In the rest we shall change the notation related to Galois correspondences, having
in  mind a subsequent comparison of them.

 As in \S 3, we set for every $A\subset\HOM$:
$$T=A'=\bigcap\limits_{\mu\in A}\Ker\mu,$$
 but we now treat $T$ as a set of equalities of the Halmos algebra $U$. Moreover,
 let
$$T=A^{\vee}=\{ u\in U,\ A\subset f*u\}. $$
 We presuppose here that  $(G,\Phi,f)$ is a model.

 In addition, we denote by $U_0$ the set of all equalities of $U$. Then, as we
 know,
$$A'=A^{\vee}\cap U_0 \quad {\rm if}\quad \Phi \quad {\rm is\quad empty}.$$
 If, furthermore, $T$ is a collection of formulas, then
$$A=T^{\vee}=\bigcap\limits_{u\in T}(f*u).$$
 If $T$ is a set of equalities, then $T^{\vee}$ conincides with
        $T'=\{\mu,\  T\subset\Ker\mu\}$.

 We are now interested in relationships between $A''$ and $A^{\vee\vee}$, and
 between $T''$ and $T^{\vee\vee}$, $\Phi =\emptyset$. We always have $A\subset
 A^{\vee\vee} \subset A''$ and $T\subset T''\subset T^{\vee\vee}$, provided $T$
 is a collection of equalities. Immediately,
$$
A''=(A')'=(A^{\vee}\cap U_0)^{\vee}\supset A^{\vee\vee},
$$
and, if  $T$ is a collection of equalities, then
$$T''=(T')'=(T^{\vee})'=(T^{\vee})^{\vee}\cap U_0 =T^{\vee\vee}\cap U_0 .$$

 We have previously introduced the notion of geometric equivalence of two
 algebras $G_1$ and $G_2$ from $\Theta$. Here we shall call this kind of
 geometric equivalence weak, and define the strong geometric equivalence by the
 condition: for every collection $T$ of formulas
$$T^{\vee\vee}_{G_1} =T^{\vee\vee}_{G_2}.$$
 The motivation for this distinction is provided by the following proposition.

 $\bf 2.1.$ Proposition.
 If two algebras, $G_1$ and $G_2$, are strongly equivalent, then they are also
 weakly equivalent.

\begin{proof}
 Assume that $G_1$ and $G_2$ are strongly equivalent, and let $T$ be a
 collection of equalities. Then
$$T''_{G_1}=T^{\vee\vee}_{G_1}\cap U_0 =T^{\vee\vee}_{G_2}\cap U_0 =T''_{G_2},$$
 and $G_1$ and $G_2$ are weakly equivalent.
\end{proof}

 We now shall consider the case $A$, the set of points, only consists of one
 point $\mu\colon A=\{\mu\}$, and find $A^{\vee}$ and $A^{\vee\vee}$--closures
 of the point.

 First of all, we shall show that $A^{\vee}$ is an ultrafilter of the Boolean
 algebra $U$. We know that $A^{\vee}$ is a filter, and we have to check that
 either $u\in A^{\vee}$ or $\neg u\in A^{\vee}$ for every $u\in U$.

 Clearly,
$$(f*u)\cup\neg (f*u)=\HOM.$$
 If $\mu\in f*u$, then $u\in A^{\vee}$. If $\mu\in\neg (f*u)=f*\neg u$, then
 $\neg u\in A^{\vee}$. Therefore, the filter $A^{\vee}$ is an ultrafilter for
 any singletone $A$.

 We move to $A^{\vee\vee}=\bigcap\limits_{u\in A^{\vee}} (f*u)$.

 As $A^{\vee\vee}\subset A''$, for every $\nu\in A^{\vee\vee}$,
 $\Ker\mu\subset\Ker\nu$.
 Now assume that we have passed from the variety $\Theta$ to the variety
 $\Theta'$
 of $G$-algebras and, consequently, from the Halmos algebra $U$ to $U'$. Then the
 point $\mu$ is weakly closed: $A=A''$ implies $A^{\vee\vee}=A$.

 Let us, furthermore, consider the following question: given a filter of Halmos
 algebra $U$, what  can we say about its closure $T^{\vee\vee}$?

 $\bf 2.2.$ Proposition.
 If $T$ is a filter, then either $T^{\vee\vee}=U$ or $T^{\vee\vee}=T(G,\Phi,f)$.

\begin{proof}
 Assume that $T$ is a filter of the Halmos algebra $U$. Accordingly to
Proposition
 1.2, the variety $T^{\vee}$ is either empty or equal to $\HOM$. In the first
 case $T^{\vee\vee}=U$, in the second -- $T^{\vee\vee}$ is $T(G,\Phi,f)$ and,
 furthermore, $T\subset T(G,\Phi,f)$. So if the inclusion does not hold,
 then $T^{\vee\vee}=U$.
\end{proof}

 Hence, if the question is about bijection between algebraic varieties for a
 given model $(G,\Phi,f)$ and closed collections $T$ of formulas, the unique
 filter occurring here is the elementary theory $T(G,\Phi,f)$.

 $\bf 2.3.$ Problem.
Assume that $H$ is a subalgebra of $G$ and $(H,\Phi,f_H)$ is the model
 induced by $(G,\Phi,f)$. Consider $\Hom (W,H)$ as a subset $A$ of $\HOM$.
 Investigate the connection between the Boolean filter $A^{\vee}$ and the filter
 (elementary theory) $T(H,\Phi,f_H)$. Compare here $A^{\vee\vee}$ with
 $T^{\vee}(H,\Phi,f_H)$ and find $T^{\vee\vee}(H,\Phi,f_H)$.

 This is a general problem, but it can also be considered in connection with
 various special situations.

 Let us make, finally, the following obvious note:
\begin{quote}
     if $T$ is an ultrafilter of a Boolean algebra $U$, then it is closed,
     $T^{\vee\vee}=T$, if and only if $\bigcap\limits_{u\in T}(f*u)$ as a
     nonempty set.
\end{quote}

\paragraph{3.  Passage to submodels.}
\stepcounter{paragraph}
 We are considering models $(G,\Phi,f)$, $G\in\Theta$. A submodel of such a
 model looks like $(H,\Phi,f_H)$. Here, $H$ is a subalgebra of $G$ and $f_H$
 is the restriction of $f$ to $H$.

 We shall specify some details. The algebras are presented in the form
 $G=(G_i,\ i\in\Gamma)$ and $H=(H_i,\ i\in\Gamma)$, respectively. For
 every $i\in\Gamma$, $H_i$ is a subset of $G_i$. Now assume that $\varphi$ is a
 relation symbol from $\Phi$ of type $\tau =(i_1,\ldots ,i_n)$. There are
 Cartesian products
$$G_{i_1}\times\cdots\times G_{i_n}\quad{\rm and}\quad H_{i_1}
\times\cdots\times H_{i_n},$$
 and the latter is considered to be a subset of the former one. So,
$$f_H(\varphi)=f(\varphi)\cap (H_{i_1}\times\cdots\times H_{i_n}).$$

 This defines  $f_H$ as a restriction of the function $f$. As usually, we also
 treat $\Hom (W,H)$ as a subset of $\HOM$. Then, for every formula $u\in U$,
 $f_H*u$ is a subset of $\Hom (W,H)$, and a subset of $\HOM$ as well.

 We now shall demonstrate that the following holds for $u\in U$ a basic formula:
$$
f_H *u=(f*u)\cap\Hom (W,H).
$$
 Assume that $u=\varphi(x_1,\ldots,x_n)$, $\varphi\in\Phi$,
 $\tau=\tau(\varphi)=(i_1,\ldots,i_n)$. Then $\mu\in f*u$ if and only if
$$
(\mu(x_1),\ldots,\mu(x_n))\in f(\varphi).
$$
 At the same time, $\mu\in f_H *u$ for $\mu\colon W\to H$ if and  only if
$$(\mu(x_1),\ldots,\mu(x_n))\in f_H
(\varphi)=f(\varphi)\cap(H_{i_1}\times\cdots\times H_{i_n}).$$

 Now let $\mu\in f_H *u$. Then $(\mu(x_1),\ldots,\mu(x_n))\in f_H
 (\varphi)=f(\varphi)\cap
 (H_{i_1}\times\cdots\times H_{i_n})$. Here $(\mu(x_1),\ldots,\mu(x_n))\in
 f(\varphi)$ implies $\mu\in f*u$. By the definition, also $\mu\in\Hom (W,H)$.
 Hence, $\mu\in (f*u)\cap\Hom (W,H)$.

 Let, conversely, $\mu\in ((f*u)\cap\Hom (W,H))$. Then
 $(\mu(x_1),\ldots,\mu(x_n))\!\!\in f(\varphi)$. Moreover, $\mu\in\Hom (W,H)$
 implies that $(\mu(x_1),\ldots,\mu(x_n))\in H_{i_1}\times\cdots\times
 H_{i_n}$. Therefore, $(\mu(x_1),\ldots,\mu(x_n))\in
 f(\varphi)\cap(H_{i_1}\times\cdots\times H_{i_n})=f_H (\varphi)$, and we come
 to $\mu\in f_H *u$. We see that for a basic formula $u$ the equality
$$f_H *u=(f*u)\cap\Hom (W,H)$$
 holds.

 We shall call this equality the {\it fundamental equality} for the formula
 $u$.
 We also shall deal with the inclusion
$$f_H *u\subset (f*u)\cap\Hom (W,H),$$
 called the {\it fundamental inclusion} (for $u$).

 A formula $u$ is said to be open, or quantifier-free, if it is built up
 from basic  formulas without using quantifiers.

 $\bf 3.1.$ Proposition.
 The fundamental equality holds for all open formulas. The fundamental inclusion
 holds for all positive formulas.

\begin{proof}
 Both the equality and the inclusion are fulfilled for basic formulas, and we
 shall proceed by induction. Let $M_0$ be the set of all formulas $u$ for which
 the fundamental equality holds, and let $M_1$ be the set of those $u$ with the
 fundamental inclusion.

 We shall show that $M_0$ is invariant relatively to Boolean operations and the
 action of the semigroups $\End W$. This will mean that all open formulas belong
 to $M_0$. We shall also show that $M_1$ is invariant relatively to $\vee$ and
 $\land$, and relatively to $\End W$ and quantifiers. Consequently, all
 positive formulas belong to $M_1$.

 Assume that $u\in M_1$ and $s\in \End W$. We have to prove that $su\in M_1$.
 Let $\mu\in f_H *su$. Then $\mu\in s(f_H *u)$, $\mu s\in f_H *u\subset
 (f*u)\cap\Hom (W,H)$, $\mu\in (f*su)\cap\Hom (W,H)$. So $su\in M_1$.

 Now assume that $u\in M_0$; we shall verify that $su\in M_0$. As $u\in M_1$,
 the fundamental inclusion holds for $su$. Let $\mu\in (f*su)\cap\Hom (W,H)$.
 Then $\mu s\in (f*u)\cap\Hom (W,H)=f_H *u$, $\mu\in f_H *su$. For $su$, the
 converse, $su\in M_0$, holds. Now let $u\in M_1$ and $Y\subset X$. We have to
 check that $\exists (Y)u\in M_1$.

 Let $\mu\in f_H *\exists (Y)u=\exists (Y)(f_H*u)$. We can select $\nu\in f_H*u$
 so, that $\mu (x)=\nu (x)$ for $x\notin Y$. Then $\nu\in f*u$; therefore,
 $\mu\in\exists (Y)(f*u)=f*\exists (Y)u$, $\mu\in (f*\exists (Y)u)\cap\Hom
 (W,H)$. Consequently, $\exists (Y)u\in M_1$.

 Assume, furthermore, $u\in M_0$; we shall check that $\neg u\in M_0$. Let
 $\mu\in f_H*\neg u=\neg (f_H*u)$. Here $\mu\colon W\to H$ does not belong to
 $f_H*u=(f*u)\cap\Hom (W,H)$. Thus, $\mu\in\neg (f*u)=f*\neg u$, $\mu\in (f*\neg
 u)\cap\Hom (W,H)$.

 To verify the converse, let $\mu\in (f*\neg u)\cap\Hom (W,H)$. Then $\mu\in\neg
 (f*u)$ and $\mu$ does not belong to $(f*u)$. Accordingly, $\mu$ does not belong
 to $(f*u)\cap\Hom (W,H)=f_H*u$. It follows that $\mu\in\neg (f_H*u)=f_H*\neg
 u$. We conclude that $\neg u\in M_0$.

 Now assume that $u_1,u_2\in M_1$. We shall show that $u_1\land u_2$ and
 $u_1\vee u_2$ also belong to $M_1$. We have:
$$
f_H*(u_1\land u_2)=
$$
$$
(f_H*u_1)\cap (f_H*u_2) \subset ((f*u_1)\cap\Hom (W,H))\cap
    ((f*u_2)\cap\Hom (W,H))=
$$
$$
(f*u_1)\cap (f*u_2)\cap\Hom (W,H) =
    f*(u_1\land u_2)\cap\Hom (W,H),
$$
 and $u_1\land u_2\in M_1$. Likewise,
$$
f*(u_1\vee u_2)\cap\Hom (W,H) =
    ((f*u_1)\cup (f*u_2))\cap\Hom (W,H) =
$$
$$
((f*u_1)\cap\Hom (W,H))\cup ((f*u_2)\cap\Hom (W,H)) \supset
    (f_H*u_1)\cup (f_H*u_2) =
$$
$$
f_H*(u_1\vee u_2),
$$
 and $u_1\vee u_2\in M_1$.

 This way we can also demonstrate that $u_1\vee u_2\in M_0$ and $u_1\land u_2\in
 M_0$ if $u_1,u_2\in M_0$. In both cases, we write equalities instead of
inclusions, and what is needed follows.

 The proposition is proved.
\end{proof}

 As to positive formulas, we note that if $u$ is positive, then $\neg u$ is
 negative and, moreover, if $u\in M_1$, then $\neg u$ may do not belong to $M_1$.
 Indeed, the inclusion $f_H*u\subset f*u$ changes to
 $f_H*\neg u\supset f*\neg u$.

 Let us make some notes regarding to the set $M_0$ in connection with
 quantifiers $\exists (Y)$. We shall show that $M_0$ does not possess the
 invariance property relatively to such quantifiers.

 Suppose that an algebra $G\in\Theta$ and a subalgebra $H\subset G$ are given.
 Let $u\in M_0$ be so selected that
    \begin{enumerate}
\item[1.] $f_H*u=\emptyset$,
\item[2.] there is $\nu\in f*u$ such that $\nu (x)\in H$ for some $x\in H$.
    \end{enumerate}
 Let $Y=X\backslash\{ x\}$. Then $\exists (Y)f_H*u=f_H*\exists (Y)u=\emptyset$.
 We choose $\mu\in\exists (Y)(f*u)=f*\exists (Y)u$ so that $\mu (x)=\nu (x)$.
 Then $\mu\in\exists (Y)(f*u)=f*\exists (Y)u$ and $\mu\in f*\exists (Y)u\cap\Hom
 (W,H)$. The fundamental equality is not fulfilled here, and $\exists (Y)u\notin
 M_0$.

 More specifically, assume that $\Theta$ is the variety of groups, $G\in\Theta$
 and $H$ is the unit subgroup. Take the formula $x\ne y$ for $u$. Then $u$ is
 open, and $u\in M_0$. Moreover, $f_H*u=\emptyset$. If $G$ is not trivial, then
 there is the needed $\nu$ such that $\nu (x)=1$ and $\nu (y)\ne 1$ for
 $y\ne x$. Here
 $f=f_G$, and $\mu =\mu_0$ is the trivial homomorphism. The necessary conditions
 are all fulfilled.

\paragraph{4.  Geometry on the level of quantifier-free logic.}
\stepcounter{paragraph}
 We have considered the levels of equational logic and pseudoequational logic,
 as well as that of first-order logic. Now the initial logic is the logic of
 open formulas. We shall also make a note regarding the case of positive
 formulas

 The Galois correspondence is assigned between subsets of $\HOM$ (the model
 $(G,\Phi,f)$ is supposed to be fixed) and collections $T$ of open formulas in a
 Halmos algebra $U$. Just as above, we denote the connection by ${}^\vee$. If
 $A$ is a subset of $\HOM$, then $T=A^{\vee}$ is
 the set of open formulas $u$  with $A\subset f*u$.
 If $T$ is a set of open formulas, then
 $A=T^{\vee}=\bigcap_{u\in T}(f*u)$. $A^{\vee\vee}$ and $T^{\vee\vee}$
 are also to be considered in this sense.

 Let us set $A=\HOM$ and find $A^{\vee}$. If $u$ is an open formula with
 $A\subset f*u$, then $u$ belongs to the open theory of the model $(G,\Phi,f)$
 which we denote here by $T(G,\Phi,f)$. On the other hand, $T(G,\Phi,f)\subset
 A^{\vee}$ for every $A$.
 Therefore,  $A^{\vee}=T(G,\Phi,f)$ in the case under consideration.

 $\bf 4.1.$ Proposition.
Suppose that $A=\Hom (W,H)$, where $H$ is a subalgebra of $G$. Then
$$A^{\vee}=T(H,\Phi,f_H).$$

\begin{proof}
 Let $u\in A^{\vee}$, $A\subset f*u$. Then
$$f_H *u=(f*u)\cap\Hom (W,H)=\Hom (W,H).$$
 Therefore, $u\in T(H,\Phi,f_H)$ and, consequently, $A^{\vee}\subset
 T(H,\Phi,f_H)$. On the other hand, if $u\in T(H,\Phi,f)$, then
 $f_H*u=\Hom(W,H)\subset f*u$, $u\in A^{\vee}$. This gives us the converse
 inclusion.
\end{proof}

 $\bf 4.2.$ Proposition.
 Under the same conditions, in positive formula logic
$$T(H,\Phi,f_H)\subset A^{\vee}.$$

\begin{proof}
 Let $u$ be a positive formula and $f_H*u=\Hom (W,H)$. Then $A=\Hom (W,H)\subset
 f*u$, $u\in A^{\vee}$.
\end{proof}

Now we pass to the main results.

 If $T$ is a collection of open formulas, then $\Theta_T$ is the class of all
 those models $(G,\Phi,f)$ which all formulas from $T$ are valid in. Here,
 $G\in\Theta$, the collection $\Phi$ is fixed for the Halmos algebra $U$,
 $\Theta_T$ is an axiomatic class.

 We again choose a model $(G,\Phi,f)$. In the following theorem, which is an
 analogue of Theorem 3.1 from \S 3, the passages $T\mapsto A=T^{\vee}$ and
 $A\mapsto  A^{\vee}=T$ (only open formulas from $U$ are taken into account)
 are considered  with respect to this model.

 $\bf 4.3.$ Theorem.
 Suppose that the set $T$ is invariant relatively to the action of
the semigroup
 $\End W$. Then
    \begin{enumerate}
\item[\rm 1.] An element $\mu\in\HOM$ belongs to the variety $A=T^{\vee}$ if and only if
 the model $(H,\Phi,f_H)$ with $H=\im\mu$ belongs to $\Theta_T$,
\item[\rm 2.] the closure $T^{\vee\vee}$ coincides with the open theory of all
$\Theta_T$-submodels of the model $(G,\Phi,f)$.
    \end{enumerate}

\begin{proof}
 Assume that $\mu\in A$ and $H= Im \mu$ and consider
the model $(H,\Phi,f_H)$. We
 have to demonstrate that, for every $u\in T$, $f_H*u=\Hom (W,H)$. Let
 $\nu\in\Hom (W,H)$. Then for some $s\in\End W$ the diagram

\begin{center}
\unitlength=0.8mm
\linethickness{0.4pt}
\begin{picture}(50,37)(20,50)
\put(29.67,79.67){\vector(1,0){30.00}}
\put(21.33,80.00){\makebox(0,0)[cc]{$W$}}
\put(69.0,79.67){\makebox(0,0)[cc]{$W$}}
\put(25.00,75.00){\vector(3,-4){18.00}}
\put(65.00,75.00){\vector(-3,-4){18.}}
\put(44.67,46.50){\makebox(0,0)[cc]{$H$}}
\put(44.67,82.67){\makebox(0,0)[cc]{$s$}}
\put(29.33,63.67){\makebox(0,0)[cc]{$\nu$}}
\put(61.00,63.67){\makebox(0,0)[cc]{$\mu$}}
\end{picture}
\end{center}

\noindent
commutes, and $\nu =\mu s$. The statement $\nu = \mu s \in f_H*u$ is equivalent
 to $\mu\in f_H*su$.

 Since the formula $su$ is open, the later statement is equivalent to $\mu\in
 f*su$. It follows from the conditions that $\mu\in f*su$ holds, and then so do
 the statements
 $\mu\in f_H*su$ and $\nu =\mu s\in f_H*u$. This is the case for any $\nu\in\Hom
 (W,H)$, and $f_H*u=\Hom (W,H)$. Here $u$ is an element of $T$, so the model
 $(H,\Phi,f_H)$ belongs to $\Theta_T$.

 Now let us check the converse. Assume that, for $\im\mu =H$, the model
 $(H,\Phi,f_H)$ belongs to $\Theta_T$. We have to derive that $\mu\in
 A=T^{\vee}$.

 For all $u\in T$, $f_H*u=\Hom (W,H)$. As $\mu\in\Hom (W,H)$, we conclude that
 $\mu\in f_H*u$. Consequently, $\mu\in f*u$. This argument remains valid for
 every $u\in T$ and $\mu\in T^{\vee}=A$.

 Note that this converse condition does not depend on the assumption that $T$
 is invariant under the action of the semigroup $\End W$. We now observe
 that
        $$A=\bigcup\Hom (W,H),$$
 where the union is taken over all subalgebras $H\subset G$ such that the model
 $(H,\Phi,f_H)$ belongs to the $\Theta_T$.

 Indeed, if $\mu\in A$, then $H=\im\mu$ determines a model in $\Theta_T$, and
 $\mu\in\Hom (W,H)$. On the other hand, if $\mu\in\Hom (W,H)$ with
 $(H,\Phi,f_H)\in\Theta_T$, then $\mu\in f_H*u$, $u\in T$ and, furthermore,
 $\mu\in  f*u$. This is so for every $u$, and $\mu\in A$.

 The corresponding open theory will be denoted by
 $\widetilde T$.
 We have to
 show that $T^{\vee\vee}=\widetilde T$.

 Let $u\in\widetilde T$. This means that, for every model $(H,\Phi,f_H)$ from
 $\Theta_T$ with $H\subset G$, $f_H*u=\Hom (W,H)$. Now take $\mu\in A$. Then
 $\mu\in\Hom (W,H)=f_H*u$ for an appropriate $H$. Moreover, $\mu\in f*u$.
 This takes
 place for every $\mu\in A$, consequently $A\subset f*u$, $u\in
 T^{\vee\vee}=A^{\vee}$.

 Let, conversely, $u\in T^{\vee\vee}$. Take $H\subset G$ so that
 $(H,\Phi,f_H)\in\Theta_T$; we have to verify that $f_H*u=\Hom (W,H)$. Whenever
 $\mu\in\Hom (W,H)$, $\mu\in A$ and $\mu\in f*u$; hence, $\mu\in f_H*u$. This
 gives the inclusion $\Hom (W,H) \subset f_H*u$; the converse follows
 immediately from definitions.  Therefore, $f_H*u=\Hom (W,H)$, indeed.
 All this remains valid for every appropriate $H$  and $u\in\widetilde T$.

 The proof is completed.
\end{proof}

 It is readily seen that $f*u=\HOM$ always implies $f*su=s(f*u)=\HOM$.
 Therefore, we can maintain that $su\in\widetilde T$ if
 $u\in\widetilde T$, i.e. the class
 $\widetilde T=T^{\vee\vee}$, as well as $T$, is invariant relatively to action
 of $\End W$.

Just as for algebras, we now shall define the notion of equivalence for models.
The question is about the geometric equivalence, and it is now considered in the
class of open formulas. We say that the models $(G_1, \Phi, f_1)$ and
$(G_2, \Phi, f_2)$ are equivalent if the respective closures $T^{\vee\vee}$ of
every collection $T \subset U$ of open formulas for the models are equal.

One easily realizes that isomorphic models are equivalent. This seems obvious,
 but we shall advance a formal proof using some useful considerations.

 The model isomorphism
$$\delta\colon (G_1,\Phi,f_1)\to (G_2,\Phi,f_2)$$
 is, first of all, an isomorphism $\delta\colon G_1\to G_2$ between algebras
 from  $\Theta$. It induces a Halmos algebra isomorphism
 $\delta_*\colon\MM_{G_2}\to\MM_{G_1}$ and a bijection
 $\delta_*\colon F_{G_2}\to F_{G_1}$ betveen realization systems,
 each realization $f$ being a database state, and in this case
$$
(f*u)^{\delta_*}=f^{\delta_*}*u,\quad f\in F_{G_2},\, u\in U.
$$
 $\delta$ is a model isomorphism if and only if $f^{\delta_*}_2 = f_1$.

 Now suppose that $T\subset U$ is a system of formulas, $A=T^{\vee}$ in the
 first model and $B=T^{\vee}$ in the second one . What are relations between
 $A$ and $B$?
 We have $A=\bigcap_{u\in T}f_1*u$, $B=\bigcap_{u\in T}f_2*u$, and
 $f_1*u=f_2*u=(f_2*u)^{\delta_*}$. Hence $A=B^{\delta_*}$, and
 $\mu\colon W\to G_1$ belongs to $A$ if and only if
 $\mu\delta = \nu\colon W\to G_2$ belongs to  $B$.

 Now
$$T_1^{\vee\vee}=A^{\vee}=\{ u\in U,\ A\subset f_1*u\},$$
$$T_2^{\vee\vee}=B^{\vee}=\{ u\in U,\ B\subset f_2*u\}.$$
 If $u\in T^{\vee\vee}_1$, then $A\subset f_1*u$, $B^{\delta_*}\subset
 f_2^{\delta_*}*u=(f_2*u)^{\delta_*}$ and $B\subset f_2*u$, $u\in
 T^{\vee\vee}_2$.

 The converse is proved analogously, and $T_1^{\vee\vee}=T_2^{\vee\vee}$. If we
 have isomorphic models, this holds for every $T$.

 We can confine ourselves here to systems of open formulas.

 Just as before, the question is how to learn to recognize equivalence of two
 models. The next theorem is similar to Theorem 4.3 of \S 3.

$\bf 4.4.$ Theorem.
 Suppose that models $(G_1,\Phi,f_1)$ and $(G_2,\Phi,f_2)$ are equivalent in the
 open formula geometry. Then their open theories coincide.

\begin{proof}
 First of all, we shall prove the
 well-known fact that if ${\mathcal X}$ is the class of
 models defined by a collection $T$ of open formulas, then, for every model
 $(G,\Phi,f)\in{\mathcal X}$, all submodels of it also are in ${\mathcal X}$.

 Assume that $(H,\Phi,f_H)$ is a submodel. We have to verify that $f_H*u=\Hom
 (W,H)$ for every formula $u\in T$. Since $u$ is open, $f_H*u=(f*u)\cap\Hom
 (W,H)$. But $f*u=\HOM$; so, $f_H*u=\Hom (W,H)$.

 Now assume that $T_1$ is an open theory of the first model and
 ${\mathcal X}_1$ is the
 class of models defined by formulas from $T_1$. Similarly, we take $T_2$ and
 ${\mathcal X}_2$ for the second model. $T_1^{\vee}$ for
 the former model is $\HOM$, and
 $T_1^{\vee\vee}$ for it is $T_1$. The same holds for the latter one; we count
 the models as equivalent.

 Further, let $(H_{\alpha},\Phi,f_{\alpha})$ be the family of all those
 submodels of $(G_2,\Phi,f_2)$ which all formulas from $T_1$ are valid in. Then
 $T_1$ is the open theory of the class of these models. All models
 $(H_{\alpha},\Phi,f_{\alpha})$ belong to ${\mathcal X}_2$. As they also belong to
 ${\mathcal X}_1$, and even generate ${\mathcal X}_1$, we obtain that
 ${\mathcal X}_1\subset{\mathcal X}_2$. Just
 in the same way we conclude that ${\mathcal X}_2\subset{\mathcal X}_1$.
 Hence the theorem.

 The theorem can also be applied in the case the set $\Phi$ is empty and the
 question is about algebras. In this case open formulas are treated as
 universal, and we then speak of universal theories and universal classes of
 algebras.

 Therefore, if two algebras are geometrically equivalent in universal logic,
 then they have the same universal theory.
\end{proof}

\paragraph{5.  Halmos algebras and Boolean algebras of varieties. Conclusion.}
\stepcounter{paragraph}
 Assume that $H$ is a Halmos algebra and $T$ is its Boolean filter. Let $T^-$
 stand for the subset of $T$ consisting of all $h$ from $T$ such that $\forall
 (X)h\in T$. Here, $X$ is the set of variables. It is well-known that $T^-$
 is a filter of the Halmos algebra $H$. It also is known that if $T$ is an
 ultrafilter, then $T^-$ is a maximal filter of $H$, and the algebra $H/T^-$ is
 simple.

 We also note that if $T_{\alpha}$, $\alpha\in I$, is a collection of Boolean
 filters of $H$, then
$$(\bigcap\limits_{\alpha}T_{\alpha})^-
=\bigcap\limits_{\alpha}(T^-_{\alpha}).$$
 Now let $H$ be the Halmos algebra $U$. For some algebra $G\in\Theta$ and set
 $A\subset\HOM$, we considered in \S 3 the algebra $W/A'$ which can be presented
 as a subdirect product of all $W/\Ker\mu$ with $\mu\in A$.

 Assume now that a model $(G,\Phi,f)$ is considered. We shall relate a Halmos
 algebra to the same set $A$. Take the Boolean filter $T=A^{\vee}$ of the Halmos
 algebra $U$ and pass to $T^-=A^{\nabla}$. The corresponding Halmos algebra is
 $U/A^{\nabla}$. If, furthermore, $A=\{\mu\}$ is a singleton, then set
 $A^{\nabla}=T_{\mu}$. As $A^{\vee}$ is an ultrafilter, $T_{\mu}$ is a maximal
 filter in $U$, $T_{\mu}\in\Spec U$.

 It is easily seen that $A^{\nabla}=\bigcap_{\mu\in A}T_{\mu}$ for $A$
 arbitrary. Indeed, we always have $A^{\vee}=\bigcap\limits_{\mu\in
 A}\{\mu\}^{\vee}$. If we apply the operation ${}^-$, we obtain
$$A^{\nabla}=(\bigcap_{\mu\in A}\{\mu\}^{\vee})^- =\bigcap\limits_{\mu\in
A}\{\mu\}^{\nabla}=\bigcap\limits_{\mu\in A}T_{\mu}.$$

 Thus we come to

 $\bf 5.1.$ Theorem.
 For every $A$, the Halmos algebra $U/A^{\nabla}$ is a subdirect product of all
 the simple algebras $U/T_{\mu}$ with $\mu\in A$.

One could take for $A$ the algebraic variety defined by a collection of
 formulas $T\subset U$: $A=T^{\vee}$. Then $A^{\vee}=T^{\vee\vee}$ and
 $A^{\nabla}=T^{\vee\nabla}$.

 We do not know, however, what role the algebra $U/A^{\nabla}$ plays for the
 variety $A$. This is one among a lot of questions the answer to which is to be
 get known.

 It is not improbable that the relation algebras, a product of the categorial
 approach to algebraic logic ( see, e.g. \cite{BPl3}), naturally shall find their
 applications in this theory along with Halmos algebras.

 We note, furthermore, that to any variety $A$ a Boolean algebra $U/A^{\vee}$
 can be related, and bring forward the following proposition in this connection.

 $\bf 5.2.$ Proposition.
 Every variety morphism $s\colon A\to B$ induces a Boolean homomorphism
$$\bar s\colon U/B^{\vee}\to U/A^{\vee}.$$

\begin{proof}
 We first make some preliminary notes. If $A$ is a subset of $\HOM$, then $As$
 stands for the set of all $\mu s$, $\mu\in A$. An endomorphism $s\in\End W$
 determines a morphism $A\to B$ if and only if $As\subset B$.

 Given a collection $T$ of elements of $U$, we define the set $Ts$ by the rule
 $u\in Ts\Leftrightarrow su\in T$. Then the following equality holds:
$$(As)^{\vee}=A^{\vee}s.$$
 Let us verify it.

 Suppose that $\in (As)^{\vee}$. Then $As\subset f*u$. Furthermore, $\mu s\in
 f*u$ whenever $\mu\in A$. Hence, $u\in s(f*u)=f*su$, $A\subset f*su$, $su\in
 A^{\vee}$, $u\in A^{\vee}s$.

 Now let $u\in A^{\vee}s$, $su\in A^{\vee}$, $A\subset f*su$. Then for every
 $\mu\in A$, $\mu\in s(f*u)$ and $\mu s\in f*u$. Therefore, $As\subset f*u$ and
 $u\in As^{\vee}$.

 Now assume that we are given a morphism $s\colon A\to B$. Then $As\subset B$
 and $(As)^{\vee}=A^{\vee}s\supset B^{\vee}$. This means that $su\in A^{\vee}$ if
 $u\in B^{\vee}$.

 Let us consider the composition Boolean homomorphism
$$U\,{\stackrel{s}{\to}}\,U\,{\stackrel{\sigma_0}{\to}}\,U/A^{\vee}.$$
 Since $u\in B^{\vee}$ implies $su\in A^{\vee}$, the filter $B^{\vee}$ is
 included into the kernel of the homomorphism $s\sigma_0$. This means that
 $s\sigma_0$ induces the homomorphism $\bar s$ we look for.
 \end{proof}

 So the proposition is proved. It admits conversion. Suppose that $B^{\vee}$ is
 included in the kernel of $s\sigma_0$. This means that, for every $u\in
 B^{\vee}$, the element $su$ belongs to $A^{\vee}$, i.e. $sB^{\vee}\subset
 A^{\vee}$. Applying ${}^{\vee}$, we get $A^{\vee\vee}=A\subset
 (sB^{\vee})^{\vee}=sB^{\vee\vee}=sB$. What does this mean is that $s$
 specifies a
 morphism $A\to B$. If, in particular, the varieties $A$ and
 $B$ are isomorphic, then so are the algebras $U/B^{\vee}$ and $U/A^{\vee}$.

 Let us verify this. Suppose that $\mu ss'=\mu$ for every $\mu\in A$. We have to
 check that, for every $u\in U$, the elements $ss'u$ and $u$ are equivalent
 modulo the filter $A^{\vee}$, i.e.
$$(ss'u\to u)\land (u\to ss'u)\in A^{\vee}$$
 or, what amounts to the same,
$$A\subset f*(ss'u\to u)\cap f*(u\to ss'u).$$
 Take $\mu\in A$ and let $\mu\in f*ss'u$. Then $\mu ss'=\mu\in f*u$. Therefore,
 $A\subset f*(ss'u\to u)$. Analogously, if $\mu=\mu ss'\in f*u$, then $\mu\in
 ss'(f*u)=f*ss'u$. So $A\subset f*(u\to ss'u)$. Thus, $ss'u$ and $u$ are
 equivalent.

 In the same fashion we check that $\nu=\nu ss'$ for $\nu\in B$ implies
 equivalence of $ss'u$ and $u$ modulo $B^{\vee}$. All this eventually means that
 the respective homomorphism $\bar s\colon U/B^{\vee}\to U/A^{\vee}$ and
 $\bar {s'}\colon U/A^{\vee}\to U/B^{\vee}$ are mutually inverse.

 The converse is not generally true, for the endomorphisms of the
 Boolean algebra
 $U$ by no means are exhausted by elements of the semigroup $\End W$.
 But we can, of course, confine ourselves to the isomorphisms between
 $U/B^{\vee}$ and $U/A^{\vee}$ induced by endomorphisms from $\End W$.

 Our further remarks are also related with the notion of an algebraic variety
 morphism. The question is of varieties specified by subsets of the Halmos
 algebra $U$.

 If $A$ and $B$ are varieties in $\HOM$ with a model $(G,\Phi,f)$, then the
 morphism $A\to B$ is determined by some $s\in\End W$ such that $\nu s\in B$ for
 every $\nu\in A$. If, for example, $A=sB$, then we have $s\colon A\to B$.

 Let us consider separately the case when the set $T$ is invariant under $s$,
 i.e. $sT\subset T$. Assume that $A=T^{\vee}$. Then
 $(sT)^{\vee}=sT^{\vee}=sA\supset T^{\vee}=A$. This means that the variety $A$
 is invariant under $s$, $\nu\in A$ implies $\nu s\in A$, and we have a morphism
 $s\colon A\to A$.

 If $T$ is invariant under $\End W$, then so is $A=T^{\vee}$ as well, and the
 semigroup $\End W$ acts on $A$ as an semigroup of endomorphisms of the variety.

 If, on the other hand, the variety $A=T^{\vee}$ is invariant relatively to the
 action of $\End W$, then always $sA\supset A=T^{\vee}$ and
 $T^{\vee\vee}\supset (sA)^{\vee}$. However, we cannot present $(sA^{\vee})$ as
 $sA^{\vee}=sT^{\vee\vee}$, so we cannot claim that the collection
 $T^{\vee\vee}$ is invariant relatively to the action of $\End W$.

 Still if the variety $A=T^{\vee}$ is invariant relatively to $\Aut W$, then the
 filter $T^{\vee\vee}$ also has the same feature.

 One can consider, for any variety $A$, the semigroup $\End A$ of those
 $s\in\End W$ acting on $A$. It is a subsemigroup of $\End W$ and characterizes
 the given $A$. $\Aut A$ is the symmetry group of the variety $A$. If $A$ and
 $B$ are isomorphic, then $\Aut A$ and $\Aut B$ are conjugated in a way. For
 $\Phi$ empty, and on the equational level, the groups $\Aut A$ and $\Aut
 (W/A')$ are well-connected; they are canonically anti-isomorphic.

 As to the group $\Aut A$, we have to note in addition that, in fact, it is not
 a subgroup of $\Aut W$. Not every automorphism of the variety $A$ is induced by
 some $s\in\Aut W$. The subset of those $s\in\Aut W$ under which $A$ is
 invariant is a subgroup of $\Aut W$ and, at the same time, of $\Aut A$. We
 denote the subsemigroup by ${\rm Aut_0}\, A$. The initial group for all the
 groups ${\rm  Aut_0}\, A$ was $\Aut W$; it acts on $\HOM$ and is
 well-coordinated with several topologies that have been studied.
 Other groups related to topologies on $\HOM$ could be taken
 instead of $\Aut W$.
 They might be either the topologies discussed above or any other ones. Along
 these lines, the geometry of varieties can be enriched.

 When $\Phi$ is empty, every algebraic variety $A$ is also invariant relatively
 to the action of $\Aut G$. Actions of $\Aut G$ and $\Aut A$ on the variety
 $A$ determine the geometry of it in many respects.

 We further make some conclusive
 remarks. We have discussed the notion of algebraic closeness of an
 algebra $G\in\Theta$  in the subsection devoted to the
 Hilbert theorem on zeros. The algebra $G$ is said to be
 algebraically closed if,
 for every finite $X$ and proper congruence $T$ of $W(X)$,
 the variety $T'=A$ in
 $\HOM$ is nonempty (equivalently, for every $T$ the set ${\rm Hom}\, (W/T,G)$
 is nonempty).

 If $G$ has a one-element subalgebra, then we must regard $T'$ be nontrivial
 rather than nonempty. It seems likely that such a version of the definition of
 algebraic closeness for algebras deserves a special discussion.

 Also, it is natural to speak of algebraic closeness for some pregiven $X$.

 We now can claim that if $G_1$ and $G_2$ are algebras without one-element
 subalgebras, and if one of them is algebraically closed while the other is
 not, then the algebras are not equivalent. Indeed, suppose that $G_1$ is
 algebraically closed and $G_2$ is not. Then there is a proper
 congruence $T$ on
 $W(X)$ for which the variety $B=T_{G_2}'$ is empty. The variety $A=T'_{G_1}$
 for the same $T$ is not empty. It follows therefrom, as we saw in \S 3, that
 $G_1$ and $G_2$ are not equivalent.

 So, we generally cannot claim that algebraic closeness of algebras $G_1$ and
 $G_2$ implies their equivalence.

 The claim is justified in the case of classical geometry, but not always. The
 algebras $G_1$ and $G_2$ can both be algebraically closed and nevertheless
 have distinct systems of identities. Then they are not equivalent.

 Of some interest is the situation when $G_1$ and $G_2$ are algebraically
 closed and do have the same identities.

 In the context in question, if an algebra $G$ is algebraically closed,
 then every algebra, containing $G$, has the same property. In particular,
 we consider extension $K$ of the field $P$ as algebraically closed if
 $K$ contains algebraic closure of the field $P$.

 Now, we return to geometric equivalence of algebras. We shall assume
 that there
 is a nullary operation $0$ among the ground operations of the variety $\Theta$
 which singles out of every $G\in\Theta$ a one-element subalgebra.

 Let us begin with few general facts and then apply them to commutative groups.

 $\bf 5.3.$ Proposition.
 Suppose that $I$ is a set and every $\alpha\in I$ is assigned an algebra
 $H_{\alpha}=\Theta$. Then, for every $G\in\Theta$,
$$
(\prod_{\alpha}H_{\alpha}\defis\Ker)(G) =
    \bigcap_{\alpha}(H_{\alpha}\defis\Ker)(G).
$$

\begin{proof}
 $\prod_{\alpha}$ is the Cartesian multiplication. Let $\tau_1$ stands
 for the left hand congruence and $\tau_2$--for the right hand one. We shall use
 standard arguments which already were in use above. Assume that $g_1$ and $g_2$
 are two elements of $G$, and that $g_1\tau_2 g_2$. This means that
 $g_1^{\nu}=g_2^{\nu}$ for all $\alpha\in I$ and $\nu$: $G\to H_{\alpha}$.
 Let us take
 $\mu\colon G\to\prod_{\alpha}H_{\alpha}$ and verify that
 $g_1^{\mu}=g_2^{\mu}$.

 The  equality means that $g_1^{\mu}(\alpha)=g_2^{\mu}(\alpha)$ for every
 $\alpha\in I$. We use the projections
 $\pi_{\alpha}\colon\prod_{\alpha}H_{\alpha}\to H_{\alpha}$
 and denote $\mu\pi_{\alpha}=\nu_{\alpha}$. Then
 $
 g_1^{\mu}(\alpha)=g_1^{\nu_{\alpha}} = g_2^{\nu_{\alpha}}=g_2^{\mu}(\alpha)
 $,
 i.e.  $g_1^{\mu}=g_2^{\mu}$ and, further, $g_1\tau_1g_2$.

 Assume that, conversely, $g_1\tau_1g_2$. Given $\alpha\in I$ and $\nu\colon
 G\to H_{\alpha}$,  we define $\mu$ by the rule: $g^{\mu}(\alpha)=g^{\nu}$, and
 $g^{\mu}(\beta)$ is the zero if $\beta\ne\alpha$. Then $\mu\colon G\to
 \prod_{\alpha}H_{\alpha}$ and $g_1^{\mu}=g_2^{\mu}$. But then
 $g_1^{\mu}(\alpha)=g_1^{\nu}=g_2^{\mu}(\alpha)=g_2^{\nu}$. Therefore,
 $g_1\tau_2g_2$.
\end{proof}

 Apart from Cartesian products $\prod_{\alpha}$, we shall also discuss
 direct products $\prod_{\alpha}^0$. Here $\prod_{\alpha}^0 H_{\alpha}$
 consists
 of all those $\alpha\in\prod_{\alpha}H_{\alpha}$ that take the zero as
 the value almost everywhere. This is a subalgebra of
 $\prod_{\alpha}H_{\alpha}$.

 Clearly, if $H_1$ is a subalgebra of $H_2$, then
$$(H_1\defis\Ker)(G)\supset (H_2\defis\Ker)(G).$$

 $\bf 5.4.$ Proposition.
Always
$$({\prod\limits_{\alpha}}^0
H_{\alpha}\defis\Ker)(G)=(\prod\limits_{\alpha}H_{\alpha}\defis\Ker)(G).$$

\begin{proof}
  Of course,
$$({\prod\limits_{\alpha}}^0 H_{\alpha}\defis\Ker)(G)\supset (\prod\limits_{\alpha}
H_{\alpha}\defis\Ker)(G)=\bigcap\limits_{\alpha}((H_{\alpha}\defis\Ker)(G)).$$
 On the other hand,
$$(H_{\alpha}\defis\Ker)(G)\supset ({\prod\limits_{\alpha}}^0 H_{\alpha}\defis\Ker
)(G),$$
 for $H_{\alpha}$ is a subalgebra of the direct product. This is so for any
 $\alpha\in I$, so that
$$\prod\limits_{\alpha}((H_{\alpha}\defis\Ker)(G))\supset
({\prod\limits_{\alpha}}^0
H_{\alpha}\defis\Ker)(G).$$
 These two inclusions justify the needed equality.
\end{proof}

 $\bf 5.5.$ Proposition.
 Suppose that algebras $H_{\alpha}$ and $H'_{\alpha}$, $\alpha\in I$, are
 geometrically equivalent. Then $\prod_{\alpha} H_{\alpha}$ and
 $\prod_{\alpha}H'_{\alpha}$ are also equivalent.

\begin{proof}
 For an arbitrary $G = W/T$,
$$
(\prod_{\alpha}H_{\alpha}\defis\Ker)(G) =
    \bigcap_{\alpha}(H_{\alpha}\defis\Ker)(G) =
        \bigcap_{\alpha}(H'_{\alpha}\defis\Ker(G)) =
            (\prod_{\alpha}H'_{\alpha}\defis\Ker)(G).
$$
\end{proof}

This gives what is desired.

 $\bf 5.6.$ Proposition.
Suppose that every algebra $H_\alpha$ in $\prod_\alpha H_\alpha$ can be regarded
as a subalgebra of one of them, call it $H$. Then $H$ and $\prod_\alpha
H_\alpha$ are geometrically equivalent.

\begin{proof}
$(\prod_\alpha H_\alpha\defis \Ker)(G) = \bigcap_\alpha(H_\alpha\defis\Ker)(G) =
(H\defis\Ker)(G)$,
for always\\
$(H_\alpha\defis\Ker)(G) \supset (H\defis\Ker)(G)$.
\end{proof}

 Now we pass to commutative groups.

 $\bf 5.7.$ Theorem.
 Two commutative groups, $A_1$ and $A_2$, with finite exponents are
 geometrically equivalent if and only if their exponents coincide or, what amounts to the
 same, $\Var A_1=\Var A_2$.

\begin{proof}
 In one direction the assertion follows from a general result in \S 3. We shall
 prove the converse.

 Assume first that the groups $A_1$ and $A_2$ are primary with respect
 to the same
 $p$, and that $p^n$ is their common exponent. Both groups are directs
 products of cyclic ones. All these cyclic groups in $A_1$ are subgroups of some
 group $A_1^0$ of order $p^n$. In the same manner we isolate a group $A_2^0$ of
 the same order. Then
$$
(A_1\defis\Ker)(G)=(A_1^0\defis\Ker)(G), \quad
        (A_2\defis\Ker)(G)=(A_2^0\defis\Ker)(G).
$$
 Since $A_1^0$ and $A_2^0$ are isomorphic,
$$(A_1\defis\Ker)(G)=(A_2\defis\Ker)(G),$$
 so $A_1$ and $A_2$ are equivalent.

 We pass to the general case. Let $A_1={\prod_p}^0 A_{1,p}$ and
 $A_2={\prod_p}^0 A_{2,p}$ are decompositions into Sylow groups. If $m$ is
 the common exponent of $A_1$ and $A_2$, $m=p_1^{n_1}\cdots p_k^{n_k}$, then in
 both cases $p$ runs over the collection $p_1,\cdots,p_k$, and the groups
 $A_{1,p}$ and $A_{2,p}$ are of the same exponent. Hence, they are equivalent,
 and so are the groups $A_1$ and $A_2$.
\end{proof}

 It is easy to prove also the following.

\smallskip
\noindent {\bf 1.} Any two commutative torsion-free groups are $X$-equivalent for every
 finite $X$.

\smallskip
\noindent {\bf 2.} Two mixed finitely generated commutative groups are equivalent if and only
 if their periodic parts are of the same exponent.

\smallskip
We now turn attention to some problem that concerns the classical situation.

 Consider a field $P$ and two its extensions $K_1$ and $K_2$. If both $K_1$ and
 $K_2$ are algebraically closed, then they are geometrically equivalent. Then
 they have the same equational theories. Actually, it is known that even their
 elementary theories coincide.

 The question is whether or not $K_1$ and $K_2$ are geometrically equivalent on
 the universal logic level if they are algebraically closed.

 Let us remind that we have proved, on the equational level, equivalence of every
 algebra to every of its Cartesian powers. In parallel, the following problem
 could
 be noticed: Whether or not every model (algebra) is equivalent, already on the
 level of open logic, to every of its ultrapowers. We now shall discuss a
 partial solution of this problem, and afterwards shall point to some
 applications to fields.

 We begin with some preliminary notes.

 Assume we are given a model $(G,\Phi,f)$, $G\in\Theta$, a set $I$ and an
 ultrafilter $D$ on it. Let $(\bar G,\Phi,\bar f)$ be the ultrapower of the
 model with respect to $D$. Here, $\bar G$ is the ultrpower of the algebra $G$,
 $\bar G=G^I/D$, and $\bar f$ is defined in a special way (see, for example,
 \cite{BPl3}). We take, furthermore, $\HOM$ and $\Hom (W,\bar G)$. For every formula
 $u\in U$, $f*u\subset\HOM$ and $\bar f*u\subset\Hom (W,\bar G)$. If $\mu\colon
 W\to\bar G$ and if $\xi\colon G^I\to\bar G$ is the natural homomorphism, then
 we obtain the commutative diagram

\begin{center}
\unitlength=1mm
\linethickness{0.4pt}
\begin{picture}(50,32)(25,6)
\put(30.00,35.00){\makebox(0,0)[cc]{$W$}}
\put(33.00,35.00){\vector(1,0){20.00}}
\put(56.33,35.60){\makebox(0,0)[cc]{$\bar G$}}
\put(33.00,32.00){\vector(1,-1){20.00}}
\put(56.00,9.00){\makebox(0,0)[cc]{$G^I$}}
\put(56.00,12.00){\vector(0,1){20.67}}
\put(43.00,38.00){\makebox(0,0)[cc]{$\mu$}}
\put(58.33,20.33){\makebox(0,0)[cc]{$\xi$}}
\put(39.5,20.33){\makebox(0,0)[cc]{$\nu$}}
\end{picture}
\end{center}

\noindent
 Here, in general, $\nu$ is not uniquely determined by $\mu$. We also consider
 the projections $\pi_{\alpha}\colon G^I\to G$, and the compositions
 $\nu_{\alpha}=\nu\pi_{\alpha}\colon W\to G$.

 The following rule holds (cf. \cite{BPl3}):
$$
\mu\in\bar f*u\Leftrightarrow\{\alpha\in I,\ \nu_{\alpha}\in f*u\}\in D.
$$
 It does not depend on the choise of $\nu$.

 $\bf 5.8.$ Theorem.
Let $T$ be a finite subset of $U$. Then $T^{\vee\vee}$ for the model
 $(G,\Phi,f)$ coincides with $T^{\vee\vee}$ for the model $(\bar G,\Phi,\bar f)$.

\begin{proof}
 We take $T^{\vee}=B$ for the initial model and $T^{\vee}=A$ for the
 ultrapower,
 and set $T^{\vee\vee}=B^{\vee}$ and $T^{\vee\vee}=A^{\vee}$. We then have to
 prove that $A^{\vee}=B^{\vee}$.

 First let $u\in B^{\vee}$, $B\subset f*u$. We shall verify that $u\in
 A^{\vee}$, $A\subset \bar f *u$, i.e. that $\mu\in\bar f*u$
 whenever $\mu\in A$. The
 conclusion $\mu\in\bar f*u$ means that,
$$\{\alpha\in I,\ \nu_{\alpha}\in f*u\}\in D.$$
 The condition $\mu\in A$ means that $\mu\in\bar f*v$ for every $v\in T$,
i.e.
$$
\{\beta\in I,\ \nu_{\beta}\in f*v\} =I_{\mu ,v}\in D.
$$

 Now set $I_{\mu}=\bigcap_{v\in T}I_{\mu ,v}$; then $I_{\mu}\in D$. By
 the definition,
$$
\alpha\in I_{\mu}\Rightarrow\nu_{\alpha}\in B=\bigcap_{v\in T}f*v.
$$
 Moreover,
$$
\{\alpha\in I,\ \nu_{\alpha}\in B\}\in D.
$$
 By the choise of $B$, $B\subset f*u$, consequently, $\nu_{\alpha}\in B$ implies
 $\nu_{\alpha}\in f*u$. Also
$$
\{\alpha\in I,\ \nu_{\alpha}\in f*u\}\in D.
$$
 But then $\mu\in\bar f*u$. Hence, $B^{\vee} \subset A^{\vee}$.

 To prove the converse inclusion, finiteness of $T$ is not needed. Let $u\in A^{\vee}$. This means that $A\subset\bar f*u$. We
 are going to demonstrate that $B\subset f*u$. Let $\nu_0\in B$, $\nu_0\colon
 W\to G$. We also take a constant $\nu\colon W\to G^I$ so, that
 $\nu_{\alpha}=\nu\pi_{\alpha}=\nu_0$ for every $\alpha$. Finally, let $\mu
 =\nu\xi$.

 We shall need to know that $\mu\in A$, i.e. that, for every $v\in T$,
 $\mu\in\bar f*u$ or, equivalentially,
$$
\{\alpha\in I,\ \nu_{\alpha}\in f*v\}\in D.
$$
 As $\nu_{\alpha}=\nu_0$ and $\nu_0\in B\subset f*v$, we conclude that
$$
\{\alpha\in I,\ v_{\alpha}\in f*v\}=I\in D.
$$
 Therefore, $\mu\in A$ and, furthermore, $\mu\in\bar f*u$, i.e.
$$
\{\alpha\in I,\ \nu_{\alpha}=\nu_0\in f*u\}\in D.
$$
 This set is not empty, and $\nu_0\in f*u$. This is so for every $\nu_0\in B$,
 and then $B\subset f*u$, $u\in B^{\vee}$.

 The theorem is proved.
\end{proof}

It can be applied in the classical situation. If $K$ is a field that is an
 extension of a field $P$, then each of its ultrapowers $\bar K$
 is also a field
 extending $P$. Therefore, $T_K^{\vee\vee}=T_{\bar K}^{\vee\vee}$ for every
 finite collection $T$.

 Moreover, we now may say that $K$ and $\bar K$ are geometrically equivalent on
 the equational level and for any finite $Y \subset X$.
 Indeed, if a finite $Y$ is selected,
 then it suffices to confine ourselves to finite collections $T$ of polynomials.
 Now $T^{\vee\vee}_K=T_{\bar K}^{\vee\vee}$ implies $T''_K=T''_{\bar K}$.

Here, we consider the collection  $T$ as a set of formulas in the
Halmos algebra $U$, and write $T=T_X$. In Section 4, the equality
$T''_X \cap W(Y)=T''_Y$, $Y \subset X$ was true for every algebra
$G \in \Theta$. In the situation when $T=T_Y$ and $W(Y)=P[Y]$ we have
$$
T''_{Y,K}= T''_{X,K} \cap P[Y] =
 T''_{X,\bar K} \cap P[Y] = T''_{Y, \bar K}.
$$
This means that $K$ and $\bar K$ are $Y$-equivalent for any finite $Y$.

The subsequent Proposition and Problem also relate to the classical case.

 $\bf 5.9.$ Proposition.
Two finite extensions $K_1$ and $K_2$ of the field $P$ are equivalent if
and only if they are isomorphic.

 $\bf 5.10.$ Problem.
Are every two really closed extensions of the field $P$ always equivalent?

The proof of Proposition 5.11 is similar to that of Proposition 5.9.

 $\bf 5.11.$ Proposition.
Let $\Theta$ be the variety of all associative algebras over the field $P$.
Then finite dimensional simple algebras $G_1$ and $G_2$ in $\Theta$ are
equivalent if and only if they are isomorphic.

The similar is true for simple Lie algebras.

Finally, we formulate a general problem concerning group representations.

 Suppose that $V$ is a $K$-module, $K$--a commutative ring with unit,
 and $\Aut V$--the automorphism group of $V$. Select $X$ and let $F=F(X)$ be
 the group of free over $X$. Consider the set of representations $\Hom (F,\Aut
 V)$. Also, select a variety $\Iks$ of group representations over $K$.
 Single out a subset $A$ of $\Hom (F,\Aut V)$ according to the rule:
    \begin{quote}
$ \mu\in A$ $\Leftrightarrow$ the representation of $(V,\ \im\mu)$ belongs to
 $\Iks$.
    \end{quote}
 What can be said about $A'$, $A^{\vee}$, $A''$ and $A^{\vee\vee}$?
 Here, of interest are both questions--what is in common for all $\Iks$, and
 what is the state of affairs for several concrete $\Iks$, e.g. for
 $\Iks ={\frak S}^n$ (see ~\cite{PlV}).

 In this problem, the corresponding $\mu$ is characterized by properties of the
 action. If $\Iks$ is the variety of groups, and if we consider group
 properties, then
 we know the answer: $A=T'$, where $T$ is the collection of identities of ${\mathfrak X}$. Therefore, in this case $A=\{\mu,\ \im\mu\in{\Iks}\}$ is an algebraic
 variety.

 We have considered in the paper a certain general scheme. As to really
 deep and
 interesting investigations, they have to be related to several special
 varieties $\Theta$. Along with  classical varieties $\Theta$, it is natural
 to admit the
 variety of modules over a fixed $K$. For myself--I would like to select
 varieties of interesting representations of groups over $K$.

 A lot of problems related to solving equations in groups use group
 representations by groups of tree automorphisms rather than linear
 representations.

\end{document}